\documentclass[twoside,11pt]{article}
\usepackage[top=1.20in,bottom=1.20in,left=1.00in,right=1.00in]{geometry}
\usepackage{amssymb} 
\usepackage{amsmath}  
\usepackage{amsthm}  
\usepackage{amscd} 
\usepackage{amsfonts} 
\usepackage{fancyhdr}       
\usepackage{bm}        
\usepackage{prettyref}  
\usepackage{algorithmic, algorithm, mathtools}   
\usepackage{tikz} 
\usepackage{pgfplots,pgfplotstable}  
\usepackage{epsfig,epsf,bbm}    
\usepackage{dsfont}
\usepackage[toc,page]{appendix}  
\usepackage{makecell}
\usepackage{caption} 
\usepackage{subcaption}
\usepackage{graphicx}
\usepackage{diagbox} 
\usepackage{natbib} 

\definecolor{darkred}{RGB}{100,0,0}
\definecolor{darkgreen}{RGB}{0,100,0} 
\definecolor{darkblue}{RGB}{0,0,150}  
\usepackage{hyperref}
\hypersetup{bookmarksnumbered, colorlinks=true, linkcolor=darkred, citecolor=darkgreen, urlcolor=darkblue}

\usepackage{threeparttable}
\usepackage[subpreambles=true]{standalone} 
\usepackage{import} 
\usetikzlibrary{pgfplots.groupplots} 

\theoremstyle{plain} 
\numberwithin{equation}{section}

\newtheorem{theorem}{Theorem}[section]
\newtheorem{lemma}[theorem]{Lemma}
\newtheorem{corollary}[theorem]{Corollary}
\newtheorem{proposition}[theorem]{Proposition}

\theoremstyle{definition}

\newtheorem{model}{Model}
\newtheorem{assumption}{Assumption}

\theoremstyle{remark}

\let\vec=\mathbf 
\let\bs=\boldsymbol

\newcommand{\vct}[1]{\bm{#1}}
\newcommand{\mtx}[1]{\bm{#1}}
\renewcommand{\vec}[1]{{\boldsymbol{#1}}}
 
\newcommand{\colspan}{\operatorname{colspan}}
\newcommand{\rowspan}{\operatorname{rowspan}}
\newcommand{\rank}{\operatorname{rank}}
\newcommand{\diag}{\operatorname{diag}}

\definecolor{xl}{RGB}{200,50,120}

\newcommand*\circled[1]{\tikz[baseline=(char.base)]{\node[shape=circle,draw,inner sep=2pt](char) {#1};}}

\begin{document}
\title{Nonconvex Matrix Completion with Linearly Parameterized Factors}
\author{Ji Chen\thanks{Department of Mathematics, University of California, Davis, Davis CA 95616}, ~~Xiaodong Li\thanks{Department of Statistics, University of California, Davis, Davis CA 95616}, ~~and Zongming Ma\thanks{Department of Statistics, The Wharton School, University of Pennsylvania, Philadelphia PA 19104}}

\date{\vspace{-5ex}}



\maketitle

\begin{abstract}
Techniques of matrix completion aim to impute a large portion of missing entries in a data matrix through a small portion of observed ones. In practice including collaborative filtering, prior information and special structures are usually employed in order to improve the accuracy of matrix completion. In this paper, we propose a unified nonconvex optimization framework for matrix completion with linearly parameterized factors. In particular, by introducing a condition referred to as Correlated Parametric Factorization, we can conduct a unified geometric analysis for the nonconvex objective by establishing uniform upper bounds for low-rank estimation resulting from any local minimum. Perhaps surprisingly, the condition of Correlated Parametric Factorization holds for important examples including subspace-constrained matrix completion and skew-symmetric matrix completion. The effectiveness of our unified nonconvex optimization method is also empirically illustrated by extensive numerical simulations. 
\end{abstract}


\section{Introduction}
\label{sec:intro}
Matrix completion techniques predict missing entries in a data matrix from partially observed ones.
Their most successful application includes collaborative filtering \citep{rennie2005fast,candes2009exact} in which unobserved user-item ratings are predicted with the available ones. 
To put the problem in mathematical terms: Let $\boldsymbol{M}^\star$ be an $n_1 \times n_2$ matrix whose rank is $r \ll \min(n_1, n_2)$, and we aim to estimate the whole matrix from a small proportion of noisy observed entries. To be specific, let $\Omega\subset[n_1]\times [n_2]$ be the index set that supports all observed entries. The observation is represented by 
\begin{equation}
\label{eq:model_observation}
\mathcal{P}_{\Omega}(\boldsymbol{M}) = \mathcal{P}_{\Omega}(\boldsymbol{M}^\star + \boldsymbol{N}), 
\end{equation}
where $\mtx{N}$ is a matrix consisting of noise, and the operator $\mathcal{P}_{\Omega}(\cdot)$ preserves the entries on $\Omega$ while zeros the entries on $\Omega^c$.

Note that any rank-$r$ matrix can be parameterized through the factorization $\mtx{X} \mtx{Y}^\top$, where both $\mtx{X}$ and $\mtx{Y}$ have $r$ columns. With this parameterization, the regularized least squares objective proposed in \citet{ge2016matrix, ge2017no} is
\begin{align}
\label{eq:obj_rect}
f(\boldsymbol{X},\boldsymbol{Y}) \coloneqq \frac{1}{2p}\|\mathcal{P}_{\Omega}(\boldsymbol{X}\boldsymbol{Y}^\top-\boldsymbol{M})\|_F^2 + \frac{1}{8} \|\boldsymbol{X}^\top\boldsymbol{X} - \boldsymbol{Y}^\top\boldsymbol{Y}\|_F^2 +\lambda (G_{\alpha}(\boldsymbol{X}) + G_{\alpha}(\boldsymbol{Y})),
\end{align}
where 
\begin{equation}
	\label{eq:G_alpha}
G_{\alpha}(\boldsymbol{X}) \coloneqq \sum_{i=1}^{n}[(\|\boldsymbol{X}_{i,\cdot}\|_2-\alpha)_+]^4
\end{equation}
and $p$ denotes the sampling rate, while both $\alpha$ and $\lambda$ can be viewed as tuning parameters.

Due to nonconvexity, standard methods such as gradient descent may converge to a local minimum of the above objective function. 
A series of papers in the literature including \citet{ge2016matrix, ge2017no} and \cite{chen2017memory} intend to understand the nonconvex geometry of \eqref{eq:obj_rect}. 
In particular, they focus on figuring out sufficient conditions for the ground-truth low rank matrix $\mtx{M}^\star$ as well as the size and pattern of $\Omega$, under which any local minimum $(\widehat{\mtx{X}}, \widehat{\mtx{Y}})$ of \eqref{eq:obj_rect} leads to an accurate estimate of $\mtx{M}^\star$ by $\widehat{\mtx{M}} = \widehat{\mtx{X}}\widehat{\mtx{Y}}^\top$.
For example, in the noiseless case where the noise matrix $\mtx{N}=\mtx{0}$, if $\mtx{M}^\star$ is a rank-$r$ well-conditioned positive semidefinite matrix and its eigenspace incoherence parameter defined in \cite{candes2009exact} is well-bounded, then it has been shown in \cite{chen2017memory} that with high probability, any local minimum of the PSD analog of \eqref{eq:obj_rect} yields exact low-rank recovery, i.e., there is no spurious local minimum, provided that the sampling pattern is i.i.d.~and the sampling rate satisfies $p\gtrsim (r^2 \log n)/n$. 

The goal of the present paper is to extend the nonconvex matrix completion objective \eqref{eq:obj_rect} as well as its geometric analysis to cases in which prior information on the low-rank matrix is known and explored. 
Two examples are subspace-constrained and skew-symmetric matrix completions. 
We aim to understand 
how to adapt the nonconvex objective \eqref{eq:obj_rect} to such cases with prior information, 
how to conduct corresponding geometrical analysis, and 
how to build a uniform framework to incorporate both cases as specific examples. 

A critical observation is that both examples can be represented in the form $\mtx{M}^\star = \boldsymbol{X}(\boldsymbol{\xi})\boldsymbol{Y}(\boldsymbol{\xi})^\top$, where the factors $\mtx{X}(\boldsymbol{\xi})$ and $\mtx{Y}(\boldsymbol{\xi})$ are linear mappings in $\vct{\xi} \in \mathbb{R}^d$. 
We elaborate on such parameterized factorizations as follows.

\begin{itemize}
\item Suppose $\mtx{M}^\star$ is known to be constrained in some pre-specified column and row spaces, with dimensions $s_1$ and $s_2$, respectively. 
Let  $\widetilde{\boldsymbol{U}}$ (and $\widetilde{\boldsymbol{V}}$) be an $n_1 \times s_1$ (and $n_2 \times s_2$) matrix whose columns form an orthogonal basis for the given column (and row) space constraint for $\mtx{M}^\star$. 
Given the rank of $\mtx{M}^\star$, we know there must exist some (not unique) $\boldsymbol{\Xi}_{A} \in \mathbb{R}^{s_1 \times r}$ and $\boldsymbol{\Xi}_{B} \in \mathbb{R}^{s_2 \times r}$, such that 
\[
\mtx{M}^\star = \left(\widetilde{\boldsymbol{U}}\boldsymbol{\Xi}_{A}\right) \left(\widetilde{\boldsymbol{V}}\boldsymbol{\Xi}_{B}\right)^\top.
\]
Denote by $\vct{\theta} = \text{vec}(\boldsymbol{\Theta}_{A}, \boldsymbol{\Theta}_{B})$ an $((s_1+s_2)r)$-dimensional vector that contains all entries in $\boldsymbol{\Theta}_{A}$ and $\boldsymbol{\Theta}_{B}$ (e.g., in the lexicographic order), and define the two linear mappings:
\begin{equation}
\label{eq:para_subspace}
 \boldsymbol{X}(\boldsymbol{\theta})=\widetilde{\boldsymbol{U}}\boldsymbol{\Theta}_{A} \in \mathbb{R}^{n_1 \times r}\quad \text{~and~} \quad \boldsymbol{Y}(\boldsymbol{\theta}) = \widetilde{\boldsymbol{V}}\boldsymbol{\Theta}_{B} \in \mathbb{R}^{n_2 \times r}. 
 \end{equation}
Without loss of generality, assume that both $\widetilde{\boldsymbol{U}}$ and $\widetilde{\boldsymbol{V}}$ consist of orthonormal basis, i.e., $\widetilde{\boldsymbol{U}}^\top \widetilde{\boldsymbol{U}}  = \boldsymbol{I}_{s_1}$ and $\widetilde{\boldsymbol{V}}^\top \widetilde{\boldsymbol{V}} = \boldsymbol{I}_{s_2}$. Then the above parameterized factorization becomes $\mtx{M}^\star = \boldsymbol{X}(\boldsymbol{\xi})\boldsymbol{Y}(\boldsymbol{\xi})^\top$ with $\vct{\xi} = \text{vec}(\boldsymbol{\Xi}_{A}, \boldsymbol{\Xi}_{B})$.

\item If $\mtx{M}^\star$ is an $n \times n$  rank-$r$ skew-symmetric matrix (which implies that $r$ is even), by the Youla decomposition \citep{youla1961normal}, it can be represented (not uniquely) as 
\[
\mtx{M}^\star = \boldsymbol{\Xi}_{A} \boldsymbol{\Xi}_{B}^\top - \boldsymbol{\Xi}_{B} \boldsymbol{\Xi}_{A}^\top, 
\]
where $\boldsymbol{\Xi}_A, \boldsymbol{\Xi}_B \in \mathbb{R}^{n \times \frac{r}{2}}$. Again, denote by $\vct{\theta} = \text{vec}(\boldsymbol{\Theta}_{A}, \boldsymbol{\Theta}_{B})$ a $(nr)$-dimensional vector that contains all entries in $\boldsymbol{\Theta}_{A}$ and $\boldsymbol{\Theta}_{B}$, and define the linear and homogeneous mappings 
\begin{equation}
\label{eq:para_skew}
 \boldsymbol{X}(\boldsymbol{\theta})= \left[\boldsymbol{\Theta}_A ,-\boldsymbol{\Theta}_B\right]  \in \mathbb{R}^{n \times r}\quad \text{~and~} \quad \boldsymbol{Y}(\boldsymbol{\theta}) = \left[\boldsymbol{\Theta}_B, \boldsymbol{\Theta}_A\right] \in \mathbb{R}^{n \times r}. 
 \end{equation}
We also have the factorization $\mtx{M}^\star = \boldsymbol{X}(\boldsymbol{\xi})\boldsymbol{Y}(\boldsymbol{\xi})^\top$ with $\vct{\xi} = \text{vec}(\boldsymbol{\Xi}_{A}, \boldsymbol{\Xi}_{B})$.
\end{itemize}

Therefore, both subspace-constrained matrix completion and skew-symmetric matrix completion are special cases of the following \emph{parameterized matrix completion} problem: 
If both factors for the low-rank matrix $\mtx{M}^\star = \mtx{X}\mtx{Y}^\top$ can be parameterized by linear mappings $\mtx{X}(\vct{\xi})$ and $\mtx{Y}(\vct{\xi})$, how can we estimate $\mtx{M}^\star$ from the partial and noisy observation $\mathcal{P}_{\Omega}(\boldsymbol{M}) = \mathcal{P}_{\Omega}(\boldsymbol{M}^\star + \boldsymbol{N})$?

Naturally, by substituting the parametric forms $\mtx{X} = \mtx{X}(\vct{\theta})$ and $\mtx{Y} = \mtx{Y}(\vct{\theta})$ into the nonconvex optimization \eqref{eq:obj_rect}, we obtain the following optimization with the argument $\vct{\theta}$:
\begin{equation}
\label{eq:obj_para}
\begin{split}
\tilde{f}(\boldsymbol{\theta}) \coloneqq & f(\boldsymbol{X}(\boldsymbol{\theta}),\boldsymbol{Y}(\boldsymbol{\theta})) \\
=&\frac{1}{2p}\|\mathcal{P}_{\Omega}(\boldsymbol{X}(\boldsymbol{\theta})\boldsymbol{Y}(\boldsymbol{\theta})^\top-\boldsymbol{M})\|_F^2 + \frac{1}{8} \|\boldsymbol{X}(\boldsymbol{\theta})^\top\boldsymbol{X}(\boldsymbol{\theta})-\boldsymbol{Y}(\boldsymbol{\theta})^\top\boldsymbol{Y}(\boldsymbol{\theta})\|_F^2 
\\
&+\lambda (G_{\alpha}(\boldsymbol{X}(\boldsymbol{\theta})) + G_{\alpha}(\boldsymbol{Y}(\boldsymbol{\theta}))).
\end{split}
\end{equation}

The nonconvex objective \eqref{eq:obj_para} can be applied to various specific parameterized matrix completion problems. For example, we can substitute the parametric forms \eqref{eq:para_subspace} into \eqref{eq:obj_para} to derive a nonconvex objective for subspace-constrained matrix completion, and substitute the parametric forms \eqref{eq:para_skew} into \eqref{eq:obj_para} to derive a nonconvex objective for skew-symmetric matrix completion. Details are elaborated on as follows:
\begin{itemize}
\item 
For subspace-constrained matrix completion, the parametric forms $\mtx{X}(\vct{\theta})$ and $\mtx{Y}(\vct{\theta})$ are defined in \eqref{eq:para_subspace} with $\vct{\theta} = \text{vec}(\boldsymbol{\Theta}_{A}, \boldsymbol{\Theta}_{B})$. Recall that we have assumed $\widetilde{\boldsymbol{U}}^\top \widetilde{\boldsymbol{U}}  = \boldsymbol{I}_{s_1}$ and $\widetilde{\boldsymbol{V}}^\top \widetilde{\boldsymbol{V}} = \boldsymbol{I}_{s_2}$. Then the parameterization \eqref{eq:para_subspace} implies the following
\[
\begin{cases}
\mtx{X}(\vct{\theta})\mtx{Y}(\vct{\theta})^\top = \widetilde{\boldsymbol{U}} \boldsymbol{\Theta}_A \boldsymbol{\Theta}_B^\top \widetilde{\boldsymbol{V}}^\top,
\\
~
\\
\boldsymbol{X}(\vct{\theta})^\top\boldsymbol{X}(\vct{\theta}) = \boldsymbol{\Theta}_A^\top \widetilde{\boldsymbol{U}}^\top \widetilde{\boldsymbol{U}} \boldsymbol{\Theta}_A = \boldsymbol{\Theta}_A^\top\boldsymbol{\Theta}_A,
\\ 
~
\\
\boldsymbol{Y}(\vct{\theta})^\top\boldsymbol{Y}(\vct{\theta}) = \boldsymbol{\Theta}_B^\top \widetilde{\boldsymbol{V}}^\top \widetilde{\boldsymbol{V}} \boldsymbol{\Theta}_B = \boldsymbol{\Theta}_B^\top\boldsymbol{\Theta}_B.
\end{cases}
\]
Substituting them into \eqref{eq:obj_para}, we have the objective function:
\begin{align}
\label{eq:obj_sub}
f_{subspace}(\boldsymbol{\Theta}_A, \boldsymbol{\Theta}_B)  \coloneqq & \frac{1}{2p}\|\mathcal{P}_{\Omega}(\widetilde{\boldsymbol{U}} \boldsymbol{\Theta}_A \boldsymbol{\Theta}_B^\top \widetilde{\boldsymbol{V}}^\top - \boldsymbol{M})\|_F^2 + \frac{1}{8} \|\boldsymbol{\Theta}_A^\top\boldsymbol{\Theta}_A-\boldsymbol{\Theta}_B^\top\boldsymbol{\Theta}_B\|_F^2 \nonumber
\\
& + \lambda (G_{\alpha}(\widetilde{\boldsymbol{U}} \boldsymbol{\Theta}_A) + G_{\alpha}(\widetilde{\boldsymbol{V}} \boldsymbol{\Theta}_B)).
\end{align}  
\item For skew-symmetric matrix completion, the parametric forms $\mtx{X}(\vct{\theta})$ and $\mtx{Y}(\vct{\theta})$ are defined in \eqref{eq:para_skew} with $\vct{\theta} = \text{vec}(\boldsymbol{\Theta}_{A}, \boldsymbol{\Theta}_{B})$. Straightforward calculation gives
\[
\begin{cases}
\mtx{X}(\vct{\theta}) \mtx{Y}(\vct{\theta})^\top = \mtx{\Theta}_A \mtx{\Theta}_B^\top - \mtx{\Theta}_B \mtx{\Theta}_A^\top,
\\
~
\\
\mtx{X}(\vct{\theta})^\top\mtx{X}(\vct{\theta}) = \begin{bmatrix} \mtx{\Theta}_A^\top\mtx{\Theta}_A & -\mtx{\Theta}_A^\top \mtx{\Theta}_B \\ - \mtx{\Theta}_B^\top \mtx{\Theta}_A & \mtx{\Theta}_B^\top \mtx{\Theta}_B \end{bmatrix},
\\
~
\\
\mtx{Y}(\vct{\theta})^\top \mtx{Y}(\vct{\theta}) = \begin{bmatrix} \mtx{\Theta}_B^\top\mtx{\Theta}_B & \mtx{\Theta}_B^\top \mtx{\Theta}_A \\  \mtx{\Theta}_A^\top \mtx{\Theta}_B & \mtx{\Theta}_A^\top \mtx{\Theta}_A \end{bmatrix}.
\end{cases}
\]
Substituting them into \prettyref{eq:obj_para}, we have the objective function
\begin{align}
\label{eq:obj_skew}
f_{\textrm{skew}}(\mtx{\Theta}_A, \mtx{\Theta}_B) =&\frac{1}{2p}\|\mathcal{P}_{\Omega}(\mtx{\Theta}_A \mtx{\Theta}_B^\top - \mtx{\Theta}_B \mtx{\Theta}_A^\top - \mtx{M})\|_F^2 +  \frac{1}{4}\|\mtx{\Theta}_A^\top\mtx{\Theta}_A-\mtx{\Theta}_B^\top\mtx{\Theta}_B\|_F^2
\\
&+\frac{1}{4}\|\mtx{\Theta}_A^\top\mtx{\Theta}_B+\mtx{\Theta}_B^\top\mtx{\Theta}_A\|_F^2 + 2\lambda G_{\alpha}(\left[\mtx{\Theta}_B, \mtx{\Theta}_A\right]).\nonumber
\end{align}
Here we use the fact $G_{\alpha}(\left[\mtx{\Theta}_B, \mtx{\Theta}_A\right]) = G_{\alpha}(\left[\mtx{\Theta}_A, -\mtx{\Theta}_B\right])$. 
\end{itemize}

Inspired by the previous studies in \citet{ge2016matrix, ge2017no} and \citet{chen2017memory} regarding standard matrix completion, 
we are interested in a unified geometric analysis for the generic nonconvex objective \eqref{eq:obj_para} for parameterized matrix completion. The crux for such extension hinges on two key assumptions on the parametric forms $(\mtx{X}(\vct{\theta}), \mtx{Y}(\vct{\theta}))$ and the ground truth $\mtx{M}^\star$. 
As suggested by the examples of subspace-constrained and skew-symmetric matrix completions, we first assume that both $\mtx{X}(\vct{\theta})$ and $\mtx{Y}(\vct{\theta})$ are linear mappings of $\vct{\theta}$. The second assumption regarding the relationship between $(\mtx{X}(\vct{\theta}), \mtx{Y}(\vct{\theta}))$ and the ground truth $\mtx{M}^\star$, referred to as \emph{correlated parametric factorization}, is novel in the literature of matrix completion, and plays the most essential role in our unified geometric analysis. Though this assumption has a neat mathematical form, it cannot be easily explained by non-mathematical intuitions. We defer its formal definition to Section \ref{sec:assumptions}.

The most intriguing part in this work is that the non-intuitive assumption of \emph{correlated parametric factorization} holds for various examples including the aforementioned subspace-constrained factorization \eqref{eq:para_subspace} and skew-symmetric factorization \eqref{eq:para_skew}. In fact, bringing both examples into a unified framework by this neat assumption is unexpected to us, and it is still unclear whether there is some geometrical explanation for such unification. Nevertheless, we verify the assumption of correlated parametric factorization for those two examples in Sections \ref{sec:rotation_positive_subspace} and \ref{sec:rotation_positive_skew} respectively with very different algebraic arguments.

\subsection{Related Work}
As with most high-dimensional problems where low-complexity structures need to be explored and exploited for statistical learning and inferences, low-rank structures are a common assumption for deriving and justifying matrix completion algorithms. 
By imposing nuclear norm regularization to recover low-rank structures \citep{RechtFazelParrilo2010}, convex optimization methods have been widely used in the literature of matrix completion, and their theoretical properties have been well studied; see, e.g., \citet{candes2009exact, candes2010power, CLMW2011, recht2011simpler, gross2011recovering, negahban2012restricted, HKZ11, SunZhang2012, Koltchinskii2010noisy, li2013compressed}.

Though convex optimization methods could have near-optimal theoretical guarantees for matrix completion under certain incoherence conditions, they could be unscalable to large data matrices whose dimensions are in hundreds of thousands. 
In contrast, nonconvex optimization methods for low-rank recovery have been proposed and analyzed in the literature due to computational convenience. Examples include optimization over Grassmann manifolds \citep{keshavan2010matrix,keshavan2010matrixnoise}, singular value projections (SVP) \citep{jain2010guaranteed, ding2020leave}, alternating minimization \citep{jain2013low}, penalized, projected, or thresholded gradient descent \citep{sun2016guaranteed, chen2015fast, cai2016optimal, zheng2016convergence, yi2016fast}, vanilla gradient descent \citep{candes2015phase, ma2020implicit, chen2020nonconvex}, etc. 
From geometric perspective, ``no spurious local minimum'' results have been established in \citet{sun2018geometric, ge2016matrix, ge2017no, chen2017memory}, and the present work belongs to this line of studies.

Prior information and special structures have also been explored and employed in the literature of matrix completion. Subspace-constrained matrix completion has been shown to be effective in improving the performances of collaborative filtering, and has various methodological and theoretical developments \citep{JainInductive, Xu2013, yi2013semi, natarajan2014inductive, chen2015incoherence,  10.1145/2939672.2939809, eftekhari2018weighted}. In fact, without the penalization terms, the least squares term in \eqref{eq:obj_sub} can be represented in form of the nonconvex objective studied in \cite{JainInductive}, though their distributional assumptions on rank-one sampling do not apply to our settings. Motivated by applications such as pairwise ranking, skew-symmetric structures have also been exploited for matrix completion in the literature. 
See, e.g., \cite{JLYY2011}, \cite{Gleich2011rank}, and \cite{chatterjee2015matrix}. For example, it has been shown in \citet[Theorem 3]{Gleich2011rank} that if the initial input is skew-symmetric, singular value projection (SVP) \citep{jain2010guaranteed} maintains the skew-symmetry in all iterations.

\subsection{Notation}

Throughout this paper, bold uppercase/lowercase characters denote matrices/vectors, respectively. 
For a given matrix $\boldsymbol{A}$, its $(i,j)$-th entry, $i$-th row, and $j$-th column are denoted by $A_{i,j}$, $\boldsymbol{A}_{i,\cdot}$, and $\boldsymbol{A}_{\cdot,j}$, respectively. Its spectral, Frobenius, and $\ell_{2,\infty}$ norms are denoted by $\|\boldsymbol{A}\|,\|\boldsymbol{A}\|_F$ and $\|\boldsymbol{A}\|_{2,\infty} := \max_{i} \|\mtx{A}_{i, \cdot}\|_2^2$, respectively. Denote by $\colspan(\boldsymbol{A})$/$\rowspan(\boldsymbol{A})$ the column/row space of $\boldsymbol{A}$. Denote by $\boldsymbol{P}_{\boldsymbol{A}}$ the Euclidean projector onto $\colspan(\boldsymbol{A})$. Denote $\boldsymbol{A}\succeq \boldsymbol{0}$ if $\boldsymbol{A}$ is a symmetric or Hermitian positive semidefinite matrix. For any two matrices $\mtx{A}$ and $\mtx{B}$ of the same dimensions, their matrix inner product is denoted by $\langle \boldsymbol{A},\boldsymbol{B} \rangle = \operatorname{trace}(\boldsymbol{A}^\top \boldsymbol{B})= \sum_i \sum_j A_{i,j}B_{i,j}$, and their Hadamard/entrywise product is denoted by $\boldsymbol{A}\circ \boldsymbol{B}$ with entries $[\boldsymbol{A}\circ \boldsymbol{B}]_{i,j} = A_{i,j}B_{i,j}$. For any two matrices $\mtx{A}$ and $\mtx{B}$, $\operatorname{vec}(\boldsymbol{A},\boldsymbol{B})$ denotes a vector consisting of all entries in $\boldsymbol{A}$ and $\boldsymbol{B}$ in some fixed order. Denote by $\boldsymbol{J}_{n_1\times n_2}$ (or $\boldsymbol{J}$ when the dimensions are clear in the context) the $n_1 \times n_2$ matrix with all entries equal to one. Denote by $\mathsf{O}(r)$ the set of $r \times r$ orthogonal matrices. 
Let $n_{\min} \coloneqq \min\{n_1,n_2\}$ and $n_{\max} \coloneqq \max\{n_1,n_2\}$. 
Finally, denote by $C_0,C_1,\ldots$ and $C_v, C_c,\ldots$ some fixed positive absolute constants. 

\subsection{Organization of the Paper}

The rest of the paper is organized as follows: In Section \ref{sec:theory}, we first introduce the key assumptions regarding the parameterized factorization formally, and then introduce our main results on the low-rank estimation with any local minimum, followed by corollaries for subspace-constrained and skew-symmetric matrix completions. 
Some results of numerical simulations are presented in Section \ref{sec:experiments} so as to illustrate our theoretical findings. We give a summary of our contributions and propose some open questions for future research in Section \ref{sec:discussion}. 
All proofs can be found in Sections \ref{sec:proof_main} with technical details deferred to the appendices.

\section{Main Results} 
\label{sec:theory}
The major contribution of this paper is on deriving theoretical properties of parameterized matrix completion with nonconvex optimization \eqref{eq:obj_para}. 
We first introduce key assumptions on the low-rank matrix $\mtx{M}^\star$, the parameterized factors $(\mtx{X}(\vct{\theta}),\mtx{Y}(\vct{\theta}))$, and the observation support $\Omega$. 
Under these assumptions, we show in Section \ref{sec:theoretical_results} that the geometric analysis in \citet{chen2017memory} for PSD matrix completion can be extended to parameterized matrix completion.
In particular, we give a uniform upper bound for the estimation error $\|\boldsymbol{X}(\hat{\boldsymbol{\xi}})\boldsymbol{Y}(\hat{\boldsymbol{\xi}})^\top - \mtx{M}^\star\|_F^2$ with any local minimum $\hat{\boldsymbol{\xi}}$ of the nonconvex objective \eqref{eq:obj_para}. Corollaries of our meta-theorem give conditions for consistent 
subspace-constrained and skew-symmetric matrix completions via nonconvex optimization.

\subsection{Assumptions}
\label{sec:assumptions}
Let us now introduce the key assumptions on the parameterized factors $(\mtx{X}(\vct{\theta}), \mtx{Y}(\vct{\theta}))$ mentioned earlier. We first assume that both factors are linear mappings in $\vct{\theta}$.

\begin{assumption}[Linearity]
\label{ass:homogeneous}
Both $\mtx{X}(\vct{\theta}) \in \mathbb{R}^{n_1 \times r}$ and $\mtx{Y}(\vct{\theta}) \in \mathbb{R}^{n_2 \times r}$ are linear mappings in $\vct{\theta}$, i.e. $\boldsymbol{X}(t_1\boldsymbol{\theta}_1 + t_2\boldsymbol{\theta}_2) = t_1\boldsymbol{X}(\boldsymbol{\theta}_1) + t_2 \boldsymbol{X}(\boldsymbol{\theta}_2)$  and  $\boldsymbol{Y}(t_1 \boldsymbol{\theta}_1 + t_2\boldsymbol{\theta}_2) = t_1\boldsymbol{Y}(\boldsymbol{\theta}_1) + t_2\boldsymbol{Y}(\boldsymbol{\theta}_2)$ for all $\boldsymbol{\theta}_1,\boldsymbol{\theta}_2\in\mathbb{R}^d$ and $t_1, t_2 \in \mathbb{R}$. Here $d$ is the intrinsic dimension of $\vct{\theta}$.
\end{assumption}

The next assumption, referred to as the \emph{Correlated Parametric Factorization}, plays an essential role in the geometric analysis of the nonconvex objective \eqref{eq:obj_para}, verification of which for subspace-constrained factorization \eqref{eq:para_subspace} and skew-symmetric factorization \eqref{eq:para_skew} will be stated in later subsections.

\begin{assumption}[Correlated Parametric Factorization of $ \mtx{M}^\star$]  
\label{ass:ideal_factorization}
The rank-$r$ matrix $\mtx{M}^\star $ and the parameterized factors $(\mtx{X}(\vct{\theta}), \mtx{Y}(\vct{\theta}))$ defined in Assumption \ref{ass:homogeneous} are said to satisfy the correlated parameterized factorization, if for any $\vct{\theta} \in \mathbb{R}^d$, there exists $\vct{\xi}  \in \mathbb{R}^d$ (not necessarily unique), such that
\begin{equation}
\label{eq:balanced parameterized factorization}
\begin{cases}
\mtx{M}^\star  = \mtx{X}(\vct{\xi} )\mtx{Y}(\vct{\xi} )^\top,
 \\
 \mtx{X}(\vct{\xi} )^\top \mtx{X}(\vct{\xi} ) = \mtx{Y}(\vct{\xi} )^\top \mtx{Y}(\vct{\xi} ), 
 \\
 \mtx{X}(\vct{\theta})^\top \mtx{X}(\vct{\xi} ) + \mtx{Y}(\vct{\theta})^\top \mtx{Y}(\vct{\xi} ) \succeq \mtx{0}.
 \end{cases}
\end{equation} 
\end{assumption}

Let us now introduce models of observation patterns $\Omega \subset [n_1] \times [n_2]$. We assume that the entries are observed independently with certain probability $p$ for the following two scenarios of rectangular and square matrices, respectively.

\begin{model}
\label{mod:sampling1}
For rectangular matrix completion, the index set $\Omega$ is assumed to follow the independent $\text{Ber}(p)$ model, i.e., each entry is sampled independently with probability $p$.
\end{model}

\begin{model}
\label{mod:sampling2}
For square matrix completion ($n_1 = n_2 := n$), the index set $\Omega$ is assumed to follow the off-diagonal symmetric independent $\text{Ber}(p)$ model, i.e., $\Omega$ is the support of symmetric off-diagonal entries that are sampled independently with probability $p$. No diagonal entries are included in $\Omega$.
\end{model}



\subsection{A Meta-Theorem}
\label{sec:theoretical_results}
Our main theorem is concerned with the conditions for accurate parametric matrix completion via the nonconvex objective \eqref{eq:obj_para}.

\begin{theorem}
\label{thm:main}
Assume that $\mtx{M}^\star$ has rank $r$ and its reduced singular value decomposition (SVD) is $\mtx{M}^\star  = \mtx{U}^{\star} \mtx{\Lambda} {\mtx{V}^{\star} }^\top$, 
where $\boldsymbol{U}^{\star} \in\mathbb{R}^{n_1 \times r}$, $\boldsymbol{V}^{\star} \in\mathbb{R}^{n_2 \times r}$ and $\mtx{\Lambda} = \text{diag}(\sigma_1, \sigma_2, \ldots, \sigma_r)$ with $\sigma_1\geqslant \sigma_2 \geqslant \ldots \geqslant \sigma_r >0$\footnote{In this paper, (reduced) SVD of matrices is repeatedly used. Note that (reduced) SVD may not be unique due to possible multiplicity in singular values. In that case, we simply choose one (reduced) SVD and keep it fixed throughout the discussion.}. The condition number of $\mtx{\Lambda}$ is denoted by $\kappa = \sigma_1/\sigma_r$. Moreover, following \citet{candes2009exact}, the incoherence parameter $\mu$ for $\mtx{M}^\star$ determined by its row and column spaces is defined as 
\begin{equation}
\label{eq:incoherence}
\mu = \max\left\{  \frac{n_1}{r} \left\| \mtx{U}^{\star} \right\|_{2, \infty}^2 ,  \frac{n_2}{r}\left\| \mtx{V}^{\star} \right\|_{2, \infty}^2 \right\}.
\end{equation}
Suppose that $\mtx{M}^\star$, $\mtx{X}(\vct{\theta})$ and $\mtx{Y}(\vct{\theta})$ satisfy Assumptions \ref{ass:homogeneous} and \ref{ass:ideal_factorization}, and that $\Omega$, the support of observed entries, follows either Model \ref{mod:sampling1} or \ref{mod:sampling2}.
Moreover, let the sampling rate $p$ and the tuning parameters $\alpha$ and $\lambda$ in \eqref{eq:obj_rect} satisfy the following inequalities:
\begin{equation} 
  \label{eq:main_rate_tuning_ass}
\begin{cases}
  p \geqslant  C_1 \max \left\{  \frac{1}{n_{\min}}\mu r \log n_{\max}, \frac{n_{\max}}{n_{\min}^2}\mu^2 r^2\kappa^2   \right\},
  \\
  ~
  \\
  C_2\sqrt{n_{\max}/{p}} \leqslant \lambda \leqslant 10C_2\sqrt{{n_{\max}}/{p}},
  \\
  ~
  \\
  C_2\sqrt{{\mu r \sigma_1}/{n_{\min}}} \leqslant \alpha \leqslant 10C_2\sqrt{{\mu r \sigma_1}/{n_{\min}}}, 
\end{cases}
\end{equation}
where $n_{\max} \coloneqq \max(n_1, n_2)$ and $n_{\min} \coloneqq \min(n_1, n_2)$. Then, on an event $E$ with probability $\mathbb{P}[E] \geqslant  1 - (n_1+n_2)^{-3}$, any local minimum $\hat{\vct{\xi}}$ of \eqref{eq:obj_para} satisfies
 \[
 \|\mtx{M}^\star - \widehat{\boldsymbol{M}}\|_F^2 \leqslant  \frac{C_3 r}{p^2} \psi^2, 
\]
where 
\begin{equation}\label{eq:psi_def}
  \psi \coloneqq \max_{\boldsymbol{\theta}_1,\boldsymbol{\theta}_2 \in\mathbb{R}^d } \|\boldsymbol{P}_{\boldsymbol{X}(\boldsymbol{\theta}_1)}\mathcal{P}_{\Omega}(\boldsymbol{N}) \boldsymbol{P}_{\boldsymbol{Y}(\boldsymbol{\theta}_2)}\|.
\end{equation}
Here $C_1,C_2,C_3$ are fixed absolute constants. Recall that $\boldsymbol{P}_{\boldsymbol{A}}$ is the Euclidean projection matrix onto $\colspan(\boldsymbol{A})$. In particular, if there is no noise, i.e., $\mtx{N} = \mtx{0}$, then on the same event any local minimum $\vct{\hat{\xi}}$ gives exact low-rank recovery $\widehat{\mtx{M}} = \mtx{M}$, that is, there is no spurious local minimum.
\end{theorem}

\subsection{Implication for PSD Matrix Completion}
Besides applications in some specific parametric matrix completion problems, 
existing geometric analysis for nonconvex matrix completion for more standard cases like general rectangular or positive semidefinite (PSD) matrices, can be implied by our meta-theorem. 

Consider, for example, PSD matrix completion where we have $n_1 = n_2 = n$, and the underlying rank-$r$ matrix can be decomposed as $\mtx{M}^\star = \mtx{\Xi}_0 \mtx{\Xi}_0^\top$ for some $\boldsymbol{\Xi}_0 \in\mathbb{R}^{n\times r}$. 
For any $\boldsymbol{\Theta}\in\mathbb{R}^{n\times r}$, denote $\boldsymbol{\theta} := \operatorname{vec}(\boldsymbol{\Theta})$ and define the linear mappings $\boldsymbol{X}(\boldsymbol{\theta}) = \boldsymbol{Y}(\boldsymbol{\theta}) = \boldsymbol{\Theta}$. 
The parameterized factorization for rank-$r$ PSD matrices becomes $\mtx{X}(\vct{\theta}) \mtx{Y}(\vct{\theta})^\top = \boldsymbol{\Theta}\boldsymbol{\Theta}^\top$, and the corresponding parameterized factorization for the ground truth is $\mtx{M}^\star = \mtx{X}(\vct{\xi}_0) \mtx{Y}(\vct{\xi}_0)^\top$ for $\vct{\xi}_0:= \operatorname{vec}(\boldsymbol{\Xi}_0)$.

It is obvious that $\boldsymbol{X}(\boldsymbol{\theta}) = \boldsymbol{Y}(\boldsymbol{\theta}) = \boldsymbol{\Theta}$ are linear mappings in $\vct{\theta}$, and so Assumption \ref{ass:homogeneous} holds.
Let us now see how to construct $\mtx{\Xi}$ based on any given $\mtx{\Theta}$, such that the requirements in Assumption \ref{ass:ideal_factorization} can hold simultaneously. 
The second equality in \eqref{eq:balanced parameterized factorization} holds since $\mtx{X}(\vct{\xi}) = \mtx{Y}(\vct{\xi}) = \mtx{\Xi}$. 
To ensure the first equality, we let $\mtx{\Xi} = \mtx{\Xi}_0 \mtx{R}$ for some underdetermined orthogonal matrix $\mtx{R}$, which gives $\mtx{\Xi} \mtx{\Xi}^\top = \mtx{\Xi}_0 \mtx{\Xi}_0^\top = \mtx{M}^\star$. The remaining question is how to choose $\mtx{R}$ such that $\mtx{\Theta}^\top \mtx{\Xi}$ is PSD. This can be easily obtained by considering the SVD of $\mtx{\Theta}^\top \mtx{\Xi}_0$. 
It is interesting that the construction of such $\mtx{R}$ plays an important role in the literature of PSD matrix completion. 
See, e.g., \citet[Lemma 1]{chen2015fast}.

The nonconvex objective \eqref{eq:obj_para} for PSD matrix completion takes the form
\begin{align}
\label{eq:obj_psd}
f_{psd}(\boldsymbol{\Theta}) \coloneqq \frac{1}{2p} \left\|\mathcal{P}_{\Omega}(\boldsymbol{\Theta}\boldsymbol{\Theta}^\top-\boldsymbol{M})\right\|_F^2  + 2\lambda G_{\alpha}(\boldsymbol{\Theta}),
\end{align}
which is exactly the objective studied in \citet{ge2016matrix,ge2017no} and  \citet{chen2017memory}. 
Implied by Theorem \ref{thm:main}, with high probability, there is no spurious local minimum of \eqref{eq:obj_psd}, provided that the tuning parameters are suitably chosen and the sampling rate satisfies $p \geqslant  \frac{C_1}{n} \max \left\{ \mu r \log n, \mu^2 r^2\kappa^2  \right\}$, which is exactly the same as the state-of-the-art result in \citet{chen2017memory}. 
Furthermore, consider the special noisy case in which the entries of noise matrix $\boldsymbol{N}$ are i.i.d.\ Gaussian random variables with mean $0$ and variance $\sigma^2$. Then $\|\mathcal{P}_{\Omega}(\boldsymbol{N})\|^2  = O( (np +\log^2 n)\sigma^2 )$ (see, e.g., \citet[Lemma 11]{chen2015fast}). Theorem \ref{thm:main} implies the estimation error bound $\|\mtx{M}^\star - \widehat{\boldsymbol{\Theta}} \widehat{\boldsymbol{\Theta}}^\top \|_F^2 = O((\frac{nr}{p}+ \frac{r\log^2 n}{p^2}) \sigma^2 )$, which matches the state-of-the-art results in the literature of noisy matrix completion. 
See, e.g., \citet{keshavan2010matrixnoise}, \citet{chen2015fast} and \citet{ma2018implicit}.

In the next two subsections, we apply Theorem \ref{thm:main} to subspace-constrained matrix completion \eqref{eq:obj_sub} and skew-symmetric matrix completion \eqref{eq:obj_skew}, respectively.

\subsection{Nonconvex Subspace-Constrained Matrix Completion}
\label{sec:subspace_constraint_matrix}
For subspace-constrained matrix completion, the parametric factors have the form \eqref{eq:para_subspace}, i.e.  $\mtx{X}(\vct{\theta})=\widetilde{\mtx{U}}\mtx{\Theta}_{A}$ and $\mtx{Y}(\vct{\theta}) = \widetilde{\mtx{V}}{\mtx{\Theta}}_B$, where $\vct{\theta} = \text{vec}(\mtx{\Theta}_{A}, \mtx{\Theta}_{B})$. Recall that we always assume $\widetilde{\boldsymbol{U}}^\top \widetilde{\boldsymbol{U}}  = \boldsymbol{I}_{s_1}$ and $\widetilde{\boldsymbol{V}}^\top \widetilde{\boldsymbol{V}} = \boldsymbol{I}_{s_2}$. 
Here $\mtx{X}(\vct{\theta})$ and $\mtx{Y}(\vct{\theta})$ are linear mappings, thereby satisfying Assumption \ref{ass:homogeneous}. 
The verification of Assumption \ref{ass:ideal_factorization} is summarized as the following lemma, the proof of which is deferred to Section \ref{sec:rotation_positive_subspace}.


\begin{lemma}
\label{lem:rotation_positive_subspace}
Let $\mtx{M}^\star  \in \mathbb{R}^{n_1\times n_2}$ be a rank-$r$ matrix whose column space and row space are constrained in $\colspan(\widetilde{\boldsymbol{U}})$ and $\colspan(\widetilde{\boldsymbol{V}})$, where $\widetilde{\boldsymbol{U}}^\top \widetilde{\boldsymbol{U}}  = \boldsymbol{I}_{s_1}$ and $\widetilde{\boldsymbol{V}}^\top \widetilde{\boldsymbol{V}} = \boldsymbol{I}_{s_2}$. Then the parameterized factors $\boldsymbol{X}(\boldsymbol{\theta})$, $\boldsymbol{Y}(\boldsymbol{\theta})$ defined in \eqref{eq:para_subspace} satisfy Assumption \ref{ass:ideal_factorization} with $\mtx{M}^\star \in \mathbb{R}^{n_1\times n_2}$. 
\end{lemma}

Theorem \ref{thm:main} implies the following corollary for nonconvex subspace-constrained matrix completion with objective \eqref{eq:obj_sub}.

\begin{corollary}
\label{cor:subspace}
Let $\mtx{M}^\star  \in \mathbb{R}^{n_1\times n_2}$ be a rank-$r$ matrix. The incoherence parameter $\mu$ and condition number $\kappa$ are defined in Theorem \ref{thm:main}. Assume that the columns of $\widetilde{\boldsymbol{U}} \in \mathbb{R}^{n_1 \times s_1}$ constitute an orthonormal basis of the column space constraint for $\mtx{M}^\star $, while the columns of $\widetilde{\boldsymbol{V}} \in \mathbb{R}^{n_2 \times s_2}$ constitute an orthonormal basis of the row space constraint. The support of observation, $\Omega$, is assumed to follow Model \ref{mod:sampling1}. The entries of the noise matrix $\boldsymbol{N}$ are i.i.d.\ centered sub-exponential random variables satisfying the Bernstein condition \citep[(2.15)]{wainwright2019high} with parameter $b$ and variance $\nu^2$. 

If the sampling rate $p$ and the tuning parameters $\alpha, \lambda$ satisfy \eqref{eq:main_rate_tuning_ass}. Then, uniformly on an event $E_{subspace}$ with probability $\mathbb{P}[E_{subspace}] \geqslant  1 - 2(n_1+n_2)^{-3}$, any local minimum $(\widehat{\boldsymbol{\Xi}}_A, \widehat{\boldsymbol{\Xi}}_B)$ of $f_{subspace}(\boldsymbol{\Theta}_A,\boldsymbol{\Theta}_B)$ in \eqref{eq:obj_sub} satisfies:
\begin{equation}
\label{eq:subspace_estimate}
  \| \widetilde{\boldsymbol{U}} \widehat{\boldsymbol{\Xi}}_A\widehat{\boldsymbol{\Xi}}_B^\top \widetilde{\boldsymbol{V}}^\top - \boldsymbol{M}^{\star}\|_F^2 \leqslant   \frac{ C_4 r}{p^2} \left( \nu^2 p (s_1 + s_2)  \log(n_1 + n_2) + b^2 \frac{\mu_{\widetilde{U}} \mu_{\widetilde{V}}s_1s_2 }{n_1n_2} \log^2(n_1 + n_2)\right).
\end{equation}
Here ${\mu}_{\tilde{U}} = \frac{n_1}{s_1} \|\tilde{\mtx{U}} \|_{2, \infty}^2$, ${\mu}_{\tilde{V}} = \frac{n_2}{s_2} \|\tilde{\mtx{V}} \|_{2, \infty}^2$, and $C_4$ is some fixed positive absolute constant. 
\end{corollary}
 

To the best of our knowledge, existing investigations in the literature on matrix completion with subspace constraints are majorly focused on the noiseless case \citep{yi2013semi,Xu2013,chen2015incoherence,JainInductive, eftekhari2018weighted}, while error rates for noisy recovery are little-studied. Consider the scenario in which the noise matrix $\boldsymbol{N}$ consists of i.i.d.\ $\mathcal{N}(0, \sigma^2)$ entries. 
This gives $b = \sigma$ and variance $\nu^2 = \sigma^2$. For simplicity of discussion, assume $s_1 = s_2 \coloneqq s$, $n_1 = n_2 \coloneqq n$, $\mu = O(1)$, $\mu_{\widetilde{U}}=O(1)$, $\mu_{\widetilde{V}}=O(1)$ and $\kappa = O(1)$. Then Corollary \ref{cor:subspace} implies that as long as $p \gtrsim \frac{1}{n} \max \left\{ r \log n, r^2  \right\}$, with high probability, $\| \widetilde{\boldsymbol{U}} \widehat{\boldsymbol{\Xi}}_A\widehat{\boldsymbol{\Xi}}_B^\top \widetilde{\boldsymbol{V}}^\top - \boldsymbol{M}^{\star}\|_F^2 \lesssim \sigma^2 s r (\log n) /p$ holds uniformly for any local minimum. As explained in the previous subsection, the error rates for standard matrix completion are $O(\sigma^2 {n r}/{p})$. Therefore, Corollary \ref{cor:subspace} indicates that the estimation error can be significantly reduced by exploring and employing subspace constraints.

In the noiseless case, we should admit that the sufficient condition on the sampling rate for exact matrix completion, $p \gtrsim \frac{1}{n} \max \{ r \log n, r^2\}$, may be suboptimal. In fact, fewer samples are sufficient for exact low-rank recovery for convex optimization methods \citep{yi2013semi,Xu2013,chen2015incoherence} as well as alternating minimization \citep{JainInductive}. 
Given our aim is to establish a unified methodological and theoretical framework for parameterized matrix completion via nonconvex optimization, improving sample size conditions for special cases is beyond the scope of the paper, but we are interested in studying this issue in future work.

\subsection{Nonconvex Skew-Symmetric Matrix Completion} 
\label{sec:skew_symmetric}
For skew-symmetric matrix completion, the parameterized factors are defined in \eqref{eq:para_skew}, which are  $\boldsymbol{X}(\boldsymbol{\theta})= \left[\boldsymbol{\Theta}_A ,-\boldsymbol{\Theta}_B\right]$ and $\boldsymbol{Y}(\boldsymbol{\theta}) = \left[\boldsymbol{\Theta}_B, \boldsymbol{\Theta}_A\right]$ with $\vct{\theta} = \text{vec}(\boldsymbol{\Theta}_{A}, \boldsymbol{\Theta}_{B})$.
The factors $\mtx{X}(\vct{\theta})$ and $\mtx{Y}(\vct{\theta})$ are linear mappings and thereby satisfying Assumption \ref{ass:homogeneous}. 
Assumption \ref{ass:ideal_factorization} is verified through the following result with the proof deferred to Section \ref{sec:rotation_positive_skew}.

\begin{lemma}
 \label{lem:rotation_positive_skew}
Let $\mtx{M}^\star$ be a rank-$r$ skew-symmetric matrix. Then, the parameterization $(\boldsymbol{X}(\boldsymbol{\theta}), \boldsymbol{Y}(\boldsymbol{\theta}))$ defined in \eqref{eq:para_skew} and $\mtx{M}^\star$ satisfy Assumption \ref{ass:ideal_factorization}.
\end{lemma}

Theorem \ref{thm:main} implies the following noisy low-rank recovery result for skew-symmetric matrix completion via nonconvex objective \eqref{eq:obj_skew}.

\begin{theorem}
\label{thm:skew}
Let $\mtx{M}^\star  \in \mathbb{R}^{n \times n}$ be a rank-$r$ skew-symmetric matrix. The incoherence parameter $\mu$ and the condition number $\kappa$ are defined in Theorem \ref{thm:main}. The support of the observed entries $\Omega$ is assumed to follow from Model \ref{mod:sampling2}. The skew-symmetric noise matrix $\boldsymbol{N}$ consists of i.i.d.\ upper triangular entries, which are centered sub-exponential random variables satisfying the Bernstein condition with parameter $b$ and variance $\nu^2$. Suppose that the sampling rate $p$ and the tuning parameters $\alpha$ and $\lambda$ satisfy \eqref{eq:main_rate_tuning_ass}. Then, uniformly on an event $E_{skew}$ with probability $\mathbb{P}[E_{skew}] \geqslant  1 - 2n^{-3}$, any local minimum $(\widehat{\boldsymbol{\Xi}}_A, \widehat{\boldsymbol{\Xi}}_B)$ of $f_{\textrm{skew}}(\boldsymbol{\Theta}_A,\boldsymbol{\Theta}_B)$ defined in \eqref{eq:obj_skew} satisfies
\[
  \| \widehat{\boldsymbol{\Xi}}_A \widehat{\boldsymbol{\Xi}}_B^\top - \widehat{\boldsymbol{\Xi}}_B \widehat{\boldsymbol{\Xi}}_A^\top - \mtx{M}^{\star}\|_F^2 \leqslant  \frac{ C_5 r}{p^2} \left(\nu^2 p n \log n+ b^2   \log^2 n \right).
 \]
Where $C_5$ is a fixed positive absolute constant. 
 \end{theorem}

As with the discussion in Section \ref{sec:subspace_constraint_matrix}, if the upper triangular part of noise matrix $\boldsymbol{N}$ consists of i.i.d.\ Gaussian random variables with mean $0$ and variance $\sigma^2$, and the sampling rate satisfis $p \gtrsim \frac{1}{n} \max \left\{ r \log n, r^2  \right\}$, then the estimation error satisfies $\| \widehat{\boldsymbol{\Xi}}_A \widehat{\boldsymbol{\Xi}}_B^\top - \widehat{\boldsymbol{\Xi}}_B \widehat{\boldsymbol{\Xi}}_A^\top - \mtx{M}^{\star}\|_F^2 = O(\sigma^2 nr (\log n) /p)$, which is comparable to the aforementioned error rate $O(\sigma^2 nr /p)$ for PSD matrix completion up to a logarithmic factor.

\section{Experiments}
\label{sec:experiments}
In this section, some numerical results are shown to demonstrate the performance of the proposed nonconvex optimization for subspace-constrained matrix completion \eqref{eq:obj_sub} and skew-symmetric matrix completion \eqref{eq:obj_skew}.

In all simulations, the sampling rate $p$ is replaced with the empirical analog $\hat{p}\coloneqq \frac{|\Omega|}{n_1n_2}$. The hyperparameters are set as $\lambda = 100\sqrt{(n_1+n_2)\hat{p}}$ and $\alpha = 100$. In both examples, the nonconvex objective \eqref{eq:obj_para} is minimized by gradient descent, with initialization for $\vct{\theta}$ consisting with i.i.d.\ standard normal entries. At each step of the gradient descent, the step size is selected through line search. 
To be specific, at each update of $\vec{\theta}$, the step size is set to be $\max\{2^{-k},10^{-10}\}$ for $k \coloneqq \min\{t \mid t = 0,1,2,3,\cdots,\tilde{f}(\boldsymbol{\theta}-2^{-t}\nabla\tilde{f}(\boldsymbol{\theta})) \leqslant  \tilde{f}(\boldsymbol{\theta}) \}$. 
The gradient descent iteration is terminated either after $500$ iterations or as soon as the update on $\vct{\theta}$ satisfied $\|\nabla\tilde{f}(\boldsymbol{\theta})\|_2^2 \leqslant 10^{-10}$. 

\subsection{Subspace-Constrained Matrix Completion}
\label{sec:sub_exp}
For simulation of nonconvex optimization for subspace-constrained matrix completion \eqref{eq:obj_sub}, the dimensions of $\boldsymbol{M}^{\star}$ are set to be $n_1 = n_2 = 500$ and its rank is set as $r=2$. We also set the dimensions of the prior column and row subspaces as $s_1 = s_2 \coloneqq s \leq 40$. In preparation for the construction of the column and row subspace constraints with various dimensions, we generate $[\boldsymbol{u}_1,\ldots,\boldsymbol{u}_{40}]$ and $[\boldsymbol{v}_1,\ldots,\boldsymbol{v}_{40}]$ according to the Haar measure on the manifold of $500 \times 40$ orthonormal basis matrices. The matrix to recover is fixed as $\boldsymbol{M}^{\star} = \boldsymbol{u}_1\boldsymbol{v}_1^{\top} + \boldsymbol{u}_2\boldsymbol{v}_2^{\top}$, which gives $\|\boldsymbol{M}^{\star}\|_F^2 = 2$. In the case of noisy matrix completion, the noise matrix $\boldsymbol{N}$ consists of i.i.d.\ $\mathcal{N}(0, \sigma^2)$ entries with $\sigma^2 = {1}/{500^2}$, which gives $\mathbb{E}\|\boldsymbol{N}\|_F^2  = 1$, and hence the SNR is 
$\|\boldsymbol{M}^{\star}\|_F^2/\mathbb{E}\|\boldsymbol{N}\|_F^2 = 2$. In the noiseless case we set $\mtx{N} = \mtx{0}$.

\begin{figure}[!h]
\centering
\begin{subfigure}[b]{\textwidth}
\centering
\includegraphics[width=0.8\textwidth]{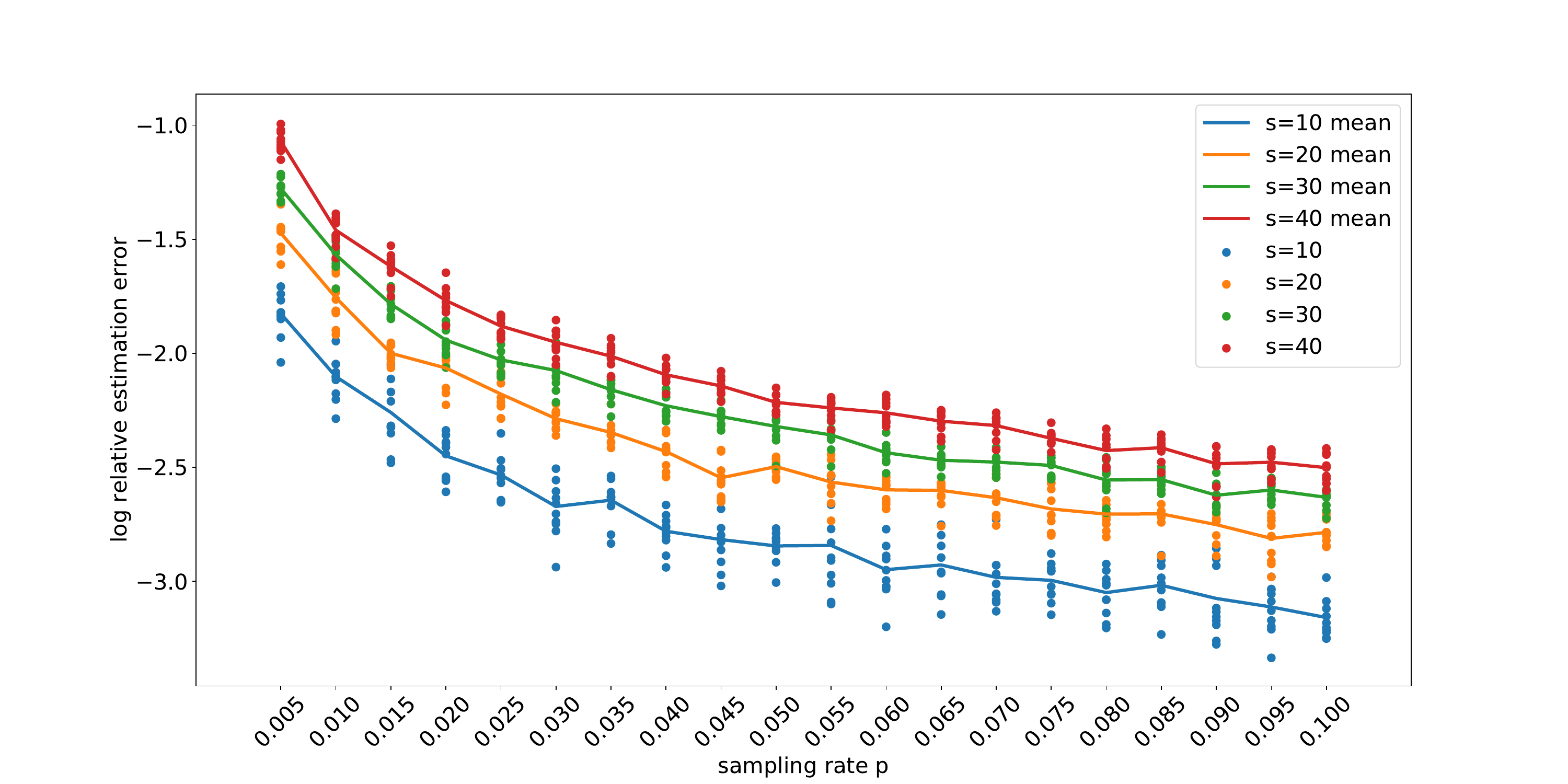}
\caption{Logarithm of relative estimation error $\log_{10} \frac{\|\widehat{\boldsymbol{M}}-\boldsymbol{M}^{\star}\|_F^2}{\|\boldsymbol{M}^{\star}\|_F^2}$ of nonconvex subspace constrained matrix completion. 
Each dot in the plot represents one trial of the numerical experiment, and the curves represent the mean of 10 independent trials for given $s$.}
\label{fig:subspace_raw}
\end{subfigure}
\begin{subfigure}[b]{\textwidth}
\centering
\includegraphics[width = 0.8\textwidth]{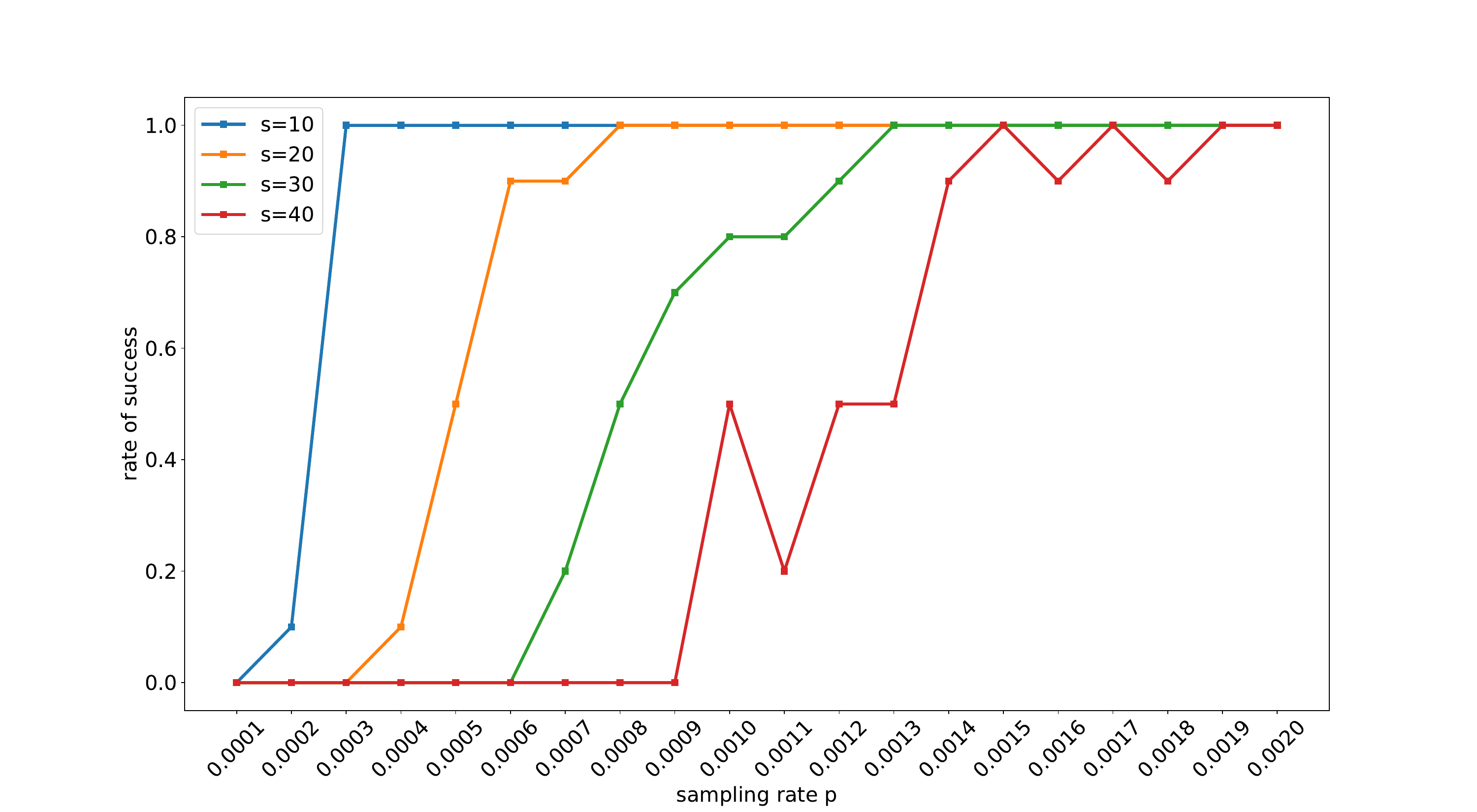}
\caption{Rates of success for noiseless nonconvex subspace constrained matrix completion. Here we set the dimension of ground truth $\boldsymbol{M}^\star \in\mathbb{R}^{n_1\times n_2}$ as $n_1 = n_2 = 500$, rank of $\boldsymbol{M}^\star$ as $r = 2$, dimension of the column/row subspace constraint as $s_1 = s_2 = s$.} 
\label{fig:subspace_phase_transition}
\end{subfigure}
\caption{Numerical results for nonconvex subspace-constrained matrix completion in both noisy and noiseless cases.}
\label{fig:subspace}
\end{figure}

In the noisy case, to implement gradient descent for solving the optimization \eqref{eq:obj_sub}, we consider different scenarios in terms of sampling rates and dimensions of subspace constraints: $p = 0.005$, $0.010$, $0.015$, ..., $0.1$ and $s = 10, 20, 30, 40$. For each fixed pair of $(p, s)$, the sampling support is generated from Model \ref{mod:sampling1}, and the subspace constraints are set by $\widetilde{\mtx{U}} = [\boldsymbol{u}_1,\ldots,\boldsymbol{u}_s]$ and $\widetilde{\mtx{V}}=[\boldsymbol{v}_1,\ldots,\boldsymbol{v}_s]$. Then, gradient descent is implemented to solve \eqref{eq:obj_sub} with the input $\mathcal{P}_{\Omega}(\boldsymbol{M}^{\star}+\boldsymbol{N})$. Logarithmic relative estimation errors $\log_{10} \frac{\|\widehat{\boldsymbol{M}}-\boldsymbol{M}^{\star}\|_F^2}{\|\boldsymbol{M}^{\star}\|_F^2}$ are averaged over $10$ independent generations of the support of observations $\Omega$ and the noise $\mtx{N}$, which are shown in Figure \ref{fig:subspace_raw}.
The results indicate that estimation errors are reduced when the the constraining subspaces have lower dimensions. This dependency illustrates our theoretical result Theorem \eqref{eq:subspace_estimate}.

Our experiments for the noiseless case are conducted in a similar manner, but here the sampling rate is chosen as $p = 0.0001, 0.0002, \ldots, 0,0020$. Rather than recording the relative errors for low-rank recovery, each experiment is viewed to be ``successful'' if and only if  $\frac{\|\widehat{\boldsymbol{M}}-\boldsymbol{M}^{\star}\|_F}{\|\boldsymbol{M}^{\star}\|_F} \leqslant 10^{-3}$, and average rates of success are plotted in Figure \ref{fig:subspace_phase_transition}. The results clearly show that with more restrictive subspace constraints, the required sample size for exact matrix completion is smaller. We admit that this phenomenon is not explained by Theorem \ref{eq:subspace_estimate}, in which the sampling size requirement is implied by the meta-theorem and hence not sufficiently adaptive to subspace-constrained models. As noted earlier, we intend to improve the sample size condition results in future work.

\subsection{Skew-Symmetric Matrix Completion}
For skew-symmetric matrix completion, we can either seek to recover the low-rank matrix by exploiting the skew-symmetric structures and solving the nonconvex optimization \eqref{eq:obj_skew}, or ignoring this structure and directly solving the standard matrix completion nonconvex objective \eqref{eq:obj_rect}. We are interested in understanding whether there is any empirical advantage to exploit the skew-symmetric structures.

\begin{figure}[htbp]
\center
\includegraphics[width=0.9\textwidth]{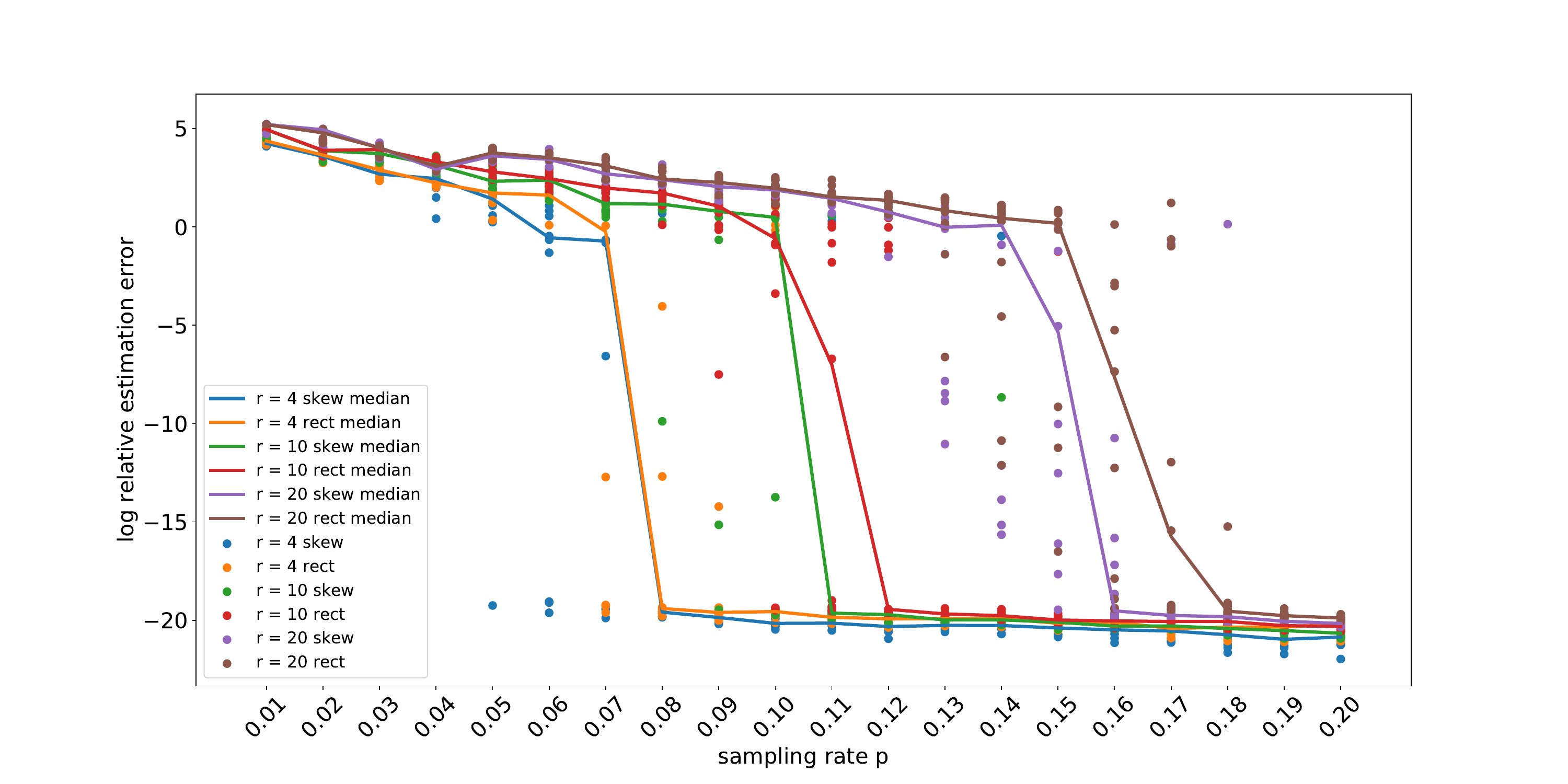}
\caption{ Logarithm of relative estimation error $\log_{10} (\|\widehat{\boldsymbol{M}}-\boldsymbol{M}^{\star}\|_F^2/\|\boldsymbol{M}^{\star}\|_F^2)$. In this plot, skew-symmetric matrix completion and regular rectangular matrix completion are used to complete noiseless skew-symmetric matrices with rank $r = 4,10,20$. Here we set the dimension of ground truth $\boldsymbol{M}^\star \in\mathbb{R}^{n\times n}$ as $n = 500$. Each dot in the plot represents one trial of the numerical experiment, and the curves represent the median of 10 independent trials.}
\label{fig:skew_transaction}
 \end{figure}

For simplicity, our experiments focused on the noiseless case. 
The matrix size is fixed as $n = 500$ and the rank is chosen to be $r = 4, 10, 20$. For each $r$, $[\boldsymbol{u}_1,\ldots, \boldsymbol{u}_{r/2},\boldsymbol{v}_1,\ldots, \boldsymbol{v}_{r/2}]$ is generated according to the Haar measure on the manifold of $500 \times r$ orthonormal basis matrices.
The low-rank matrix is then constructed as 
$\boldsymbol{M}^{\star} = \boldsymbol{u}_1\boldsymbol{v}_1^\top - \boldsymbol{v}_1\boldsymbol{u}_1^\top+\ldots +\boldsymbol{u}_{r/2}\boldsymbol{v}_{r/2}^\top - \boldsymbol{v}_{r/2}\boldsymbol{u}_{r/2}^\top$.
The sampling rate was fixed at $p = 0.01, 0.02,\ldots, 0.20$. For each fixed pair $(r, p)$, 10 independent copies of $\Omega\in[500]\times [500]$ are generated from Model \ref{mod:sampling2}. Figure \ref{fig:skew_transaction} plots the relative estimation errors as well as the corresponding medians in logarithmic scale by implementing \eqref{eq:obj_skew} and \eqref{eq:obj_rect} with gradient descent, respectively. The comparison
indicates that when the rank $r$ is not too small, exploiting skew-symmetric structures in nonconvex optimization is helpful in reducing the required sample size for low-rank recovery.

\section{Proofs}
\label{sec:proof_main}

In this section, 
we give a proof of Theorem \ref{thm:main}, followed by proofs of the two corollaries on subspace-constrained and skew-symmetric matrix completion by verifying the assumption of correlated parametric factorization in respective settings. 

In the proof of Theorem \ref{thm:main}, for any local minimum $\hat{\boldsymbol{\xi}}$ of \prettyref{eq:obj_para}, we aim to control the estimation error $\|\mtx{M}^\star  - \boldsymbol{X}(\hat{\boldsymbol{\xi}}) \boldsymbol{Y}(\hat{\boldsymbol{\xi}})^\top\|_F^2$. 
Given $\tilde{f}(\vct{\theta})$ is smooth, it is natural to study the first and second order optimality conditions, i.e., $\nabla \tilde{f}(\vct{\theta}) = \vct{0}$ and $\nabla^2 \tilde{f}(\vct{\theta}) \succeq \mtx{0}_{d \times d}$. The strategy on employing these two conditions in the control of estimation error is key. 
To this end, we start with introducing an auxiliary function that connects the optimality conditions and the estimation error.

\subsection{An Auxiliary Function}

We first recall a key component in geometric analysis for nonconvex matrix completion: the auxiliary function associated to the gradient and Hessian of $f(\boldsymbol{X},\boldsymbol{Y})$ defined in \eqref{eq:obj_rect}, which has been used in \citet{ge2017no, jin2017escape, chen2017memory}. For any $\boldsymbol{X} \in \mathbb{R}^{n_1 \times r}$, $\boldsymbol{Y} \in \mathbb{R}^{n_2 \times r}$, $\boldsymbol{D}_{\boldsymbol{X}} \in \mathbb{R}^{n_1 \times r}$ and $\boldsymbol{D}_{\boldsymbol{Y}} \in \mathbb{R}^{n_2 \times r}$, the auxiliary function associated to $f(\mtx{X}, \mtx{Y})$ is defined as
\begin{equation}
\label{eq:auxiliary_rect}
\begin{split}
K_f(\boldsymbol{X}, \boldsymbol{Y}; \boldsymbol{D}_{\boldsymbol{X}}, \boldsymbol{D}_{\boldsymbol{Y}}) 
	\coloneqq   \operatorname{vec}(\boldsymbol{D}_{\boldsymbol{X}}, \boldsymbol{D}_{\boldsymbol{Y}})^\top \nabla^2 f(\boldsymbol{X}, \boldsymbol{Y})\operatorname{vec}(\boldsymbol{D}_{\boldsymbol{X}}, \boldsymbol{D}_{\boldsymbol{Y}})~&\\
- 4\left\langle\nabla f(\boldsymbol{X}, \boldsymbol{Y}), \operatorname{vec}(\boldsymbol{D}_{\boldsymbol{X}}, \boldsymbol{D}_{\boldsymbol{Y}})  \right\rangle.&
\end{split}
\end{equation}
Because any local minimum $(\widehat{\boldsymbol{X}}, \widehat{\boldsymbol{Y}})$ of $f(\boldsymbol{X}, \boldsymbol{Y})$ satisfies the first and second order conditions, i.e., $\nabla f(\widehat{\boldsymbol{X}}, \widehat{\boldsymbol{Y}}) = \vct{0}$ and $\nabla^2 f(\widehat{\boldsymbol{X}}, \widehat{\boldsymbol{Y}}) \succeq \mtx{0}$,
we have
\[
K_f(\widehat{\boldsymbol{X}}, \widehat{\boldsymbol{Y}}; \boldsymbol{D}_{\boldsymbol{X}}, \boldsymbol{D}_{\boldsymbol{Y}}) \geqslant  0, \quad \forall \boldsymbol{D}_{\boldsymbol{X}} \in \mathbb{R}^{n_1 \times r}, \boldsymbol{D}_{\boldsymbol{Y}} \in \mathbb{R}^{n_2 \times r}.
\]
The ratio of coefficients $1:(-4)$ for the two terms on the right hand side of \eqref{eq:auxiliary_rect} turns out to be essential for the nonconvex geometric analysis; see \citet{jin2017escape, ge2017no, chen2017memory}.

Following this idea, we define an analogous auxiliary function associated to $\tilde{f}(\boldsymbol{\theta})$ in general parametric matrix completion \eqref{eq:obj_para}. For any $\boldsymbol{\theta}, \boldsymbol{\delta}_{\boldsymbol{\theta}} \in \mathbb{R}^d$, define the auxiliary function associated with $\tilde{f}$ as
\begin{equation}
\label{eq:auxiliary_para}
K_{\tilde{f}}(\boldsymbol{\theta}; \boldsymbol{\delta}_{\boldsymbol{\theta}})\coloneqq \boldsymbol{\delta}_{\boldsymbol{\theta}}^\top \nabla^2 \tilde{f}(\boldsymbol{\theta}) \boldsymbol{\delta}_{\boldsymbol{\theta}} - 4\boldsymbol{\delta}_{\boldsymbol{\theta}}^\top \nabla \tilde{f}(\boldsymbol{\theta}).
\end{equation}
Again, any minimum $\hat{\boldsymbol{\xi}}$ of $\tilde{f}$ satisfies $\nabla\tilde{f}(\hat{\boldsymbol{\xi}}) = \boldsymbol{0}$ and $\nabla^2 \tilde{f}(\hat{\boldsymbol{\xi}}) \succeq \boldsymbol{0}$, which gives
\[
K_{\tilde{f}}(\hat{\boldsymbol{\xi}}; \boldsymbol{\delta}_{\boldsymbol{\theta}}) \geqslant  0, \quad \forall \boldsymbol{\delta}_{\boldsymbol{\theta}} \in \mathbb{R}^d.
\]
Control of estimation error $\|\mtx{M}^\star  - \boldsymbol{X}(\hat{\boldsymbol{\xi}}) \boldsymbol{Y}(\hat{\boldsymbol{\xi}})^\top\|_F^2$ relies on controlling $K_{\tilde{f}}(\boldsymbol{\theta}; \boldsymbol{\delta}_{\boldsymbol{\theta}})$, which in turn depends crucially on a careful choice of $\boldsymbol{\delta}_{\boldsymbol{\theta}}$. To this end, Assumption \ref{ass:ideal_factorization} plays the key role and serves as the guiding principle. 

We first give a lemma that characterizes the relationship between $K_{\tilde{f}}(\boldsymbol{\theta}; \boldsymbol{\delta}_{\boldsymbol{\theta}})$ and $K_f(\boldsymbol{X}, \boldsymbol{Y}; \boldsymbol{D}_{\boldsymbol{X}}, \boldsymbol{D}_{\boldsymbol{Y}})$:
\begin{lemma}
\label{lem:key_identity}
For any $\boldsymbol{\theta}, \boldsymbol{\delta}_{\boldsymbol{\theta}} \in \mathbb{R}^p$, it holds that 
\begin{equation}
\label{eq:aux_relation}
K_{\tilde{f}}(\boldsymbol{\theta}; \boldsymbol{\delta}_{\boldsymbol{\theta}})= K_f(\boldsymbol{X}(\boldsymbol{\theta}), \boldsymbol{Y}(\boldsymbol{\theta}); \boldsymbol{X}(\boldsymbol{\delta}_{\boldsymbol{\theta}}), \boldsymbol{Y}(\boldsymbol{\delta}_{\boldsymbol{\theta}})).
\end{equation} 
\end{lemma}

\begin{proof}
Assumption \ref{ass:homogeneous} implies
\begin{equation}
\label{eq:starting_equality}
\tilde{f}(\boldsymbol{\theta} + \boldsymbol{\delta}_{\boldsymbol{\theta}}) = f(\boldsymbol{X}(\boldsymbol{\theta} + \boldsymbol{\delta}_{\boldsymbol{\theta}}), \boldsymbol{Y}(\boldsymbol{\theta} + \boldsymbol{\delta}_{\boldsymbol{\theta}})) = f(\boldsymbol{X}(\boldsymbol{\theta}) + \boldsymbol{X}(\boldsymbol{\delta}_{\boldsymbol{\theta}}), \boldsymbol{Y}(\boldsymbol{\theta}) + \boldsymbol{Y}(\boldsymbol{\delta}_{\boldsymbol{\theta}})).
\end{equation}
Due to the linearity and homogeneity of $\boldsymbol{X}(\boldsymbol{\theta})$ and $\boldsymbol{Y}(\boldsymbol{\theta})$, considering the Taylor expansions of both sides in \eqref{eq:starting_equality} at $\boldsymbol{\theta}$, we get
\begin{equation}
\label{eq:gradient_eq_fact}
\boldsymbol{\delta}_{\boldsymbol{\theta}}^\top \nabla \tilde{f}(\boldsymbol{\theta}) = \text{vec}(\boldsymbol{X}(\boldsymbol{\delta}_{\boldsymbol{\theta}}), \boldsymbol{Y}(\boldsymbol{\delta}_{\boldsymbol{\theta}}))^\top \nabla f(\boldsymbol{X}(\boldsymbol{\theta}), \boldsymbol{Y}(\boldsymbol{\theta}))
\end{equation}
and
\begin{equation}
\label{eq:hessian_eq_fact}
\boldsymbol{\delta}_{\boldsymbol{\theta}}^\top \nabla^2 \tilde{f}(\boldsymbol{\theta}) \boldsymbol{\delta}_{\boldsymbol{\theta}} = \text{vec}(\boldsymbol{X}(\boldsymbol{\delta}_{\boldsymbol{\theta}}), \boldsymbol{Y}(\boldsymbol{\delta}_{\boldsymbol{\theta}}))^\top \nabla^2 f(\boldsymbol{X}(\boldsymbol{\theta}), \boldsymbol{Y}(\boldsymbol{\theta})) \text{vec}(\boldsymbol{X}(\boldsymbol{\delta}_{\boldsymbol{\theta}}), \boldsymbol{Y}(\boldsymbol{\delta}_{\boldsymbol{\theta}}))
\end{equation}
By combining \eqref{eq:auxiliary_rect}, \eqref{eq:auxiliary_para}, \eqref{eq:gradient_eq_fact} and \eqref{eq:hessian_eq_fact}, the equality \eqref{eq:aux_relation} is obtained.
\end{proof}

Our choice of $\boldsymbol{\delta}_{\boldsymbol{\theta}}$ is inspired by the choices of $\mtx{D}_{\mtx{X}}$ and $\mtx{D}_{\mtx{Y}}$ for $K_f$ in \citet{jin2017escape, ge2017no, chen2017memory}. For any $\vct{\theta} \in \mathbb{R}^d$, let $\boldsymbol{\xi}\in\mathbb{R}^d$ be selected (not necessarily uniquely) such that Assumption \ref{ass:ideal_factorization} is satisfied. Then let 
\[
\boldsymbol{\delta}_{\boldsymbol{\theta}} \coloneqq \boldsymbol{\theta} - \boldsymbol{\xi}.
\]

In the rest of this subsection, we will explain why our choice of $\boldsymbol{\delta}_{\boldsymbol{\theta}}$ is helpful in controlling the auxiliary function. For simplicity, we introduce the following abbreviations and notations:  
\begin{equation}
\label{eq:simplified_notation}
\begin{aligned}
\boldsymbol{X} = \boldsymbol{X}(\boldsymbol{\theta})    \in \mathbb{R}^{n_1 \times r},\quad
&\boldsymbol{U}=\boldsymbol{X}(\boldsymbol{\xi}) \in \mathbb{R}^{n_1 \times r},\quad
\boldsymbol{\Delta}_{\boldsymbol{X}} = \boldsymbol{X}(\boldsymbol{\delta}_{\boldsymbol{\theta}}) = \boldsymbol{X} - \boldsymbol{U} \in \mathbb{R}^{n_1 \times r},
\\
\boldsymbol{Y} = \boldsymbol{Y}(\boldsymbol{\theta})  \in \mathbb{R}^{n_2 \times r},\quad
&\boldsymbol{V} = \boldsymbol{Y}(\boldsymbol{\xi})  \in \mathbb{R}^{n_2 \times r},\quad
\boldsymbol{\Delta}_{\boldsymbol{Y}}= \boldsymbol{Y}(\boldsymbol{\delta}_{\boldsymbol{\theta}}) = \boldsymbol{Y} - \boldsymbol{V} \in \mathbb{R}^{n_2 \times r}. 
\end{aligned}
\end{equation}
Throughout our proof, $\boldsymbol{X} ,\boldsymbol{U},\boldsymbol{\Delta}_{\boldsymbol{X}}, \boldsymbol{Y} ,\boldsymbol{V}$, and $\boldsymbol{\Delta}_{\boldsymbol{Y}}$ will refer to the matrices defined in \eqref{eq:simplified_notation} if not specified otherwise. Now Lemma \ref{lem:key_identity} gives
\[
K_{\tilde{f}}(\boldsymbol{\theta}; \boldsymbol{\delta}_{\boldsymbol{\theta}}) = K_f(\mtx{X}, \mtx{Y}; \mtx{\Delta}_{\mtx{X}}, \mtx{\Delta}_{\mtx{Y}}).
\]
Note that Assumption \ref{ass:ideal_factorization} implies
\begin{equation}
\label{eq:rotation_rect}
\boldsymbol{M}^{\star}  = \mtx{U}\mtx{V}^\top,~\mtx{U}^\top \mtx{U} = \mtx{V}^\top \mtx{V},~\text{and}~\mtx{X}^\top \mtx{U} + \mtx{Y}^\top \mtx{V} \succeq \boldsymbol{0}.
\end{equation}
The following lemma gives an explicit upper bound for $K_{\tilde{f}}(\boldsymbol{\theta}; \boldsymbol{\delta}_{\boldsymbol{\theta}})$.
\begin{lemma}
\label{lem:ge_aux_formula}
For any $\boldsymbol{\theta} \in \mathbb{R}^d$, with $\boldsymbol{\delta}_{\boldsymbol{\theta}}= \boldsymbol{\theta} - \boldsymbol{\xi}$ where $\boldsymbol{\xi}$ satisfies the conditions in \eqref{eq:balanced parameterized factorization}, and $\boldsymbol{X},\boldsymbol{U},\boldsymbol{\Delta}_{\boldsymbol{X}},\boldsymbol{Y},\boldsymbol{V},\boldsymbol{\Delta}_{\boldsymbol{Y}}$ defined as in \eqref{eq:simplified_notation}, denote 
\begin{equation*}
\boldsymbol{Z} = \begin{bmatrix} \boldsymbol{X} \\ \boldsymbol{Y} \end{bmatrix}, ~\boldsymbol{W} = \begin{bmatrix}\boldsymbol{U} \\ \boldsymbol{V} \end{bmatrix}, \text{~and~} \boldsymbol{\Delta}_{\boldsymbol{Z}} = \begin{bmatrix}\boldsymbol{\Delta}_{\boldsymbol{X}} \\ \boldsymbol{\Delta}_{\boldsymbol{Y}} \end{bmatrix} = \mtx{Z} - \mtx{W}. 
\end{equation*}
Then the auxiliary function $K_{\tilde{f}}(\boldsymbol{\theta}; \boldsymbol{\delta}_{\boldsymbol{\theta}})$ defined in \eqref{eq:auxiliary_para} can be upper bounded as following:
\begin{equation}
\label{eq:K_noisy_decomp_ks}
\begin{split}
      K_{\tilde{f}}(\boldsymbol{\theta}; \boldsymbol{\delta}_{\boldsymbol{\theta}}) \leqslant &  K_1 + K_2+K_3 +K_4,\\
      \end{split} 
    \end{equation}
    where 
    \begin{equation}
    \label{eq:K_noisy_decomp}
     \begin{split}
    	K_1 \coloneqq &\frac{1}{4} \left( \left\| \boldsymbol{\Delta}_{\boldsymbol{Z}}\boldsymbol{\Delta}_{\boldsymbol{Z}}^\top  \right\|_F^2 - 3\left\| \boldsymbol{Z}\boldsymbol{Z}^\top -\boldsymbol{W}\boldsymbol{W}^\top \right\|_F^2 \right),\\
	K_2 \coloneqq& \left(\frac{1}{p}\left\| \mathcal{P}_{\Omega}\left(\boldsymbol{\Delta}_{\boldsymbol{X}}\boldsymbol{\Delta}_{\boldsymbol{Y}}^\top \right) \right\|_F^2 - \|\boldsymbol{\Delta}_{\boldsymbol{X}}\boldsymbol{\Delta}_{\boldsymbol{Y}}^\top \|_F^2 \right) \\
	&- \left( \frac{3}{p}\left\| \mathcal{P}_{\Omega}\left( \boldsymbol{X}\boldsymbol{Y}^\top - \boldsymbol{U}\boldsymbol{V}^\top  \right) \right\|_F^2 - 3 \|\boldsymbol{X}\boldsymbol{Y}^\top - \boldsymbol{U}\boldsymbol{V}^\top\|_F^2 \right),\\
	K_3 \coloneqq& \lambda \left[\operatorname{vec}(\boldsymbol{\Delta}_{\boldsymbol{X}})^\top \nabla^2 G_{\alpha}(\boldsymbol{X}) \operatorname{vec}(\boldsymbol{\Delta}_{\boldsymbol{X}})- 4\left\langle \nabla G_{\alpha}(\boldsymbol{X}),\boldsymbol{\Delta}_{\boldsymbol{X}} \right\rangle  \right]\\
	&+\lambda \left[\operatorname{vec}(\boldsymbol{\Delta}_{\boldsymbol{Y}})^\top \nabla^2 G_{\alpha}(\boldsymbol{Y}) \operatorname{vec}(\boldsymbol{\Delta}_{\boldsymbol{Y}})- 4\left\langle \nabla G_{\alpha}(\boldsymbol{Y}),\boldsymbol{\Delta}_{\boldsymbol{Y}} \right\rangle  \right],\\
	K_4 \coloneqq  &\frac{6}{p}\left\langle \boldsymbol{\Delta}_{\boldsymbol{X}}\boldsymbol{\Delta}_{\boldsymbol{Y}}^\top, \mathcal{P}_{\Omega}(\boldsymbol{N}) \right\rangle + \frac{4}{p} \left\langle \boldsymbol{U}\boldsymbol{\Delta}_{\boldsymbol{Y}}^\top + \boldsymbol{\Delta}_{\boldsymbol{X}}\boldsymbol{V}^\top, \mathcal{P}_{\Omega}(\boldsymbol{N}) \right\rangle.
	\end{split}
    \end{equation}
\end{lemma}
The proof of Lemma \ref{lem:ge_aux_formula} follows closely that of \citet[Lemma 16]{ge2017no}. We give a complete proof in Appendix \ref{sec:proof_ge_aux_formula} for being self-contained.

\subsection{Proof of Theorem \ref{thm:main}}
\label{sec:comp_proof_main}

Given the decomposition in Lemma \ref{lem:ge_aux_formula}, the proof of Theorem \ref{thm:main} depends on the following two auxiliary lemmas. Note that all results are established on the event $E$ defined in Lemma \ref{lemma:eq_of_rate}.

\begin{lemma}
\label{lem:K2K3}
Let the quantities $K_2$ and $K_3$ be defined as in Lemma \ref{lem:ge_aux_formula}. Then the sum $K_2+K_3$ is bounded as
\begin{equation*}
K_2 + K_3 \leqslant   0.001 \left(\| \boldsymbol{\Delta}_{\boldsymbol{X}}\boldsymbol{V}^\top \|_F^2  + \|\boldsymbol{U}\boldsymbol{\Delta}_{\boldsymbol{Y}}^\top \|_F^2\right).
\end{equation*}
\end{lemma}
\begin{proof}
See Section \ref{sec:pf-lem-K2K3}.
\end{proof}

\begin{lemma}
\label{lem:K4}
Let the quantity $K_4$ be defined as in Lemma \ref{lem:ge_aux_formula}. 
If $K_1 + K_2 + K_3 + K_4 \geqslant 0$, the following upper bounds hold
\begin{enumerate}
	\item $\|\boldsymbol{\Delta}_{\boldsymbol{Z}} \boldsymbol{\Delta}_{\boldsymbol{Z}}^\top\|_F = \|\boldsymbol{\Delta}_{\boldsymbol{Z}}^\top \boldsymbol{\Delta}_{\boldsymbol{Z}}\|_F\leqslant   3\sqrt{|K_4|}, \quad \|\boldsymbol{\Delta}_{\boldsymbol{Z}}^\top \boldsymbol{W}\|_F \leqslant   2\sqrt{|K_4|}$.
	
	\item 
	$\left\| \boldsymbol{\Delta}_{\boldsymbol{X}}\boldsymbol{V}^\top\right\|_F^2  + \left\|\boldsymbol{U}\boldsymbol{\Delta}_{\boldsymbol{Y}}^\top \right\|_F^2 \leqslant 13|K_4|$.
\end{enumerate}
\end{lemma}
\begin{proof}
See Section \ref{sec:pf-lem-K4}.
\end{proof}

By the definition of $\boldsymbol{X}, \boldsymbol{Y}, \boldsymbol{U}$ and $\boldsymbol{V}$ in \eqref{eq:simplified_notation}, we note that the difference between any local minimizer $\widehat{\boldsymbol{M}}$ and the matrix $\boldsymbol{M}^\star$ can be written as 
\begin{equation}
	\label{eq:pf-main-1}
\|\widehat{\boldsymbol{M}} - \boldsymbol{M}^\star  \|_F^2 = 
\| \boldsymbol{X}\boldsymbol{Y}^\top - \boldsymbol{U}\boldsymbol{V}^\top \|_F^2
\end{equation}
Next, by the definition of $\boldsymbol{Z}$ and $\boldsymbol{W}$ in Lemma \ref{lem:ge_aux_formula}, we expand $\boldsymbol{Z}\boldsymbol{Z}^\top - \bs{W}\bs{W}^\top$ to obtain the inequality
\begin{equation}
		\label{eq:pf-main-2}
\| \boldsymbol{X}\boldsymbol{Y}^\top - \boldsymbol{U}\boldsymbol{V}^\top \|_F^2
\leqslant 
\| \boldsymbol{Z}\boldsymbol{Z}^\top - \bs{W}\bs{W}^\top \|_F^2.
\end{equation}
Since $\widehat{\bs{M}} = \bs{X}(\hat{\bs{\xi}}) \bs{Y}(\hat{\bs{\xi}})^\top$ is a local minimizer,
$\hat{\boldsymbol{\xi}}$ is a local minimum of $\tilde{f}$ and so it holds that $K_{\tilde{f}}(\hat{\boldsymbol{\xi}}; \boldsymbol{\delta}_{\hat{\boldsymbol{\xi}}}) \geqslant  0$. By Lemma \ref{lem:ge_aux_formula}, it holds that
\begin{equation}
\label{eq:k1234geq0}
K_1 + K_2 + K_3 + K_4 \geqslant 0.
\end{equation}
By Lemma \ref{lem:K2K3} and the second claim of Lemma \ref{lem:K4}, $K_2 + K_3 + K_4 \leqslant 2|K_4|$. 
Together with the definition of $K_1$ in \eqref{eq:K_noisy_decomp}, \eqref{eq:k1234geq0} and the first claim of Lemma \ref{lem:K4}, this leads to
\begin{equation}
	\label{eq:pf-main-3}
\frac{3}{4} \|\boldsymbol{Z}\boldsymbol{Z}^\top - \boldsymbol{W}\boldsymbol{W}^\top\|_F^2 \leqslant \frac{1}{4} \|\boldsymbol{\Delta}_{\boldsymbol{Z}} \boldsymbol{\Delta}_{\boldsymbol{Z}}^\top\|_F^2 + K_2 + K_3 + K_4 \leqslant \frac{17}{4}|K_4|.
\end{equation}
If we can show that 
\begin{equation}
\label{eq:K4abs}
  |K_4| \leqslant \frac{10^4r}{p^2 }\psi^2,
\end{equation}
then the desired conclusion of Theorem \ref{thm:main} is a direct consequence of \eqref{eq:pf-main-1}, \eqref{eq:pf-main-2}, \eqref{eq:pf-main-3}, and \eqref{eq:K4abs}.
It remains only to prove \eqref{eq:K4abs}, to which we now turn.

\paragraph{Proof of \eqref{eq:K4abs}}
Recall the fact that $\boldsymbol{\Delta}_{\boldsymbol{Z}} = \begin{bmatrix}\boldsymbol{\Delta}_{\boldsymbol{X}}^\top & \boldsymbol{\Delta}_{\boldsymbol{Y}}^\top \end{bmatrix}^\top$, then $\|\boldsymbol{\Delta}_{\boldsymbol{X}}\boldsymbol{\Delta}_{\boldsymbol{Y}}^\top\|_F\leqslant \|\boldsymbol{\Delta}_{\boldsymbol{Z}} \boldsymbol{\Delta}_{\boldsymbol{Z}}^\top\|_F$. Therefore, by \prettyref{lem:K4}, and the definition of $K_4$ in \eqref{eq:K_noisy_decomp}, we have 
\begin{equation}
\label{eq:k4_upper}
    \begin{split}
      |K_4| \leqslant & \frac{6}{p} \|\boldsymbol{\Delta}_{\boldsymbol{X}}\boldsymbol{\Delta}_{\boldsymbol{Y}}^\top\|_F \|\boldsymbol{P}_{\boldsymbol{\Delta}_{\boldsymbol{X}}}\mathcal{P}_{\Omega}(\boldsymbol{N})\boldsymbol{P}_{\boldsymbol{\Delta}_{\boldsymbol{Y}}}\|_F +\frac{4}{p}\|\boldsymbol{U}\boldsymbol{\Delta}_{\boldsymbol{Y}}^\top \|_F\|\boldsymbol{P}_{\boldsymbol{U}}\mathcal{P}_{\Omega}(\boldsymbol{N})\boldsymbol{P}_{\boldsymbol{\Delta}_{\boldsymbol{Y}}}\|_F  \\
      &+ \frac{4}{p}\|\boldsymbol{\Delta}_{\boldsymbol{X}}\boldsymbol{V}^\top\|_F \|\boldsymbol{P}_{\boldsymbol{\Delta}_{\boldsymbol{X}}}\mathcal{P}_{\Omega}(\boldsymbol{N})\boldsymbol{P}_{\boldsymbol{V}}\|_F\\
      \leqslant & \frac{100\sqrt{|K_4|}}{p}\max\{ \|\boldsymbol{P}_{\boldsymbol{\Delta}_{\boldsymbol{X}}}\mathcal{P}_{\Omega}(\boldsymbol{N})\boldsymbol{P}_{\boldsymbol{\Delta}_{\boldsymbol{Y}}}\|_F, \|\boldsymbol{P}_{\boldsymbol{U}}\mathcal{P}_{\Omega}(\boldsymbol{N})\boldsymbol{P}_{\boldsymbol{\Delta}_{\boldsymbol{Y}}}\|_F, \|\boldsymbol{P}_{\boldsymbol{\Delta}_{\boldsymbol{X}}}\mathcal{P}_{\Omega}(\boldsymbol{N})\boldsymbol{P}_{\boldsymbol{V}}\|_F\}.
    \end{split}
\end{equation}
Due to the fact that $\boldsymbol{U},\boldsymbol{\Delta}_{\boldsymbol{X}}\in\mathbb{R}^{n_1\times r}$ and $\boldsymbol{V},\boldsymbol{\Delta}_{\boldsymbol{Y}}\in\mathbb{R}^{n_2\times r}$, $\boldsymbol{P}_{\boldsymbol{\Delta}_{\boldsymbol{X}}}\mathcal{P}_{\Omega}(\boldsymbol{N})\boldsymbol{P}_{\boldsymbol{\Delta}_{\boldsymbol{Y}}}$, $\boldsymbol{P}_{\boldsymbol{U}}\mathcal{P}_{\Omega}(\boldsymbol{N})\boldsymbol{P}_{\boldsymbol{\Delta}_{\boldsymbol{Y}}}$ and $\boldsymbol{P}_{\boldsymbol{\Delta}_{\boldsymbol{X}}}\mathcal{P}_{\Omega}(\boldsymbol{N})\boldsymbol{P}_{\boldsymbol{V}}$ are matrices with rank at most $r$. Therefore,
\[
  \begin{split}
    & \max\{ \|\boldsymbol{P}_{\boldsymbol{\Delta}_{\boldsymbol{X}}}\mathcal{P}_{\Omega}(\boldsymbol{N})\boldsymbol{P}_{\boldsymbol{\Delta}_{\boldsymbol{Y}}}\|_F, \|\boldsymbol{P}_{\boldsymbol{U}}\mathcal{P}_{\Omega}(\boldsymbol{N})\boldsymbol{P}_{\boldsymbol{\Delta}_{\boldsymbol{Y}}}\|_F, \|\boldsymbol{P}_{\boldsymbol{\Delta}_{\boldsymbol{X}}}\mathcal{P}_{\Omega}(\boldsymbol{N})\boldsymbol{P}_{\boldsymbol{V}}\|_F\}\\
    \leqslant& \sqrt{r} \max\{ \|\boldsymbol{P}_{\boldsymbol{\Delta}_{\boldsymbol{X}}}\mathcal{P}_{\Omega}(\boldsymbol{N})\boldsymbol{P}_{\boldsymbol{\Delta}_{\boldsymbol{Y}}}\|, \|\boldsymbol{P}_{\boldsymbol{U}}\mathcal{P}_{\Omega}(\boldsymbol{N})\boldsymbol{P}_{\boldsymbol{\Delta}_{\boldsymbol{Y}}}\|, \|\boldsymbol{P}_{\boldsymbol{\Delta}_{\boldsymbol{X}}}\mathcal{P}_{\Omega}(\boldsymbol{N})\boldsymbol{P}_{\boldsymbol{V}}\|\}\\
     = &\sqrt{r}\psi.
  \end{split} 
\]
Where the last line follows from \eqref{eq:simplified_notation} and \eqref{eq:psi_def}. Therefore, \eqref{eq:k4_upper} gives
\[
  |K_4|  \leqslant \frac{100\sqrt{|K_4|r}}{p}\psi.
\]
Rearranging terms in the last display yields \eqref{eq:K4abs}.

%
%
%

%
%
%

\subsection{Analysis for Subspace-Constrained Matrix Completion}
\label{sec:subspace_proof}
We give a proof of Lemma \ref{lem:rotation_positive_subspace} followed by that of Corollary \ref{cor:subspace}.

\subsubsection{Proof of Lemma \ref{lem:rotation_positive_subspace}}
\label{sec:rotation_positive_subspace}
\begin{proof}
In order to show the parameterization \eqref{eq:para_subspace} satisfies Assumption \ref{ass:ideal_factorization}, we want to show that for any $\boldsymbol{\theta} = \operatorname{vec}(\boldsymbol{\Theta}_{A}, \boldsymbol{\Theta}_{B}) \in \mathbb{R}^{r(s_1+s_2)}$, there exits a $\boldsymbol{\xi} = \operatorname{vec}(\boldsymbol{\Xi}_{A}, \boldsymbol{\Xi}_{B}) \in \mathbb{R}^{r(s_1+s_2)}$ that satisfies \eqref{eq:balanced parameterized factorization}. Let $\boldsymbol{S} \coloneqq  \widetilde{\boldsymbol{U}}^\top \boldsymbol{M}^\star \widetilde{\boldsymbol{V}} \in \mathbb{R}^{s_1 \times s_2}$. Recall that $\widetilde{\boldsymbol{U}}$ consists of an orthonormal basis of the column space constraint, and $\widetilde{\boldsymbol{V}}$ consists of an orthonormal basis of the column row constraint of $\boldsymbol{M}^\star$. Therefore, $\boldsymbol{P}_{\widetilde{\boldsymbol{U}}} = \widetilde{\boldsymbol{U}}\widetilde{\boldsymbol{U}}^\top$,  $\boldsymbol{P}_{\widetilde{\boldsymbol{V}}} = \widetilde{\boldsymbol{V}}\widetilde{\boldsymbol{V}}^\top$ and $\mtx{M}^\star $ can be represented as $\mtx{M}^\star  = \widetilde{\boldsymbol{U}} \boldsymbol{S} \widetilde{\boldsymbol{V}}^\top$. Since $\boldsymbol{M}^\star$ is of rank $r$, by the orthogonality of $\widetilde{\boldsymbol{U}}$ and $\widetilde{\boldsymbol{V}}$, $\rank(\mtx{S})=r$. Let the reduced SVD of $\mtx{S}$ be 
\begin{equation}\label{eq:H_SVD}
  \mtx{S} = \mtx{S}_L \mtx{\Lambda} \mtx{S}_R^\top,
\end{equation}
where $\mtx{S}_L \in \mathbb{R}^{s_1 \times r}$, $\mtx{S}_R \in \mathbb{R}^{s_2 \times r}$, $\mtx{S}_L^\top \mtx{S}_L = \boldsymbol{I}_{s_1}, \mtx{S}_R^\top \mtx{S}_R = \boldsymbol{I}_{s_2}$ and $\mtx{\Lambda} = \diag(\sigma_1, \ldots, \sigma_r)$ is a $r\times r$ diagonal matrix with $\sigma_1\geqslant \sigma_2\geqslant  \ldots \geqslant \sigma_r$. Moreover, by letting $\boldsymbol{U}^{\star}  \coloneqq \widetilde{\mtx{U}}\mtx{S}_L \in\mathbb{R}^{n_1\times r}, \boldsymbol{V}^{\star} \coloneqq \widetilde{\mtx{V}}\mtx{S}_R   \in\mathbb{R}^{n_2\times r}$, we can verify that $\mtx{M}^\star  = \mtx{U}^{\star}  \mtx{\Lambda}  {\mtx{V}^{\star} }^\top$ is a reduced SVD of $\mtx{M}^\star $.

Define
\[
\boldsymbol{\Xi}_A^{\star} \coloneqq \mtx{S}_L \mtx{\Lambda}^{1/2}  \in \mathbb{R}^{s_1 \times r},~\quad \boldsymbol{\Xi}_B^{\star} \coloneqq \mtx{S}_R \mtx{\Lambda}^{1/2} \in \mathbb{R}^{s_2 \times r}.
\]
For any $\boldsymbol{\Theta}_{A} \in \mathbb{R}^{s_1 \times r}$ and $, \boldsymbol{\Theta}_{B} \in \mathbb{R}^{s_2 \times r}$, by considering the SVD of $(\boldsymbol{\Theta}_{A}^\top \boldsymbol{\Xi}_A^{\star} + \boldsymbol{\Theta}_{B}^\top \boldsymbol{\Xi}_B^{\star})$, we know there exits an $r \times r$ orthogonal matrix $\mtx{T} \in O(r)$ \citep[Lemma 1]{chen2015fast}, such that 
\[
(\boldsymbol{\Theta}_{A}^\top \boldsymbol{\Xi}_A^{\star} + \boldsymbol{\Theta}_{B}^\top \boldsymbol{\Xi}_B^{\star})\boldsymbol{T} \succeq \mtx{0}.
\]

Let $\boldsymbol{\xi} = \operatorname{vec}(\boldsymbol{\Xi}_A^{\star}\boldsymbol{T}, \boldsymbol{\Xi}_B^{\star}\boldsymbol{T})$, then
\[
    \boldsymbol{X}(\boldsymbol{\xi}) =  \widetilde{\boldsymbol{U}}\boldsymbol{\Xi}_A^{\star}\boldsymbol{T} \in \mathbb{R}^{n_1 \times r}\quad \text{~and~} \quad \boldsymbol{Y}(\boldsymbol{\xi}) = \widetilde{\boldsymbol{V}}\boldsymbol{\Xi}_B^{\star}\boldsymbol{T} \in \mathbb{R}^{n_2 \times r}.
\] 
Keeping in mind that both $\widetilde{\mtx{U}}$ and $\widetilde{\mtx{V}}$ are orthonormal basis matrices, the conditions in \eqref{eq:balanced parameterized factorization} can be verified one by one:
\[
 \boldsymbol{X}(\boldsymbol{\xi})\boldsymbol{Y}(\boldsymbol{\xi})^\top = \widetilde{\boldsymbol{U}}\boldsymbol{\Xi}_A^{\star}\boldsymbol{T}(\widetilde{\boldsymbol{V}} \boldsymbol{\Xi}_B^{\star}\boldsymbol{T})^\top = \widetilde{\boldsymbol{U}} \boldsymbol{\Xi}_A^{\star} {\boldsymbol{\Xi}_B^{\star}}^\top \widetilde{\boldsymbol{V}}^\top= \widetilde{\boldsymbol{U}} \mtx{S}_L \mtx{\Lambda} \mtx{S}_R^\top \widetilde{\boldsymbol{V}}^\top =\boldsymbol{M}^{\star}.
 \]
 The last equality is by \eqref{eq:H_SVD}.
 \begin{align*}
 \boldsymbol{X}(\boldsymbol{\xi})^\top \boldsymbol{X}(\boldsymbol{\xi}) &= (\widetilde{\boldsymbol{U}} \boldsymbol{\Xi}_A^{\star}\boldsymbol{T})^\top \widetilde{\boldsymbol{U}} \boldsymbol{\Xi}_A^{\star}\boldsymbol{T} = \boldsymbol{T}^\top {\boldsymbol{\Xi}_A^{\star}}^\top \boldsymbol{\Xi}_A^{\star} \boldsymbol{T} = \boldsymbol{T}^\top \mtx{\Lambda}^{1/2}\mtx{S}_L^\top \mtx{S}_L\mtx{\Lambda}^{1/2}\boldsymbol{T} = \boldsymbol{T}^\top \mtx{\Lambda} \boldsymbol{T}
 \\
 &= \boldsymbol{T}^\top {\boldsymbol{\Xi}_B^{\star}}^\top \boldsymbol{\Xi}_B^{\star} \boldsymbol{T} = (\widetilde{\boldsymbol{V}} \boldsymbol{\Xi}_B^{\star}\boldsymbol{T})^\top \widetilde{\boldsymbol{V}} \boldsymbol{\Xi}_B^{\star}\boldsymbol{T} = \boldsymbol{Y}(\boldsymbol{\xi})^\top \boldsymbol{Y}(\boldsymbol{\xi}).
\end{align*}
Here we use the fact $\mtx{S}_L^\top \mtx{S}_L = \mtx{S}_R^\top \mtx{S}_R = \boldsymbol{I}_r$. Moreover, 
\begin{align*}
\boldsymbol{X}(\boldsymbol{\theta})^\top \boldsymbol{X}(\boldsymbol{\xi}) + \boldsymbol{Y}(\boldsymbol{\theta})^\top \boldsymbol{Y}(\boldsymbol{\xi}) &=  (\widetilde{\boldsymbol{U}} \boldsymbol{\Theta}_A)^\top \widetilde{\boldsymbol{U}} \boldsymbol{\Xi}_A^{\star}\boldsymbol{T} + (\widetilde{\boldsymbol{V}} \boldsymbol{\Theta}_B)^\top \widetilde{\boldsymbol{V}}\boldsymbol{\Xi}_B^{\star}\boldsymbol{T} = (\boldsymbol{\Theta}_A^\top \boldsymbol{\Xi}_A^{\star} + \boldsymbol{\Theta}_B^\top \boldsymbol{\Xi}_B^{\star})\boldsymbol{T} \succeq \boldsymbol{0}.
\end{align*}
Therefore, the parameterization $(\boldsymbol{X}(\boldsymbol{\theta}),\boldsymbol{Y}(\boldsymbol{\theta}))$ satisfies Assumption \ref{ass:ideal_factorization}.
\end{proof}

\subsubsection{Proof of Corollary \ref{cor:subspace}}\label{sec:pf_subspace}
Since the assumptions of Theorem \ref{thm:main} are satisfied, therefore, in the event $E$ defined in Theorem \ref{thm:main},
\[
	\|\mtx{M}^\star - \widehat{\boldsymbol{M}}\|_F^2 \leqslant  \frac{C_3 r}{p^2} \psi^2.
\]
Therefore, it suffices to show that 
\begin{equation}\label{eq:subspace_proof_psi}
\psi^2 \leqslant  \frac{C_4}{C_3}\left( p s_{\max}\log n_{\max}\nu^2+ b^2 \frac{\mu_{\widetilde{U}} \mu_{\widetilde{V}}s_1s_2 }{n_1n_2} \log^2 n_{\max}\right). 
\end{equation}

By \eqref{eq:simplified_notation} and \eqref{eq:para_subspace}, we have for any $\boldsymbol{\theta}\in\mathbb{R}^d$,
\[
\colspan(\boldsymbol{X}(\boldsymbol{\theta})) \subset \colspan(\widetilde{\mtx{U}}), \quad \colspan(\boldsymbol{Y}(\boldsymbol{\theta})) \subset \colspan(\widetilde{\mtx{V}}).
\]
Therefore, for any $\boldsymbol{\theta}_1,\boldsymbol{\theta}_2\in\mathbb{R}^d$,
\[
	 \|\boldsymbol{P}_{\boldsymbol{X}(\boldsymbol{\theta}_1)}\mathcal{P}_{\Omega}(\boldsymbol{N}) \boldsymbol{P}_{\boldsymbol{Y}(\boldsymbol{\theta}_2)}\| = \|\boldsymbol{P}_{\boldsymbol{X}(\boldsymbol{\theta}_1)}\boldsymbol{P}_{\widetilde{\boldsymbol{U}}}\mathcal{P}_{\Omega}(\boldsymbol{N})\boldsymbol{P}_{\widetilde{\boldsymbol{V}}} \boldsymbol{P}_{\boldsymbol{Y}(\boldsymbol{\theta}_2)}\| \leqslant \|\boldsymbol{P}_{\widetilde{\boldsymbol{U}}}\mathcal{P}_{\Omega}(\boldsymbol{N})\boldsymbol{P}_{\widetilde{\boldsymbol{V}}}\|.
\]
So we have 
\[
	\psi \leqslant \|\boldsymbol{P}_{\widetilde{\boldsymbol{U}}}\mathcal{P}_{\Omega}(\boldsymbol{N})\boldsymbol{P}_{\widetilde{\boldsymbol{V}}}\|.
\]
Therefore, \eqref{eq:subspace_proof_psi} can be proved by the following lemma.
 \begin{lemma}
 \label{lem:cons_noise}
 Assume that the support of observation $\Omega$ follows from Model \ref{mod:sampling1}. We assume that the entries of the noise matrix $\boldsymbol{N}$ are i.i.d.\ centered sub-exponential random variables satisfying the Bernstein condition with parameter $b$ and variance $\nu^2$. $\widetilde{\mtx{U}}$ and $\widetilde{\mtx{V}}$ are defined in Corollary \ref{cor:subspace}. Then in an event $E_{subspace\_noise}$ with probability $\mathbb{P}[E_{subspace\_noise}]\geqslant  1- (n_1+n_2)^{-3}$, we have 
 \[
 \|\boldsymbol{P}_{\widetilde{\boldsymbol{U}}}\mathcal{P}_{\Omega}(\boldsymbol{N})\boldsymbol{P}_{\widetilde{\boldsymbol{V}}}\| \leqslant C_w \left( \sqrt{p (s_1+s_2)\log (n_1+n_2)}\nu+ b\sqrt{\frac{\mu_{\widetilde{U}} \mu_{\widetilde{V}}s_1s_2 }{n_1n_2}} \log (n_1+n_2)\right).
 \] 
for some absolute constant $C_w$ defined in the proof. 
\end{lemma}

The proof of Lemma \ref{lem:cons_noise} is mainly following the discussion in \citet[Example 6.18]{wainwright2019high} as well as \citet[Example 6.14]{wainwright2019high} and is deferred to Appendix \ref{sec:proof_lemma_subspace_noise}.

Letting $E_{subspace} = E\cap E_{subspace\_noise}$, and $C_4 = 2C_3C_w^2$ finishes the proof.

\subsection{Analysis for Skew-Symmetric Matrix Completion}
\label{sec:skew_proof}
In this section, we first give a proof of Lemma \ref{lem:rotation_positive_skew}. Then we give a proof of Theorem \ref{thm:skew}. 

\subsubsection{Proof of Lemma \ref{lem:rotation_positive_skew}}
\label{sec:rotation_positive_skew}
\begin{proof}
Recall that $\mtx{M}^\star$ is a rank-$r$ skew-symmetric matrix, where $r$ is even. 
Then its Youla decomposition \citep{youla1961normal} can be written as
\[
\mtx{M}^\star := \lambda_1 \vct{\phi}_1 \vct{\psi}_1^\top - \lambda_1 \vct{\psi}_1 \vct{\phi}_1^\top + \lambda_2 \vct{\phi}_2 \vct{\psi}_2^\top - \lambda_2 \vct{\psi}_2 \vct{\phi}_2^\top + \ldots + \lambda_{r/2} \vct{\phi}_{r/2} \vct{\psi}_{r/2}^\top - \lambda_{r/2} \vct{\psi}_{r/2} \vct{\phi}_{r/2}^\top,
\]
where $\lambda_1 \geqslant \lambda_2 \geqslant \ldots \geqslant \lambda_{r/2} > 0$ and $\vct{\phi}_1, \ldots, \vct{\phi}_{r/2}, \vct{\psi}_1, \ldots, \vct{\psi}_{r/2}$ are unit vectors in $\mathbb{R}^n$. Moreover, $\vct{\phi}_i$'s and $\boldsymbol{\psi}_i$'s are pairwise perpendicular to each other, 
i.e., 
for any $i,j\in[r/2]$, $\boldsymbol{\phi}_i^\top \boldsymbol{\psi}_j = 0$, $\boldsymbol{\phi}_i^\top \boldsymbol{\phi}_j = 0$ if $i\neq j$, and $\boldsymbol{\psi}_i^\top \boldsymbol{\psi}_j = 0$ if $i\neq j$.


  Let 
\[
\mtx{\Xi}_A^\star = [\sqrt{\lambda_1}\vct{\phi}_1, \ldots, \sqrt{\lambda_{r/2}} \vct{\phi}_{r/2}] \in \mathbb{R}^{n \times \frac{r}{2}} 
\quad\text{~and~}\quad 
\mtx{\Xi}_B^\star = [\sqrt{\lambda_1}\vct{\psi}_1, \ldots, \sqrt{\lambda_{r/2}} \vct{\psi}_{r/2}] \in \mathbb{R}^{n \times \frac{r}{2}}.
\]
It is straightforward to verify that 
\begin{equation*}
\mtx{M}^\star  =\mtx{\Xi}_A^\star {\mtx{\Xi}_B^\star}^\top - \mtx{\Xi}_B^\star {\mtx{\Xi}_A^\star}^\top. 
\end{equation*}
Recall the fact that for any $i,j\in[r/2]$, $\boldsymbol{\phi}_i^\top \boldsymbol{\psi}_j = 0$; $\boldsymbol{\phi}_i^\top \boldsymbol{\phi}_j = 0$ if $i\neq j$ and $\boldsymbol{\phi}_i^\top \boldsymbol{\phi}_j = 1$ if $i = j$; $\boldsymbol{\psi}_i^\top \boldsymbol{\psi}_j = 0$ if $i\neq j$ and $\boldsymbol{\psi}_i^\top \boldsymbol{\psi}_j = 1$ if $i = j$. Therefore, 
\begin{equation}
\label{eq:par_balance}
{\mtx{\Xi}_A^\star}^\top \mtx{\Xi}_B^\star= \mtx{0} \text{~~and~~} {\mtx{\Xi}_A^\star}^\top \mtx{\Xi}_A^\star = {\mtx{\Xi}_B^\star}^\top \mtx{\Xi}_B^\star = \text{diag}(\lambda_1, \ldots, \lambda_{r/2}).
\end{equation}

For any $\vct{\theta} = \text{vec}(\mtx{\Theta}_A, \mtx{\Theta}_B)$ with $\mtx{\Theta}_A, \mtx{\Theta}_B \in \mathbb{R}^{n \times \frac{r}{2}}$, consider the singular value decomposition of the complex matrix $(\mtx{\Theta}_A  + \sqrt{-1}\mtx{\Theta}_B)^H (\mtx{\Xi}_A^\star + \sqrt{-1}\mtx{\Xi}_B^\star)$ ($\mtx{A}^H$ is conjugate transpose of complex matrix $\mtx{A}$), 
\[
(\mtx{\Theta}_A  + \sqrt{-1} \mtx{\Theta}_B)^H (\mtx{\Xi}_A^\star + \sqrt{-1}\mtx{\Xi}_B^\star) = \boldsymbol{A}\boldsymbol{D}\boldsymbol{B}^H, 
\]
where $\boldsymbol{A},\boldsymbol{B}\in\mathbb{C}^{\frac{r}{2}\times \frac{r}{2}}$ are complex unitary matrices and $\boldsymbol{D}\in\mathbb{R}^{\frac{r}{2}\times \frac{r}{2}}$ is a real diagonal matrix. Therefore, $\boldsymbol{B}\boldsymbol{A}^H$ is also a complex unitary matrix, decompose it as 
\[
\boldsymbol{B}\boldsymbol{A}^H = \mtx{R}_1 + \sqrt{-1}\mtx{R}_2
\]
with $\mtx{R}_1,\mtx{R}_2\in \mathbb{R}^{\frac{r}{2}\times \frac{r}{2}}$. Therefore, 
\[
  (\mtx{\Theta}_A  + \sqrt{-1}\mtx{\Theta}_B)^H (\mtx{\Xi}_A^\star + \sqrt{-1}\mtx{\Xi}_B^\star)(\mtx{R}_1 + \sqrt{-1}\mtx{R}_2) = \boldsymbol{A}\boldsymbol{D}\boldsymbol{B}^H \boldsymbol{B}\boldsymbol{A}^H = \boldsymbol{A}\boldsymbol{D} \boldsymbol{A}^H \succeq \mtx{0},
\]
 that is, it is a Hermitian positive semidefinite matrix. Let 
 \[
 {\mtx{\Xi}_A} = \mtx{\Xi}_A^\star \mtx{R}_1 - \mtx{\Xi}_B^\star \mtx{R}_2 \text{~~and~~} {\mtx{\Xi}_B} = \mtx{\Xi}_A^\star \mtx{R}_2 + \mtx{\Xi}_B^\star \mtx{R}_1.
 \]
 Then it holds that 
 \[
 (\mtx{\Xi}_A^\star + \sqrt{-1} \mtx{\Xi}_B^\star)(\mtx{R}_1 + \sqrt{-1}\mtx{R}_2) = \mtx{\Xi}_A + \sqrt{-1} \mtx{\Xi}_B
 \]
 and 
 \[
 (\mtx{\Theta}_A + \sqrt{-1} \mtx{\Theta}_B)^H (\mtx{\Xi}_A + \sqrt{-1} \mtx{\Xi}_B) \succeq \mtx{0}, 
 \]
 which is equivalent to the following $r$-by-$r$ real matrix is positive semidefinite:
\[
\begin{bmatrix}
\mtx{\Theta}_A^\top {\mtx{\Xi}_A} + \mtx{\Theta}_B^\top {\mtx{\Xi}_B} & \mtx{\Theta}_B^\top {\mtx{\Xi}_A} - \mtx{\Theta}_A^\top {\mtx{\Xi}_B}
\\
\mtx{\Theta}_A^\top {\mtx{\Xi}_B} - \mtx{\Theta}_B^\top \mtx{\Xi}_A & \mtx{\Theta}_A^\top {\mtx{\Xi}_A} + \mtx{\Theta}_B^\top {\mtx{\Xi}_B}
\end{bmatrix}
 \succeq \mtx{0}.
\]
Also, since $\mtx{R}_1 + \sqrt{-1} \mtx{R}_2$ is unitary, we have 
\[
 \mtx{R}:= \begin{bmatrix} \mtx{R}_1 & -\mtx{R}_2 \\ \mtx{R}_2 & \mtx{R}_1 \end{bmatrix} \in \mathsf{O}(r).
 \]

Let $\vct{\xi} = \text{vec}(\mtx{\Xi}_A, \mtx{\Xi}_B)$. Then we have 
\begin{align*}
\mtx{X}(\vct{\xi}) &= \left[\mtx{\Xi}_A, -\mtx{\Xi}_B\right] = \left[\mtx{\Xi}_A^\star \mtx{R}_1 - \mtx{\Xi}_B^\star \mtx{R}_2, -\mtx{\Xi}_A^\star \mtx{R}_2 - \mtx{\Xi}_B^\star \mtx{R}_1\right] 
\\
&= \left[\mtx{\Xi}_A^\star, -\mtx{\Xi}_B^\star\right] \begin{bmatrix} 
\mtx{R}_1 & -\mtx{R}_2 \\ 
\mtx{R}_2 & \mtx{R}_1 
\end{bmatrix} 
= \left[\mtx{\Xi}_A^\star, -\mtx{\Xi}_B^\star\right] \mtx{R}, 
\end{align*}
and similarly
\[
\mtx{Y}(\vct{\xi}) =\left[\mtx{\Xi}_B, \mtx{\Xi}_A\right] = \left[\mtx{\Xi}_A^\star \mtx{R}_2 + \mtx{\Xi}_B^\star \mtx{R}_1, \mtx{\Xi}_A^\star \mtx{R}_1 - \mtx{\Xi}_B^\star \mtx{R}_2\right] = \left[\mtx{\Xi}_B^\star, \mtx{\Xi}_A^\star\right] \begin{bmatrix} \mtx{R}_1 & -\mtx{R}_2 \\ \mtx{R}_2 & \mtx{R}_1 \end{bmatrix} = \left[\mtx{\Xi}_B^\star, \mtx{\Xi}_A^\star\right]  \mtx{R}.
\]
It is then straightforward to verify that
\[
\mtx{X}(\vct{\xi})\mtx{Y}(\vct{\xi})^\top = \left[\mtx{\Xi}_A, -\mtx{\Xi}_B\right]\left[\mtx{\Xi}_B, \mtx{\Xi}_A\right]^\top = \left[\mtx{\Xi}_A^\star, -\mtx{\Xi}_B^\star\right]\left[\mtx{\Xi}_B^\star, \mtx{\Xi}_A^\star \right]^\top = \mtx{M}^\star .
\]
In order to further verify $\mtx{X}(\vct{\xi})^\top \mtx{X}(\vct{\xi}) = \mtx{Y}(\vct{\xi})^\top \mtx{Y}(\vct{\xi})$, it suffices to prove 
\[
\left[\mtx{\Xi}_A^\star, -\mtx{\Xi}_B^\star\right]^\top\left[\mtx{\Xi}_A^\star, -\mtx{\Xi}_B^\star\right] =  \left[\mtx{\Xi}_B^\star, \mtx{\Xi}_A^\star\right]^\top \left[\mtx{\Xi}_B^\star, \mtx{\Xi}_A^\star\right], 
\]
which is guaranteed by ${\mtx{\Xi}_A^\star}^\top \mtx{\Xi}_B^\star= \mtx{0}$ and ${\mtx{\Xi}_A^\star}^\top \mtx{\Xi}_A^\star = {\mtx{\Xi}_B^\star}^\top \mtx{\Xi}_B^\star$ as was shown in \eqref{eq:par_balance}.

Finally, straightforward calculation gives
\begin{align*}
\mtx{X}(\vct{\theta})^\top \mtx{X}(\vct{\xi}) + \mtx{Y}(\vct{\theta})^\top \mtx{Y}(\vct{\xi}) &= [\mtx{\Theta}_A, -\mtx{\Theta}_B]^\top \left[\mtx{\Xi}_A, -\mtx{\Xi}_B\right] + [\mtx{\Theta}_B, \mtx{\Theta}_A]^\top \left[\mtx{\Xi}_B, \mtx{\Xi}_A\right] \\
&=\begin{bmatrix}
\mtx{\Theta}_A^\top {\mtx{\Xi}_A} + \mtx{\Theta}_B^\top {\mtx{\Xi}_B} & \mtx{\Theta}_B^\top {\mtx{\Xi}_A} - \mtx{\Theta}_A^\top {\mtx{\Xi}_B}
\\
\mtx{\Theta}_A^\top \mtx{\Xi}_B - \mtx{\Theta}_B^\top \mtx{\Xi}_A & \mtx{\Theta}_A^\top {\mtx{\Xi}_A} + \mtx{\Theta}_B^\top {\mtx{\Xi}_B}
\end{bmatrix}
 \succeq \mtx{0}.
\end{align*}
\end{proof}

\subsubsection{Proof of Theorem \ref{thm:skew}}
Following the lines in Section \ref{sec:pf_subspace}, it suffices to show that
\[
	\psi^2 \leqslant \frac{ C_5 }{C_3} \left( p n \log n\nu^2+ b^2   \log^2 n \right).
\]
Recall the fact that $\psi \leqslant \|\mathcal{P}_{\Omega}(\boldsymbol{N})\|$, then the proof can be done by employing the following Lemma.
\begin{lemma}
	Let the support of the observed entries $\Omega$ satisfy Model \ref{mod:sampling2}. We assume that the noise matrix $\boldsymbol{N}$ is a skew-symmetric matrix, and upper triangular part of $\boldsymbol{N}$ consists of i.i.d.\ centered sub-exponential random variables satisfying the Bernstein condition with parameter $b$ and variance $\nu^2$. Then in an event $E_{skew\_noise}$ with probability $\mathbb{P}[E_{skew\_noise}]\geqslant  1- n^{-3}$, we have 
 \[
 \| \mathcal{P}_{\Omega}(\boldsymbol{N}) \| \leqslant C_{w'} \left( \sqrt{p n \log n }\nu+ b  \log n\right).
 \] 
for some absolute constant $C_{w'}$.

\end{lemma}
The proof is almost exactly the same with proof of Lemma \ref{lem:cons_noise}, thus omitted here. We can finish the proof of Theorem \ref{thm:skew} by letting $E_{skew} = E\cap E_{skew\_noise}$ and $C_5 = 2C_3C_{w'}^2$.

\section{Discussion}
\label{sec:discussion}
This paper proposes a unified nonconvex optimization framework for matrix completion with linearly parameterized factors. Examples include subspace-constrained and skew-symmetric matrix completion, with applications in collaborative filtering with side information and pairwise comparisons. We also conduct a unified geometric analysis for this nonconvex optimization framework, where a novel assumption referred to as \emph{correlated parametric factorization} plays a key role. In particular, for the noiseless case, we give a sufficient condition on the sampling rate to guarantee that with high probability no spurious local minimum exists. Our meta-theorem is applicable to examples including subspace-constrained and skew-symmetric matrix completion, given in either case the assumption of correlated parametric factorization is verified. 

For future work, we are particularly interested in extending the current framework to more general parameterized factorization, such as extension from linear mappings to affine mappings or even nonlinear mappings. The geometrical interpretation of the assumption of correlated parametric factorization is still intriguing. Its verification for subspace-constrained matrix completion and skew-symmetric matrix completion relies on very different algebraic arguments. Are there any underlying geometric reasons behind such coincidence? Due to its applications in practice, we are particularly interested in the geometric analysis for nonconvex subspace-constrained matrix completion. As explained in Section \ref{sec:subspace_constraint_matrix} and illustrated in our simulations, our current sampling rate result is suboptimal and we are interested in improving upon it in future.



\section*{Acknowledgements}
J. Chen and X. Li acknowledge support from the NSF CAREER Award DMS-1848575. J. Chen also acknowledges support from the NSF HDR TRIPODS grant CCF-1934568. We would like to thank an anonymous reviewer who provides ideas for reformulating the proof of Theorem \ref{thm:main}.

\bibliography{cite} 
\bibliographystyle{plainnat}
\newpage

\begin{appendix}

\section{Proofs of Auxiliary Results used in the Proof of Theorem \ref{thm:main}}  
\subsection{Preliminaries} 

First, we list two useful lemmas controlling the difference between the random sampled matrix inner product (or Frobenius norm) and its expectation.
 
\begin{lemma}[{\citealt{chen2017memory}}]
\label{lem:determin}
It holds that uniformly for all $\boldsymbol{A},\boldsymbol{B}\in\mathbb{R}^{n_1\times r}$, $\boldsymbol{C},\boldsymbol{D}\in\mathbb{R}^{n_2\times r}$
\begin{equation}\label{eq:determin}
  \begin{split}
    & | \langle \mathcal{P}_{\Omega}(\boldsymbol{A}\boldsymbol{C}^\top),\mathcal{P}_{\Omega}(\boldsymbol{B}\boldsymbol{D}^\top)  \rangle - p \langle\boldsymbol{A}\boldsymbol{C}^\top,\boldsymbol{B}\boldsymbol{D}^\top \rangle| \\
    \leqslant  & \|\boldsymbol{\Omega}-p\boldsymbol{J}\|\sqrt{\sum_{k=1}^{n_1}\|\boldsymbol{A}_{k,\cdot}\|_2^2 \|\boldsymbol{B}_{k,\cdot}\|_2^2 }\sqrt{\sum_{k=1}^{n_2}\|\boldsymbol{C}_{k,\cdot}\|_2^2 \|\boldsymbol{D}_{k,\cdot}\|_2^2 }\\
\leqslant  & \frac{1}{2}\|\boldsymbol{\Omega}-p\boldsymbol{J}\|
\left(\sum_{k=1}^{n_1}\|\boldsymbol{A}_{k,\cdot}\|_2^2 \|\boldsymbol{B}_{k,\cdot}\|_2^2 
+ \sum_{k=1}^{n_2}\|\boldsymbol{C}_{k,\cdot}\|_2^2 \|\boldsymbol{D}_{k,\cdot}\|_2^2\right).
  \end{split} 
\end{equation}
Here $\boldsymbol{J}$ denotes the matrix with all entries equal to one.
\end{lemma}

\begin{lemma}[{\citealt{candes2009exact}}]
\label{lem:rip}
Denote 
\[
\mathcal{T} \coloneqq \left\{ \boldsymbol{M}\in\mathbb{R}^{n_1\times n_2} \mid \boldsymbol{M} = \boldsymbol{U}^{\star} \boldsymbol{A}^\top + \boldsymbol{B}{\boldsymbol{V}^{\star} }^\top, \boldsymbol{A}\in\mathbb{R}^{n_2\times r},\boldsymbol{B}\in\mathbb{R}^{n_1\times r} \right\}. 
\]
If $p\geqslant  C_c\frac{ \mu r\log n_{\max}}{n_{\min}}$ for some absolute constant $C_c$, then in an event $E_c$ with probability $\mathbb{P}[E_c] \geqslant  1- (n_1+n_2)^{-5}$,
\[
\left\|\frac{1}{p}\mathcal{P}_{\mathcal{T}}\mathcal{P}_{\Omega}\mathcal{P}_{\mathcal{T}} - \mathcal{P}_{\mathcal{T}}\right\|\leqslant  10^{-4}.
\]
Here $\mathcal{P}_{\Omega}$ is determined as in Models \ref{mod:sampling1} and \ref{mod:sampling2}, where $\mathcal{P}_{\mathcal{T}}$ is the orthogonal projector on $\mathcal{T}$ with the Euclidean spaced defined by the standard matrix inner product. In other words, uniformly for all $\mtx{M} \in \mathcal{T}$,
\begin{equation}\label{eq:rip}
 \left| \left\|\frac{1}{p}\mathcal{P}_{\Omega}(\mtx{M}) \right\|_F^2 -  \left\| \mtx{M} \right\|_F^2 \right| \leqslant  10^{-4} \left\| \mtx{M}\right\|_F^2.
\end{equation}
\end{lemma}
The following lemma on $\|\boldsymbol{\Omega}-p\boldsymbol{J}\|$ is well known in the literature, see, e.g., \citet{vu2014simple} and \citet{bandeira2016sharp}:
\begin{lemma}
\label{lem:vu14_psd}
There is an absolute constant $C_v>0$, such that if $p\geqslant C_v\frac{\log n_{\max}}{n_{\min}}$, on an event with probability at least $1 - (n_1+n_2)^{-5}$, there holds
\begin{equation}
\label{eq:omega_operator}
\|\boldsymbol{\Omega}-p\boldsymbol{J}\|<   C_v \sqrt{n_{\max}p}.  
\end{equation}
\end{lemma} 
We use the following lemma to specify the event $E$ in Theorem \ref{thm:main}.
\begin{lemma}
\label{lemma:eq_of_rate}
  Under the assumptions in Theorem \ref{thm:main}, in an event $E$ with probability $\mathbb{P}[E] \geqslant 1-(n_1+n_2)^{-3}$,
\begin{equation}
  \label{eq:rate_tuning_ass}
  \begin{cases}
    p\geqslant  C_0 \frac{\mu r \log n_{\max}}{ n_{\min}},
    \\
    \alpha \geqslant  100 \sqrt{\frac{\mu r \sigma_1}{n_{\min}}},
    \\
    \sigma_r  \geqslant C_0(\frac{\mu r \sigma_1 \|\boldsymbol{\Omega}-p\boldsymbol{J}\|}{p n_{\min}} + \lambda\alpha^2) ,
    \\
    p \lambda \geqslant C_0 \|\boldsymbol{\Omega}-p\boldsymbol{J}\|
    \end{cases}
  \end{equation}
as well as \eqref{eq:rip} in Lemma \ref{lem:rip} and \eqref{eq:omega_operator} in Lemma \ref{lem:vu14_psd} hold. Here $C_0 = 5 \times 10^5$.
\end{lemma}

%

Finally, we collect some useful properties of $\mtx{U}= \mtx{X}(\vct{\xi})$ and $\mtx{V} = \mtx{Y}(\vct{\xi})$. The proof is left to Section \ref{sec:proof_of_prop_basic_properties}.
\begin{proposition}
\label{prop:basic_properties}
For any $\vct{\theta}$, the matrices $\mtx{U}= \mtx{X}(\vct{\xi})$ and $\mtx{V} = \mtx{Y}(\vct{\xi})$ defined in \eqref{eq:simplified_notation} satisfy the following basic properties:
\begin{itemize}
\item $\colspan(\boldsymbol{U}) = \colspan(\boldsymbol{U}^{\star} )$ and $\colspan(\boldsymbol{V}) = \colspan(\boldsymbol{V}^{\star} )$;
\item The largest singular values of both $\mtx{U}$ and $\mtx{V}$ are $\sqrt{\sigma_1}$;
\item The $r$-th singular values of both $\mtx{U}$ and $\mtx{V}$ are $\sqrt{\sigma_r}$. 
\item $\|\boldsymbol{U}\|_{2,\infty}^2  \leqslant   \frac{\mu r }{n_1}\sigma_1$ and $\|\boldsymbol{V}\|_{2,\infty}^2 \leqslant  \frac{\mu r }{n_2}\sigma_1$. 
\end{itemize}
\end{proposition}

%
%

\subsection{Proof of Lemma \ref{lem:ge_aux_formula}}
\label{sec:proof_ge_aux_formula}
\begin{proof} 
Lemma \ref{lem:ge_aux_formula} is essentially \citet[Lemma 16]{ge2017no} (with noise). Here we give a sketch of the proof for the purpose of being self-contained.

First, denote $f_{\textrm{clean}}(\boldsymbol{X}(\boldsymbol{\theta}),\boldsymbol{Y}(\boldsymbol{\theta}))$ as 
\[ 
  \begin{split}
    f_{\textrm{clean}}(\boldsymbol{X}(\boldsymbol{\theta}),\boldsymbol{Y}(\boldsymbol{\theta})) =& \frac{1}{2p}\|\mathcal{P}_{\Omega}(\boldsymbol{X}(\boldsymbol{\theta})\boldsymbol{Y}(\boldsymbol{\theta})^\top-\boldsymbol{M}^\star)\|_F^2 + \frac{1}{8} \|\boldsymbol{X}(\boldsymbol{\theta})^\top\boldsymbol{X}(\boldsymbol{\theta})-\boldsymbol{Y}(\boldsymbol{\theta})^\top\boldsymbol{Y}(\boldsymbol{\theta})\|_F^2\\
    & +\lambda (G_{\alpha}(\boldsymbol{X}(\boldsymbol{\theta})) + G_{\alpha}(\boldsymbol{Y}(\boldsymbol{\theta}))).
  \end{split} 
\]
Compare with \eqref{eq:obj_para}, and use the simplified notations introduced in \eqref{eq:simplified_notation}. We can see
\[
  f(\boldsymbol{X},\boldsymbol{Y}) = f_{\textrm{clean}}(\boldsymbol{X},\boldsymbol{Y}) - \frac{1}{p} \langle \mathcal{P}_{\Omega}(\boldsymbol{X}\boldsymbol{Y}^\top - \boldsymbol{M}^\star),\mathcal{P}_{\Omega}(\boldsymbol{N}) \rangle + \frac{1}{2p}\|\mathcal{P}_{\Omega}(\boldsymbol{N})\|_F^2.
\]

Therefore, 
\begin{equation*}
    \langle \nabla f(\boldsymbol{X},\boldsymbol{Y}),\operatorname{vec}(\boldsymbol{\Delta}_{\boldsymbol{X}},\boldsymbol{\Delta}_{\boldsymbol{Y}}) \rangle = \langle \nabla f_{\textrm{clean}}(\boldsymbol{X},\boldsymbol{Y}),\operatorname{vec}(\boldsymbol{\Delta}_{\boldsymbol{X}},\boldsymbol{\Delta}_{\boldsymbol{Y}}) \rangle - \frac{1}{p}\langle \mathcal{P}_{\Omega}(\boldsymbol{\Delta}_{\boldsymbol{X}}\boldsymbol{Y}^\top + \boldsymbol{X}\boldsymbol{\Delta}_{\boldsymbol{Y}}^\top ) ,\mathcal{P}_{\Omega}(\boldsymbol{N})\rangle
  \end{equation*}
and 
\begin{equation*}
    \begin{split}
      & \operatorname{vec}(\boldsymbol{\Delta}_{\boldsymbol{X}},\boldsymbol{\Delta}_{\boldsymbol{Y}})^\top \nabla^2 f(\boldsymbol{X},\boldsymbol{Y}) \operatorname{vec}(\boldsymbol{\Delta}_{\boldsymbol{X}},\boldsymbol{\Delta}_{\boldsymbol{Y}})\\
       =& \operatorname{vec}(\boldsymbol{\Delta}_{\boldsymbol{X}},\boldsymbol{\Delta}_{\boldsymbol{Y}})^\top \nabla^2 f_{\textrm{clean}}(\boldsymbol{X},\boldsymbol{Y}) \operatorname{vec}(\boldsymbol{\Delta}_{\boldsymbol{X}},\boldsymbol{\Delta}_{\boldsymbol{Y}}) - \frac{2}{p} \langle \mathcal{P}_{\Omega}(\boldsymbol{\Delta}_{\boldsymbol{X}}\boldsymbol{\Delta}_{\boldsymbol{Y}}^\top),\mathcal{P}_{\Omega}(\boldsymbol{N}) \rangle.
    \end{split}
  \end{equation*}

Therefore, we only need to concern about $f_{\textrm{clean}}(\boldsymbol{X},\boldsymbol{Y})$ now, which has already been discussed in \citet{ge2017no}. Interested readers can refer to \citet{ge2017no} for the detail.  

By \citet[Lemma 16]{ge2017no}, we have 
\[
	\begin{split}
		& \operatorname{vec}(\boldsymbol{\Delta}_{\boldsymbol{X}},\boldsymbol{\Delta}_{\boldsymbol{Y}})^\top \nabla^2 f_{\textrm{clean}}(\boldsymbol{X},\boldsymbol{Y}) \operatorname{vec}(\boldsymbol{\Delta}_{\boldsymbol{X}},\boldsymbol{\Delta}_{\boldsymbol{Y}})  -4\langle \nabla f_{\textrm{clean}}(\boldsymbol{X},\boldsymbol{Y}),\operatorname{vec}(\boldsymbol{\Delta}_{\boldsymbol{X}},\boldsymbol{\Delta}_{\boldsymbol{Y}}) \rangle  \\ 
		 \leqslant  & \frac{1}{4} \left\{ \left\|\boldsymbol{\Delta}_{\boldsymbol{Z}}\boldsymbol{\Delta}_{\boldsymbol{Z}}^\top\right\|_F^2 - 3\left\| \boldsymbol{Z}\boldsymbol{Z}^\top -\boldsymbol{W}\boldsymbol{W}^\top \right\|_F^2 \right\}  + \left(\frac{1}{p}\left\| \mathcal{P}_{\Omega}\left(\boldsymbol{\Delta}_{\boldsymbol{X}}\boldsymbol{\Delta}_{\boldsymbol{Y}}^\top \right) \right\|_F^2 - \|\boldsymbol{\Delta}_{\boldsymbol{X}}\boldsymbol{\Delta}_{\boldsymbol{Y}}^\top \|_F^2 \right) \\
    &- \left( \frac{3}{p}\left\| \mathcal{P}_{\Omega}\left( \boldsymbol{X}\boldsymbol{Y}^\top - \boldsymbol{U}\boldsymbol{V}^\top  \right) \right\|_F^2 - 3 \|\boldsymbol{X}\boldsymbol{Y}^\top - \boldsymbol{U}\boldsymbol{V}^\top\|_F^2 \right)\\
    &+ \lambda \left[\operatorname{vec}(\boldsymbol{\Delta}_{\boldsymbol{X}})^\top \nabla^2 G_{\alpha}(\boldsymbol{X}) \operatorname{vec}(\boldsymbol{\Delta}_{\boldsymbol{X}})- 4\left\langle \nabla G_{\alpha}(\boldsymbol{X}),\boldsymbol{\Delta}_{\boldsymbol{X}} \right\rangle  \right]\\
    &+ \lambda \left[\operatorname{vec}(\boldsymbol{\Delta}_{\boldsymbol{Y}})^\top \nabla^2 G_{\alpha}(\boldsymbol{Y}) \operatorname{vec}(\boldsymbol{\Delta}_{\boldsymbol{Y}})- 4\left\langle \nabla G_{\alpha}(\boldsymbol{Y}),\boldsymbol{\Delta}_{\boldsymbol{Y}} \right\rangle  \right].
	\end{split}
\]

Therefore,
\[
  \begin{split}
    &\operatorname{vec}(\boldsymbol{\Delta}_{\boldsymbol{X}},\boldsymbol{\Delta}_{\boldsymbol{Y}})^\top \nabla^2 f(\boldsymbol{X},\boldsymbol{Y}) \operatorname{vec}(\boldsymbol{\Delta}_{\boldsymbol{X}},\boldsymbol{\Delta}_{\boldsymbol{Y}})   - 4\langle \nabla f(\boldsymbol{X},\boldsymbol{Y}),\operatorname{vec}(\boldsymbol{\Delta}_{\boldsymbol{X}},\boldsymbol{\Delta}_{\boldsymbol{Y}}) \rangle\\
    =& \operatorname{vec}(\boldsymbol{\Delta}_{\boldsymbol{X}},\boldsymbol{\Delta}_{\boldsymbol{Y}})^\top \nabla^2 f_{\textrm{clean}}(\boldsymbol{X},\boldsymbol{Y}) \operatorname{vec}(\boldsymbol{\Delta}_{\boldsymbol{X}},\boldsymbol{\Delta}_{\boldsymbol{Y}}) - \frac{2}{p} \langle \mathcal{P}_{\Omega}(\boldsymbol{\Delta}_{\boldsymbol{X}}\boldsymbol{\Delta}_{\boldsymbol{Y}}^\top),\mathcal{P}_{\Omega}(\boldsymbol{N}) \rangle \\
    & -4\langle \nabla f_{\textrm{clean}}(\boldsymbol{X},\boldsymbol{Y}),\operatorname{vec}(\boldsymbol{\Delta}_{\boldsymbol{X}},\boldsymbol{\Delta}_{\boldsymbol{Y}}) \rangle + \frac{4}{p}\langle \mathcal{P}_{\Omega}(\boldsymbol{\Delta}_{\boldsymbol{X}}\boldsymbol{Y}^\top + \boldsymbol{X}\boldsymbol{\Delta}_{\boldsymbol{Y}}^\top ) ,\mathcal{P}_{\Omega}(\boldsymbol{N})\rangle \\  
    \leqslant  & \frac{4}{p}\langle \mathcal{P}_{\Omega}(\boldsymbol{\Delta}_{\boldsymbol{X}}\boldsymbol{V}^\top + \boldsymbol{U}\boldsymbol{\Delta}_{\boldsymbol{Y}}^\top ) ,\mathcal{P}_{\Omega}(\boldsymbol{N})\rangle+\frac{6}{p} \langle \mathcal{P}_{\Omega}(\boldsymbol{\Delta}_{\boldsymbol{X}}\boldsymbol{\Delta}_{\boldsymbol{Y}}^\top),\mathcal{P}_{\Omega}(\boldsymbol{N}) \rangle \\
    &  +\frac{1}{4} \left\{ \left\|\boldsymbol{\Delta}_{\boldsymbol{Z}}\boldsymbol{\Delta}_{\boldsymbol{Z}}^\top\right\|_F^2 - 3\left\| \boldsymbol{Z}\boldsymbol{Z}^\top -\boldsymbol{W}\boldsymbol{W}^\top \right\|_F^2 \right\}  + \left(\frac{1}{p}\left\| \mathcal{P}_{\Omega}\left(\boldsymbol{\Delta}_{\boldsymbol{X}}\boldsymbol{\Delta}_{\boldsymbol{Y}}^\top \right) \right\|_F^2 - \|\boldsymbol{\Delta}_{\boldsymbol{X}}\boldsymbol{\Delta}_{\boldsymbol{Y}}^\top \|_F^2 \right) \\
    &- \left( \frac{3}{p}\left\| \mathcal{P}_{\Omega}\left( \boldsymbol{X}\boldsymbol{Y}^\top - \boldsymbol{U}\boldsymbol{V}^\top  \right) \right\|_F^2 - 3 \|\boldsymbol{X}\boldsymbol{Y}^\top - \boldsymbol{U}\boldsymbol{V}^\top\|_F^2 \right)\\
    &+ \lambda \left[\operatorname{vec}(\boldsymbol{\Delta}_{\boldsymbol{X}})^\top \nabla^2 G_{\alpha}(\boldsymbol{X}) \operatorname{vec}(\boldsymbol{\Delta}_{\boldsymbol{X}})- 4\left\langle \nabla G_{\alpha}(\boldsymbol{X}),\boldsymbol{\Delta}_{\boldsymbol{X}} \right\rangle  \right]\\
    &+ \lambda \left[\operatorname{vec}(\boldsymbol{\Delta}_{\boldsymbol{Y}})^\top \nabla^2 G_{\alpha}(\boldsymbol{Y}) \operatorname{vec}(\boldsymbol{\Delta}_{\boldsymbol{Y}})- 4\left\langle \nabla G_{\alpha}(\boldsymbol{Y}),\boldsymbol{\Delta}_{\boldsymbol{Y}} \right\rangle  \right].
  \end{split}  
\]
Combining with \eqref{eq:auxiliary_rect} and Lemma \ref{lem:key_identity} finishes the proof.
\end{proof}

\subsection{Proof of Lemma \ref{lem:K2K3}}
\label{sec:pf-lem-K2K3}
This section is meant to control $K_2+K_3$.
The arguments follow closely those used in \citet{chen2017memory} with modifications pertinent to the correlated parametric factorization structure in Assumption \ref{ass:ideal_factorization}. 
Recall that we assume event $E$ defined in Lemma \ref{lemma:eq_of_rate} holds, so that \eqref{eq:determin}, \eqref{eq:rip}, \eqref{eq:omega_operator} and \eqref{eq:rate_tuning_ass} hold simultaneously.
In what follows, we shall control $K_2$ and $K_3$ respectively first before piecing the bounds together.


\subsubsection{Control of $K_2$}
By the way we define $\boldsymbol{\Delta}_{\boldsymbol{X}}, \boldsymbol{\Delta}_{\boldsymbol{Y}}$ in \eqref{eq:simplified_notation},
\begin{align*}
	\boldsymbol{X}\boldsymbol{Y}^\top - \boldsymbol{U}\boldsymbol{V}^\top & = (\boldsymbol{U}+ \boldsymbol{\Delta}_{\boldsymbol{X}}) (\boldsymbol{V}+ \boldsymbol{\Delta}_{\boldsymbol{Y}})^\top - \boldsymbol{U}\boldsymbol{V}^\top  \\
& =  \boldsymbol{\Delta}_{\boldsymbol{X}}\boldsymbol{V}^\top + \boldsymbol{U}\boldsymbol{\Delta}_{\boldsymbol{Y}}^\top + \boldsymbol{\Delta}_{\boldsymbol{X}}\boldsymbol{\Delta}_{\boldsymbol{Y}}^\top.	
\end{align*}
Therefore,
\begin{align*}
&~~~\,\left| \frac{1}{p}\left\| \mathcal{P}_{\Omega}\left( \boldsymbol{X}\boldsymbol{Y}^\top - \boldsymbol{U}\boldsymbol{V}^\top  \right) \right\|_F^2 -  \|\boldsymbol{X}\boldsymbol{Y}^\top - \boldsymbol{U}\boldsymbol{V}^\top\|_F^2 \right|
\\
& =\left| \frac{1}{p}\left\|\mathcal{P}_{\Omega}\left( \boldsymbol{\Delta}_{\boldsymbol{X}}\boldsymbol{V}^\top + \boldsymbol{U}\boldsymbol{\Delta}_{\boldsymbol{Y}}^\top + \boldsymbol{\Delta}_{\boldsymbol{X}}\boldsymbol{\Delta}_{\boldsymbol{Y}}^\top  \right)\right\|_F^2 -  \left\| \boldsymbol{\Delta}_{\boldsymbol{X}}\boldsymbol{V}^\top + \boldsymbol{U}\boldsymbol{\Delta}_{\boldsymbol{Y}}^\top + \boldsymbol{\Delta}_{\boldsymbol{X}}\boldsymbol{\Delta}_{\boldsymbol{Y}}^\top  \right\|_F^2\right|
\\
& \leqslant  \circled{1} + \circled{2}+\circled{3}+\circled{4}.
\end{align*}
Where
\begin{align*}
\circled{1} & \coloneqq \left| \frac{1}{p} \left\|\mathcal{P}_{\Omega}\left( \boldsymbol{\Delta}_{\boldsymbol{X}}\boldsymbol{V}^\top + \boldsymbol{U}\boldsymbol{\Delta}_{\boldsymbol{Y}}^\top  \right) \right\|_F^2  - \left\| \boldsymbol{\Delta}_{\boldsymbol{X}}\boldsymbol{V}^\top + \boldsymbol{U}\boldsymbol{\Delta}_{\boldsymbol{Y}}^\top \right\|_F^2\right|,
\\
\circled{2} & \coloneqq \left| \frac{1}{p}\left\| \mathcal{P}_{\Omega}\left(\boldsymbol{\Delta}_{\boldsymbol{X}}\boldsymbol{\Delta}_{\boldsymbol{Y}}^\top \right) \right\|_F^2 - \|\boldsymbol{\Delta}_{\boldsymbol{X}}\boldsymbol{\Delta}_{\boldsymbol{Y}}^\top \|_F^2 \right|,	
\\
\circled{3} & \coloneqq \left| \frac{2}{p} \left\langle\mathcal{P}_{\Omega}\left( \boldsymbol{\Delta}_{\boldsymbol{X}}\boldsymbol{V}^\top \right), \mathcal{P}_{\Omega}\left(\boldsymbol{\Delta}_{\boldsymbol{X}}\boldsymbol{\Delta}_{\boldsymbol{Y}}^\top \right)\right\rangle - 2 \left\langle  \boldsymbol{\Delta}_{\boldsymbol{X}}\boldsymbol{V}^\top,  \boldsymbol{\Delta}_{\boldsymbol{X}}\boldsymbol{\Delta}_{\boldsymbol{Y}}^\top  \right\rangle\right|,~~~ \mbox{and}\\
\circled{4} & \coloneqq \left| \frac{2}{p} \left\langle\mathcal{P}_{\Omega}\left(\boldsymbol{U}\boldsymbol{\Delta}_{\boldsymbol{Y}}^\top  \right), \mathcal{P}_{\Omega}\left(\boldsymbol{\Delta}_{\boldsymbol{X}}\boldsymbol{\Delta}_{\boldsymbol{Y}}^\top \right)\right\rangle - 2 \left\langle  \boldsymbol{U}\boldsymbol{\Delta}_{\boldsymbol{Y}}^\top,  \boldsymbol{\Delta}_{\boldsymbol{X}}\boldsymbol{\Delta}_{\boldsymbol{Y}}^\top  \right\rangle\right|.
\end{align*}

By Proposition \ref{prop:basic_properties}, the matrix $\boldsymbol{\Delta}_{\boldsymbol{X}}\boldsymbol{V}^\top + \boldsymbol{U}\boldsymbol{\Delta}_{\boldsymbol{Y}}^\top$ belongs to the subspace $\mathcal{T}$ defined in Lemma \ref{lem:rip}. 
Therefore, Lemma \ref{lem:rip} implies
\[
\circled{1} \leqslant   0.0001 \left\| \boldsymbol{\Delta}_{\boldsymbol{X}}\boldsymbol{V}^\top + \boldsymbol{U}\boldsymbol{\Delta}_{\boldsymbol{Y}}^\top \right\|_F^2 \leqslant  0.0002 \left(\left\| \boldsymbol{\Delta}_{\boldsymbol{X}}\boldsymbol{V}^\top\right\|_F^2 + \left\|\boldsymbol{U}\boldsymbol{\Delta}_{\boldsymbol{Y}}^\top \right\|_F^2\right).
\]
Applying Lemma \ref{lem:determin}, we further have
\begin{align*}
\circled{2} & \leqslant   \frac{\|\boldsymbol{\Omega}-p\boldsymbol{J}\|}{2p} \left(\sum_{k=1}^{n_1}\|(\boldsymbol{\Delta}_{\boldsymbol{X}})_{k,\cdot}\|_2^4 + \sum_{k=1}^{n_2}\|(\boldsymbol{\Delta}_{\boldsymbol{Y}})_{k,\cdot}\|_2^4  \right),
\\
\circled{3} & \leqslant   \frac{\|\boldsymbol{\Omega}-p\boldsymbol{J}\|}{p} \left(\sum_{k=1}^{n_1}\|(\boldsymbol{\Delta}_{\boldsymbol{X}})_{k,\cdot}\|_2^4 + \sum_{k=1}^{n_2} \|\boldsymbol{V}_{k,\cdot}\|_2^2\|(\boldsymbol{\Delta}_{\boldsymbol{Y}})_{k,\cdot}\|_2^2  \right),
\\
\circled{4} &\leqslant   \frac{\|\boldsymbol{\Omega}-p\boldsymbol{J}\|}{p} \left(\sum_{k=1}^{n_2}\|(\boldsymbol{\Delta}_{\boldsymbol{Y}})_{k,\cdot}\|_2^4 + \sum_{k=1}^{n_1} \|\boldsymbol{U}_{k,\cdot}\|_2^2\|(\boldsymbol{\Delta}_{\boldsymbol{X}})_{k,\cdot}\|_2^2  \right).
\end{align*}
By using Proposition \ref{prop:basic_properties}, $\|\boldsymbol{U}\|_{2,\infty}^2   \leqslant   \frac{\mu r }{n_1}\sigma_1$ and $\|\boldsymbol{V}\|_{2,\infty}^2 \leqslant  \frac{\mu r }{n_2}\sigma_1$. 
Hence, we further bound $\circled{3}$ and $\circled{4}$ as
\begin{align*}
\circled{3} &\leqslant   \frac{\|\boldsymbol{\Omega}-p\boldsymbol{J}\|}{p} \left(\sum_{k=1}^{n_1}\|(\boldsymbol{\Delta}_{\boldsymbol{X}})_{k,\cdot}\|_2^4 + \frac{\mu r }{n_2}\sigma_1 \|\boldsymbol{\Delta}_{\boldsymbol{Y}}\|_F^2  \right),\\
\circled{4} &\leqslant   \frac{\|\boldsymbol{\Omega}-p\boldsymbol{J}\|}{p} \left(\sum_{k=1}^{n_2}\|(\boldsymbol{\Delta}_{\boldsymbol{Y}})_{k,\cdot}\|_2^4 + \frac{\mu r}{n_1}\sigma_1 \|\boldsymbol{\Delta}_{\boldsymbol{X}}\|_F^2  \right).
\end{align*}
We 
combine the foregoing inequalities to obtain
\begin{align*}
K_2 & \leqslant  \left| \frac{1}{p}\left\| \mathcal{P}_{\Omega}\left(\boldsymbol{\Delta}_{\boldsymbol{X}}\boldsymbol{\Delta}_{\boldsymbol{Y}}^\top \right) \right\|_F^2 - \|\boldsymbol{\Delta}_{\boldsymbol{X}}\boldsymbol{\Delta}_{\boldsymbol{Y}}^\top \|_F^2 \right| \\
&~~~~ + \left| \frac{3}{p}\left\| \mathcal{P}_{\Omega}\left( \boldsymbol{X}\boldsymbol{Y}^\top - \boldsymbol{U}\boldsymbol{V}^\top  \right) \right\|_F^2 - 3 \|\boldsymbol{X}\boldsymbol{Y}^\top - \boldsymbol{U}\boldsymbol{V}^\top\|_F^2 \right|\\ 
& \leqslant   \circled{2} + 3\big( \circled{1}  +\circled{2} +   \circled{3} +  \circled{4}\big)
\\
& \leqslant  0.0006 \left(\left\| \boldsymbol{\Delta}_{\boldsymbol{X}}\boldsymbol{V}^\top\right\|_F^2 + \left\|\boldsymbol{U}\boldsymbol{\Delta}_{\boldsymbol{Y}}^\top \right\|_F^2\right)
\\
&~~~~ + \frac{\|\boldsymbol{\Omega}-p\boldsymbol{J}\|}{p} \left(5\sum_{k=1}^{n_1}\|(\boldsymbol{\Delta}_{\boldsymbol{X}})_{k,\cdot}\|_2^4 + 5\sum_{k=1}^{n_2}\|(\boldsymbol{\Delta}_{\boldsymbol{Y}})_{k,\cdot}\|_2^4 + 3\frac{\mu r}{n_1}\sigma_1 \|\boldsymbol{\Delta}_{\boldsymbol{X}}\|_F^2 + 3\frac{\mu r}{n_2}\sigma_1 \|\boldsymbol{\Delta}_{\boldsymbol{Y}}\|_F^2\right).
\end{align*}

\subsubsection{Control of $K_3$}
For $K_3$, when $\alpha \geqslant  100 \sqrt{\frac{\mu r \sigma_1}{n_{\min}}}$, we have 
\begin{equation*}
K_3 \leqslant  200 \lambda\alpha^2 (\|\boldsymbol{\Delta}_{\boldsymbol{X}}\|_F^2 + \|\boldsymbol{\Delta}_{\boldsymbol{Y}}\|_F^2) - 0.3 \lambda \left(\sum_{k=1}^{n_1}\|(\boldsymbol{\Delta}_{\boldsymbol{X}})_{k,\cdot}\|_2^4 + \sum_{k=1}^{n_2}\|(\boldsymbol{\Delta}_{\boldsymbol{Y}})_{k,\cdot}\|_2^4\right).
\end{equation*}
The remaining argument is the same as that in \citet[Lemma 11]{ge2017no}, which has also been employed in \citet[Lemma 14]{chen2017memory}.  


\subsubsection{Putting $K_2$ and $K_3$ together}
\label{sec:put_k2_k3_together}
Combining the upper bounds on $K_2$ and $K_3$, we obtain that
\begin{align*}
K_2 + K_3 & \leqslant   0.0006 \left(\left\| \boldsymbol{\Delta}_{\boldsymbol{X}}\boldsymbol{V}^\top\right\|_F^2  + \left\|\boldsymbol{U}\boldsymbol{\Delta}_{\boldsymbol{Y}}^\top \right\|_F^2\right)
\\
&~~~~ + \left(\frac{ 3\mu r \sigma_1 \|\boldsymbol{\Omega}-p\boldsymbol{J}\|}{p n_{\min}} + 200 \lambda\alpha^2\right)
 \left(\|\boldsymbol{\Delta}_{\boldsymbol{X}}\|_F^2 +  \|\boldsymbol{\Delta}_{\boldsymbol{Y}}\|_F^2\right)
\\
&~~~~ + \left(\frac{5\|\boldsymbol{\Omega}-p\boldsymbol{J}\|}{p}  - 0.3\lambda\right) \left(\sum_{k=1}^{n_1}\|(\boldsymbol{\Delta}_{\boldsymbol{X}})_{k,\cdot}\|_2^4 + \sum_{k=1}^{n_2}\|(\boldsymbol{\Delta}_{\boldsymbol{Y}})_{k,\cdot}\|_2^4\right).
\end{align*} 
By the third inequality in \eqref{eq:rate_tuning_ass}, 
\[
  \left(\frac{ 3\mu r \sigma_1 \|\boldsymbol{\Omega}-p\boldsymbol{J}\|}{p n_{\min}} + 200 \lambda\alpha^2\right) \leqslant \frac{200\sigma_r}{C_0}.
\]
In addition, by the fourth inequality in \eqref{eq:rate_tuning_ass},
\[
  \left(\frac{5\|\boldsymbol{\Omega}-p\boldsymbol{J}\|}{p}  - 0.3\lambda\right) \leqslant \left(\frac{5}{C_0}-0.3\right)\lambda. 
\]
Therefore, with sufficiently large $C_0$ (e.g., $C_0 = 200/0.0004 = 5\times 10^5$), 
\[
\begin{split}
K_2 + K_3 & \leqslant  0.0006 \left(\left\| \boldsymbol{\Delta}_{\boldsymbol{X}}\boldsymbol{V}^\top\right\|_F^2  + \left\|\boldsymbol{U}\boldsymbol{\Delta}_{\boldsymbol{Y}}^\top \right\|_F^2\right) + 0.0004\sigma_r \left(\|\boldsymbol{\Delta}_{\boldsymbol{X}}\|_F^2 +  \|\boldsymbol{\Delta}_{\boldsymbol{Y}}\|_F^2\right) \\
&~~~~~+ (10^{-5}-0.3) \lambda \left(\sum_{k=1}^{n_1}\|(\boldsymbol{\Delta}_{\boldsymbol{X}})_{k,\cdot}\|_2^4 + \sum_{k=1}^{n_2}\|(\boldsymbol{\Delta}_{\boldsymbol{Y}})_{k,\cdot}\|_2^4\right)\\
& \leqslant   0.0006 \left(\left\| \boldsymbol{\Delta}_{\boldsymbol{X}}\boldsymbol{V}^\top\right\|_F^2  + \left\|\boldsymbol{U}\boldsymbol{\Delta}_{\boldsymbol{Y}}^\top \right\|_F^2\right) + 0.0004\sigma_r \left(\|\boldsymbol{\Delta}_{\boldsymbol{X}}\|_F^2 +  \|\boldsymbol{\Delta}_{\boldsymbol{Y}}\|_F^2\right),
\end{split}
\]
where the last inequality uses the facts that $\lambda\geqslant 0$ and that $\sum_{k=1}^{n_1}\|(\boldsymbol{\Delta}_{\boldsymbol{X}})_{k,\cdot}\|_2^4 + \sum_{k=1}^{n_2}\|(\boldsymbol{\Delta}_{\boldsymbol{Y}})_{k,\cdot}\|_2^4  \geqslant 0$.
By Proposition \ref{prop:basic_properties}, we have $\left\|\boldsymbol{U}\boldsymbol{\Delta}_{\boldsymbol{Y}}^\top \right\|_F^2 \geqslant  \sigma_r^2(\boldsymbol{U}) \left\|\boldsymbol{\Delta}_{\boldsymbol{Y}} \right\|_F^2 =  \sigma_r \left\|\boldsymbol{\Delta}_{\boldsymbol{Y}}  \right\|_F^2$ and $\left\|\boldsymbol{V}\boldsymbol{\Delta}_{\boldsymbol{X}}^\top \right\|_F^2 \geqslant \sigma_r^2(\boldsymbol{V}) \left\|\boldsymbol{\Delta}_{\boldsymbol{X}}  \right\|_F^2 =  \sigma_r \left\|\boldsymbol{\Delta}_{\boldsymbol{X}}\right\|_F^2$. 
Together with the last display, they lead to
\begin{align}
\label{eq:k234_simple}
K_2 + K_3 \leqslant  & 0.001 \left(\left\| \boldsymbol{\Delta}_{\boldsymbol{X}}\boldsymbol{V}^\top\right\|_F^2  + \left\|\boldsymbol{U}\boldsymbol{\Delta}_{\boldsymbol{Y}}^\top \right\|_F^2\right).
\end{align}
This completes the proof of Lemma \ref{lem:K2K3}.

\subsection{Proof of Lemma \ref{lem:K4}}
\label{sec:pf-lem-K4}
\begin{proof}

By the way we define $\boldsymbol{\Delta}_{\boldsymbol{Z}}$,
\begin{equation} 
\label{eq:Jin_equality}
	\begin{split}
& \left\| \boldsymbol{Z}\boldsymbol{Z}^\top -\boldsymbol{W}\boldsymbol{W}^\top \right\|_F^2\\
 = &\left\| \boldsymbol{\Delta}_{\boldsymbol{Z}} \boldsymbol{W}^\top + \boldsymbol{W}\boldsymbol{\Delta}_{\boldsymbol{Z}}^\top  + \boldsymbol{\Delta}_{\boldsymbol{Z}}\boldsymbol{\Delta}_{\boldsymbol{Z}}^\top\right\|_F^2
\\
=& \|\boldsymbol{\Delta}_{\boldsymbol{Z}}\boldsymbol{\Delta}_{\boldsymbol{Z}}^\top\|_F^2 + 2\|\boldsymbol{\Delta}_{\boldsymbol{Z}} \boldsymbol{W}^\top\|_F^2  + 2 \langle \boldsymbol{\Delta}_{\boldsymbol{Z}} \boldsymbol{W}^\top,\boldsymbol{W}\boldsymbol{\Delta}_{\boldsymbol{Z}}^\top \rangle  + 4 \langle  \boldsymbol{\Delta}_{\boldsymbol{Z}} \boldsymbol{W}^\top, \boldsymbol{\Delta}_{\boldsymbol{Z}}\boldsymbol{\Delta}_{\boldsymbol{Z}}^\top \rangle\\
 = & \|\boldsymbol{\Delta}_{\boldsymbol{Z}}^\top \boldsymbol{\Delta}_{\boldsymbol{Z}}\|_F^2 +2\langle \boldsymbol{\Delta}_{\boldsymbol{Z}}^\top \boldsymbol{\Delta}_{\boldsymbol{Z}}, \boldsymbol{W}^\top \boldsymbol{W} \rangle  +2\langle \boldsymbol{\Delta}_{\boldsymbol{Z}}^\top \boldsymbol{W}, \boldsymbol{W}^\top \boldsymbol{\Delta}_{\boldsymbol{Z}}\rangle +4\langle \boldsymbol{\Delta}_{\boldsymbol{Z}}^\top \boldsymbol{\Delta}_{\boldsymbol{Z}}, \boldsymbol{\Delta}_{\boldsymbol{Z}}^\top \boldsymbol{W} \rangle.
\end{split}
\end{equation}
Here we have used the facts that $\langle \boldsymbol{A},\boldsymbol{B}\rangle = \operatorname{trace}(\boldsymbol{A}^\top \boldsymbol{B})$ and $\langle \boldsymbol{A},\boldsymbol{B}\rangle = \langle \boldsymbol{B},\boldsymbol{A}\rangle$.
By the definition of $K_1$ in \eqref{eq:K_noisy_decomp}, \eqref{eq:Jin_equality} further implies  
\begin{equation}
\label{eq:K_factorization_noisy2}
  \begin{split}
K_1 & =  - \frac{1}{2}\|\boldsymbol{\Delta}_{\boldsymbol{Z}}^\top \boldsymbol{\Delta}_{\boldsymbol{Z}}\|_F^2 - \frac{3}{2} \langle \boldsymbol{\Delta}_{\boldsymbol{Z}}^\top \boldsymbol{\Delta}_{\boldsymbol{Z}}, \boldsymbol{W}^\top \boldsymbol{W}  \rangle - \frac{3}{2} \langle \boldsymbol{\Delta}_{\boldsymbol{Z}}^\top \boldsymbol{W}, \boldsymbol{W}^\top \boldsymbol{\Delta}_{\boldsymbol{Z}} \rangle \\
&~~~ - 3\langle \boldsymbol{\Delta}_{\boldsymbol{Z}}^\top \boldsymbol{\Delta}_{\boldsymbol{Z}}, \boldsymbol{\Delta}_{\boldsymbol{Z}}^\top \boldsymbol{W} \rangle.
  \end{split}
\end{equation} 
Moreover, condition \eqref{eq:rotation_rect} leads to
\begin{equation}
\label{eq:psd_z}
\boldsymbol{Z}^\top \boldsymbol{W} = \mtx{X}^\top \mtx{U} + \mtx{Y}^\top \mtx{V} \succeq \boldsymbol{0},
\end{equation}
which implies the symmetricity of $\boldsymbol{W}^\top \boldsymbol{\Delta}_{\boldsymbol{Z}} = \mtx{W}^\top \mtx{Z} - \mtx{Z}^\top \mtx{Z}$ (this is a crucial step for the analysis in \citet{jin2017escape} and \citet{ge2017no}). 
Hence,
\begin{equation}
\label{eq:symm_implication}
 \langle \boldsymbol{\Delta}_{\boldsymbol{Z}}^\top \boldsymbol{W}, \boldsymbol{W}^\top \boldsymbol{\Delta}_{\boldsymbol{Z}} \rangle = \|\boldsymbol{\Delta}_{\boldsymbol{Z}}^\top \boldsymbol{W}\|_F^2. 
\end{equation}      
Combining \eqref{eq:symm_implication} with \eqref{eq:K_factorization_noisy2} we have 
\begin{equation*}
	K_1 =  - 0.5 \|\boldsymbol{\Delta}_{\boldsymbol{Z}}^\top \boldsymbol{\Delta}_{\boldsymbol{Z}}\|_F^2 - 1.5\langle \boldsymbol{\Delta}_{\boldsymbol{Z}}^\top \boldsymbol{\Delta}_{\boldsymbol{Z}}, \boldsymbol{W}^\top \boldsymbol{W} \rangle - 1.5 \|\boldsymbol{\Delta}_{\boldsymbol{Z}}^\top \boldsymbol{W}\|_F^2 - 3\langle \boldsymbol{\Delta}_{\boldsymbol{Z}}^\top \boldsymbol{\Delta}_{\boldsymbol{Z}}, \boldsymbol{\Delta}_{\boldsymbol{Z}}^\top \boldsymbol{W} \rangle.
\end{equation*}
Therefore, based on \eqref{eq:k234_simple}, we are able to upper bound the righthand side of \eqref{eq:K_noisy_decomp_ks} as
\begin{equation*}
\begin{split}
& K_1+K_2+K_3+K_4 \\
& \leqslant    - 0.5\|\boldsymbol{\Delta}_{\boldsymbol{Z}}^\top \boldsymbol{\Delta}_{\boldsymbol{Z}}\|_F^2 - 1.5\langle \boldsymbol{\Delta}_{\boldsymbol{Z}}^\top \boldsymbol{\Delta}_{\boldsymbol{Z}}, \boldsymbol{W}^\top \boldsymbol{W} \rangle- 1.5 \|\boldsymbol{\Delta}_{\boldsymbol{Z}}^\top \boldsymbol{W}\|_F^2 - 3\langle \boldsymbol{\Delta}_{\boldsymbol{Z}}^\top \boldsymbol{\Delta}_{\boldsymbol{Z}}, \boldsymbol{\Delta}_{\boldsymbol{Z}}^\top \boldsymbol{W} \rangle\\
&~~~\,  +0.001 \left(\| \boldsymbol{\Delta}_{\boldsymbol{X}}\boldsymbol{V}^\top\|_F^2  + \|\boldsymbol{U}\boldsymbol{\Delta}_{\boldsymbol{Y}}^\top \|_F^2\right)+ |K_4|.
	\end{split}
\end{equation*}
%
Furthermore, due to cyclic property of trace,
\begin{equation}
\label{eq:zxy2}
\begin{split}
\langle \boldsymbol{\Delta}_{\boldsymbol{Z}}^\top \boldsymbol{\Delta}_{\boldsymbol{Z}}, \boldsymbol{W}^\top \boldsymbol{W}  \rangle & = \operatorname{trace}(\boldsymbol{\Delta}_{\boldsymbol{Z}}^\top \boldsymbol{\Delta}_{\boldsymbol{Z}} \boldsymbol{W}^\top \boldsymbol{W}) = \|\boldsymbol{W}\boldsymbol{\Delta}_{\boldsymbol{Z}}^\top\|_F^2\\
& \geqslant  \| \boldsymbol{\Delta}_{\boldsymbol{X}}\boldsymbol{V}^\top \|_F^2  + \|\boldsymbol{U}\boldsymbol{\Delta}_{\boldsymbol{Y}}^\top \|_F^2.
\end{split}
\end{equation}
Together with the second last display, \eqref{eq:zxy2}, we further have
\begin{equation}
\label{eq:k1234_simple}
  \begin{split}
& K_1+K_2+K_3+K_4\\
& \leqslant    - 0.5\|\boldsymbol{\Delta}_{\boldsymbol{Z}}^\top \boldsymbol{\Delta}_{\boldsymbol{Z}}\|_F^2 - 1.499\langle \boldsymbol{\Delta}_{\boldsymbol{Z}}^\top \boldsymbol{\Delta}_{\boldsymbol{Z}}, \boldsymbol{W}^\top \boldsymbol{W} \rangle - 1.5 \|\boldsymbol{\Delta}_{\boldsymbol{Z}}^\top \boldsymbol{W}\|_F^2\\
&~~~\, - 3\langle \boldsymbol{\Delta}_{\boldsymbol{Z}}^\top \boldsymbol{\Delta}_{\boldsymbol{Z}}, \boldsymbol{\Delta}_{\boldsymbol{Z}}^\top \boldsymbol{W} \rangle + |K_4|.
  \end{split}
\end{equation}     
Now \eqref{eq:psd_z} further implies that
\begin{equation*}
  \langle \boldsymbol{\Delta}_{\boldsymbol{Z}}^\top \boldsymbol{\Delta}_{\boldsymbol{Z}}, \boldsymbol{W}^\top \boldsymbol{W}  \rangle +  \langle \boldsymbol{\Delta}_{\boldsymbol{Z}}^\top \boldsymbol{\Delta}_{\boldsymbol{Z}}, \boldsymbol{\Delta}_{\boldsymbol{Z}}^\top \boldsymbol{W}  \rangle =  \langle \boldsymbol{\Delta}_{\boldsymbol{Z}}^\top \boldsymbol{\Delta}_{\boldsymbol{Z}}, \boldsymbol{Z}^\top \boldsymbol{W}  \rangle \geqslant  0,
\end{equation*}        
in which we use the fact that the inner product of two PSD matrices is nonnegative. 
Thus,
\begin{align}
& \hskip -0.75em K_1+K_2+K_3+K_4 \nonumber \\
\leqslant & -0.5\|\boldsymbol{\Delta}_{\boldsymbol{Z}}^\top \boldsymbol{\Delta}_{\boldsymbol{Z}}\|_F^2 - 1.499\left(\langle \boldsymbol{\Delta}_{\boldsymbol{Z}}^\top \boldsymbol{\Delta}_{\boldsymbol{Z}}, \boldsymbol{W}^\top \boldsymbol{W} \rangle +  \langle \boldsymbol{\Delta}_{\boldsymbol{Z}}^\top \boldsymbol{\Delta}_{\boldsymbol{Z}}, \boldsymbol{\Delta}_{\boldsymbol{Z}}^\top \boldsymbol{W}  \rangle\right)- 1.5\|\boldsymbol{\Delta}_{\boldsymbol{Z}}^\top \boldsymbol{W}\|_F^2 \nonumber\\
&- 1.501  \langle \boldsymbol{\Delta}_{\boldsymbol{Z}}^\top \boldsymbol{\Delta}_{\boldsymbol{Z}}, \boldsymbol{\Delta}_{\boldsymbol{Z}}^\top \boldsymbol{W}  \rangle + |K_4| \nonumber\\
\leqslant &  -0.5\|\boldsymbol{\Delta}_{\boldsymbol{Z}}^\top \boldsymbol{\Delta}_{\boldsymbol{Z}}\|_F^2 - 1.5\|\boldsymbol{\Delta}_{\boldsymbol{Z}}^\top \boldsymbol{W}\|_F^2 - 1.501  \langle \boldsymbol{\Delta}_{\boldsymbol{Z}}^\top \boldsymbol{\Delta}_{\boldsymbol{Z}}, \boldsymbol{\Delta}_{\boldsymbol{Z}}^\top \boldsymbol{W}  \rangle + |K_4| \nonumber\\
\leqslant & -0.5\|\boldsymbol{\Delta}_{\boldsymbol{Z}}^\top \boldsymbol{\Delta}_{\boldsymbol{Z}}\|_F^2 - 1.5\|\boldsymbol{\Delta}_{\boldsymbol{Z}}^\top \boldsymbol{W}\|_F^2 +1.501 \|\boldsymbol{\Delta}_{\boldsymbol{Z}}^\top \boldsymbol{\Delta}_{\boldsymbol{Z}}\|_F\|\boldsymbol{\Delta}_{\boldsymbol{Z}}^\top \boldsymbol{W}\|_F  + |K_4|.
\label{eq:k1234_final}
\end{align}

Note that \eqref{eq:k1234_final} implies
\[
\begin{split}
0.5\|\boldsymbol{\Delta}_{\boldsymbol{Z}}^\top \boldsymbol{\Delta}_{\boldsymbol{Z}}\|_F^2 + 1.5\|\boldsymbol{\Delta}_{\boldsymbol{Z}}^\top \boldsymbol{W}\|_F^2 - 1.501 \| \boldsymbol{\Delta}_{\boldsymbol{Z}}^\top \boldsymbol{\Delta}_{\boldsymbol{Z}}\|_F \|\boldsymbol{\Delta}_{\boldsymbol{Z}}^\top \boldsymbol{W} \|_F \leqslant  |K_4|.
\end{split}
\]
We apply Young's inequality to the LHS of the inequality and subsequently deduce the upper bounds
\begin{equation}
\label{eq:dzdz}
\|\boldsymbol{\Delta}_{\boldsymbol{Z}} \boldsymbol{\Delta}_{\boldsymbol{Z}}^\top\|_F = \|\boldsymbol{\Delta}_{\boldsymbol{Z}}^\top \boldsymbol{\Delta}_{\boldsymbol{Z}}\|_F\leqslant   3\sqrt{|K_4|}, \quad \|\boldsymbol{\Delta}_{\boldsymbol{Z}}^\top \boldsymbol{W}\|_F \leqslant   2\sqrt{|K_4|}. 
\end{equation}
By \eqref{eq:k1234_simple} as well as \eqref{eq:k1234geq0}, we have
\begin{equation*}
  \begin{split}
 & 1.499 \langle \boldsymbol{\Delta}_{\boldsymbol{Z}}^\top \boldsymbol{\Delta}_{\boldsymbol{Z}}, \boldsymbol{W}^\top \boldsymbol{W}  \rangle\\
   \leqslant &  - 0.5 \|\boldsymbol{\Delta}_{\boldsymbol{Z}}^\top \boldsymbol{\Delta}_{\boldsymbol{Z}}\|_F^2 - 1.5\|\boldsymbol{\Delta}_{\boldsymbol{Z}}^\top \boldsymbol{W}\|_F^2 - 3 \langle \boldsymbol{\Delta}_{\boldsymbol{Z}}^\top \boldsymbol{\Delta}_{\boldsymbol{Z}}, \boldsymbol{\Delta}_{\boldsymbol{Z}}^\top \boldsymbol{W} \rangle + |K_4| \\
\leqslant & 3 \|\boldsymbol{\Delta}_{\boldsymbol{Z}}^\top \boldsymbol{\Delta}_{\boldsymbol{Z}}\|_F \|\boldsymbol{\Delta}_{\boldsymbol{Z}}^\top \boldsymbol{W}\|_F+ |K_4|\\
\leqslant& 19|K_4|.
  \end{split}
\end{equation*}   
Combining with \eqref{eq:zxy2} we have claim 2:
\begin{equation}\label{eq:dxv_udy}
\| \boldsymbol{\Delta}_{\boldsymbol{X}}\boldsymbol{V}^\top \|_F^2  + \|\boldsymbol{U}\boldsymbol{\Delta}_{\boldsymbol{Y}}^\top \|_F^2 \leqslant 13|K_4|.
\end{equation} 
This concludes the proof of Lemma \ref{lem:K4}.
  
\end{proof}

\subsection{Proof of Lemma \ref{lemma:eq_of_rate}}
\label{sec:proof_lemma_eq_of_rate}
\begin{proof}
  First of all, for the first line and second line of \eqref{eq:rate_tuning_ass}, by first line and third line of assumption \eqref{eq:main_rate_tuning_ass}, we have 
  \[
    p\geqslant C_1\frac{\mu r \log n_{\max}}{n_{\min}}
  \]
  and 
  \[
    \alpha \geqslant C_2 \sqrt{\frac{\mu r \sigma_1}{n_{\min}}}.
    \]
  
  Furthermore, for the fourth line of \eqref{eq:rate_tuning_ass}, by second line of \eqref{eq:main_rate_tuning_ass}, $p\lambda \geqslant C_2\sqrt{n_{\max}p}$. From Lemma \ref{lem:vu14_psd}, if $p\geqslant C_v\frac{\log n_{\max}}{n_{\min}}$, in an event $E_v$ with probability $\mathbb{P}[E_v]\geqslant 1-(n_1+n_2)^{-5}$, $\|\boldsymbol{\Omega}-p\boldsymbol{J}\|\leqslant C_v\sqrt{n_{\max}p}$, where $C_v$ is defined in Lemma \ref{lem:vu14_psd}. Therefore, in the event $E_v$, $p\lambda\geqslant \frac{C_2}{C_v}\|\boldsymbol{\Omega}-p\boldsymbol{J}\|$.  
  
  Finally, by the fact that $\lambda\leqslant 10C_2\sqrt{\frac{n_{\max}}{p}}$, $\alpha \leqslant 10C_2\sqrt{\frac{\mu r \sigma_1}{n_{\min}}}$ and $p\lambda\geqslant \frac{C_2}{C_v}\|\boldsymbol{\Omega}-p\boldsymbol{J}\|$, we have 
  \[
    \frac{\mu r \sigma_1 \|\boldsymbol{\Omega}-p\boldsymbol{J}\|}{p n_{\min}} \leqslant \frac{C_v \mu r \sigma_1  p\lambda}{C_2p n_{\min}} \leqslant \frac{C_v}{C_2} \frac{\mu r \sigma_1}{n_{\min}} 10 C_2\sqrt{\frac{n_{\max}}{p}} = 10 C_v \sigma_r \sqrt{\frac{\mu^2 r^2 \kappa^2 n_{\max}}{n_{\min}^2p}}
  \]
  and 
  \[
      \lambda \alpha^2 \leqslant 10^3  C_2^3 \sigma_r \sqrt{\frac{\mu^2 r^2 \kappa^2 n_{\max}}{n_{\min}^2p}}.
  \]
  By the first line of \eqref{eq:main_rate_tuning_ass},  
  \[
    p\geqslant  \frac{C_1}{n_{\min}^2} n_{\max}\mu^2 r^2\kappa^2.
\]
Therefore,
\[
  \frac{\mu r \sigma_1 \|\boldsymbol{\Omega}-p\boldsymbol{J}\|}{p n_{\min}} + \lambda \alpha^2 \leqslant \frac{10C_v + 10^3  C_2^3}{\sqrt{C_1}} \sigma_r.
\]
In other words, 
\[
  \sigma_r \geqslant \frac{\sqrt{C_1}}{10C_v + 10^3  C_2^3} \left(\frac{\mu r \sigma_1 \|\boldsymbol{\Omega}-p\boldsymbol{J}\|}{p n_{\min}} + \lambda \alpha^2\right)   .
\]
Therefore, by choosing
\[
    C_1 = \max\{ C_0,C_c,C_0^2 (10C_v+10^3C_2^3)^2 \}
\]
and
\[
    C_2 = \max\{100,C_0C_v\} 
\]
finishes the proof of the first part of the Lemma.
 
Recall by the way we define $C_1$, if \eqref{eq:main_rate_tuning_ass} is satisfied, by Lemma \ref{lem:rip}, in an event $E_c$ with probability $\mathbb{P}[E_c] \geqslant 1-(n_1+n_2)^{-5}$, \eqref{eq:rip} holds. Therefore, let $E = E_v\cap E_c$, then by union bound, $\mathbb{P}[E]\geqslant 1-(n_1+n_2)^{-3}$.

\end{proof}

\subsection{Proof of Proposition \ref{prop:basic_properties}}\label{sec:proof_of_prop_basic_properties}
\begin{proof}
  First, since $\boldsymbol{M}^{\star}$ has SVD $\boldsymbol{M}^{\star} =  {\boldsymbol{U}^{\star} }\boldsymbol{\Lambda} {\boldsymbol{V}^{\star} }^\top$, we have 
  \[
      \colspan( \boldsymbol{U}^{\star} ) =  \colspan( \boldsymbol{M}^{\star}) \quad \textrm{and} \quad \colspan( \boldsymbol{V}^{\star} ) =  \rowspan( \boldsymbol{M}^{\star})
  \]
  as well as
  \[
      \dim(\colspan( \boldsymbol{M}^{\star})) = \dim(\rowspan( \boldsymbol{M}^{\star})) = r.
  \]
  From \eqref{eq:balanced parameterized factorization}, we also have 
  \[
    \colspan( \boldsymbol{M}^{\star}) \subset \colspan( \boldsymbol{U})\quad \textrm{and} \quad \rowspan( \boldsymbol{M}^{\star}) \subset \colspan( \boldsymbol{V}).
  \]
  By the way we define $\boldsymbol{U}$ and $\boldsymbol{V}$, we have $\dim(\colspan( \boldsymbol{U}))\leqslant  r$ and $\dim(\colspan( \boldsymbol{V}))\leqslant  r$. Therefore, $\colspan(\boldsymbol{U}) = \colspan( \boldsymbol{U}^{\star} )$ and $\colspan(\boldsymbol{V}) = \colspan( \boldsymbol{V}^{\star} )$. 

  From second equation in \eqref{eq:balanced parameterized factorization}, $\boldsymbol{U}^\top\boldsymbol{U} = \boldsymbol{V}^\top\boldsymbol{V}$, therefore,
  \[
      \sigma_i(\boldsymbol{U}) = \sqrt{\lambda_i(\boldsymbol{U}^\top\boldsymbol{U})} = \sqrt{\lambda_i(\boldsymbol{V}^\top\boldsymbol{V})}= \sigma_i(\boldsymbol{V}), \; i = 1,2,\dots,r.
  \]
  Moreover, suppose $\boldsymbol{U}^\top \boldsymbol{U} = \boldsymbol{V}^\top\boldsymbol{V} = \boldsymbol{B}\boldsymbol{D}^2\boldsymbol{B}^\top$ be a fixed eigenvalue decomposition of $\boldsymbol{U}^\top \boldsymbol{U}$, with $\boldsymbol{B}\in \mathsf{O}(r)$ and $\boldsymbol{D}\in\mathbb{R}^{r\times r}$ diagonal matrix. Then the reduced SVD of $\boldsymbol{U}$ and $\boldsymbol{V}$ can be written as 
  \[
    \boldsymbol{U} = \boldsymbol{A}_{\boldsymbol{U}}\boldsymbol{D}\boldsymbol{B}^\top ,\quad \boldsymbol{V} = \boldsymbol{A}_{\boldsymbol{V}}\boldsymbol{D}\boldsymbol{B}^\top 
  \]
  with $\boldsymbol{A}_{\boldsymbol{U}}\in\mathbb{R}^{n_1\times r},\boldsymbol{A}_{\boldsymbol{V}}\in \mathbb{R}^{n_2\times r}$ satisfying $\boldsymbol{A}_{\boldsymbol{U}}^\top \boldsymbol{A}_{\boldsymbol{U}} = \boldsymbol{I}$ and $\boldsymbol{A}_{\boldsymbol{V}}^\top \boldsymbol{A}_{\boldsymbol{V}} = \boldsymbol{I}$. Therefore, $\boldsymbol{M}^{\star} = \boldsymbol{U}\boldsymbol{V}^\top = \boldsymbol{A}_{\boldsymbol{U}}\boldsymbol{D}^2 \boldsymbol{A}_{\boldsymbol{V}}^\top$. It is a reduced SVD of $\boldsymbol{M}^{\star}$ by the way we define $\boldsymbol{A}_{\boldsymbol{U}}$, $\boldsymbol{A}_{\boldsymbol{V}}$ and $\boldsymbol{D}$. Therefore, $\sigma_1(\boldsymbol{U}) = \sigma_1(\boldsymbol{V}) = \sqrt{\sigma_1}$ and $\sigma_r(\boldsymbol{U}) = \sigma_r(\boldsymbol{V}) = \sqrt{\sigma_r}$.

  Moreover, there is $\boldsymbol{R}_{\boldsymbol{U}},\boldsymbol{R}_{\boldsymbol{V}}\in\mathsf{O}(r)$ such that $\boldsymbol{A}_{\boldsymbol{U}} = {\boldsymbol{U}^{\star} }\boldsymbol{R}_{\boldsymbol{U}}, \boldsymbol{A}_{\boldsymbol{V}} = {\boldsymbol{V}^{\star} }\boldsymbol{R}_{\boldsymbol{V}}$. Therefore,
  \begin{align*}
      \| \boldsymbol{U} \|_{2,\infty}^2 &= \|\boldsymbol{A}_{\boldsymbol{U}}\boldsymbol{D}\boldsymbol{B}^\top\|_{2,\infty}^2 = \|\boldsymbol{A}_{\boldsymbol{U}}\boldsymbol{D} \|_{2,\infty}^2 
      \\
      &\leqslant  \|\boldsymbol{A}_{\boldsymbol{U}} \|_{2,\infty}^2\| \boldsymbol{D} \|_{\ell_\infty}^2 = \sigma_1 \| {\boldsymbol{U}^{\star} }\boldsymbol{R}_{\boldsymbol{U}}\|_{2,\infty}^2 = \sigma_1 \| {\boldsymbol{U}^{\star} } \|_{2,\infty}^2 \leqslant  \frac{\mu r}{n_1} \sigma_1.
  \end{align*} 
  Similarly, we also have $\| \boldsymbol{V} \|_{2,\infty}^2  \leqslant  \frac{\mu r}{n_2} \sigma_1$.

\end{proof}

\section{Proof of Lemma \ref{lem:cons_noise}}
\label{sec:proof_lemma_subspace_noise}
\begin{proof}

Recall $\widetilde{\boldsymbol{U}}$ and $\widetilde{\boldsymbol{V}}$ are orthonormal basis matrices, $\boldsymbol{P}_{\widetilde{\boldsymbol{U}}} = \widetilde{\boldsymbol{U}}\widetilde{\boldsymbol{U}}^\top$, $\boldsymbol{P}_{\widetilde{\boldsymbol{V}}} =  \widetilde{\boldsymbol{V}} \widetilde{\boldsymbol{V}}^\top$. Therefore,
  \[
       \begin{split}
         & \|\boldsymbol{P}_{\widetilde{\boldsymbol{U}}}\mathcal{P}_{\Omega}(\boldsymbol{N})\boldsymbol{P}_{\widetilde{\boldsymbol{V}}}\| 
         = \| \widetilde{\boldsymbol{U}} \widetilde{\boldsymbol{U}}^\top \mathcal{P}_{\Omega}(\boldsymbol{N}) \widetilde{\boldsymbol{V}} \widetilde{\boldsymbol{V}}^\top \|
         = \|  \widetilde{\boldsymbol{U}}^\top \mathcal{P}_{\Omega}(\boldsymbol{N}) \widetilde{\boldsymbol{V}} \|.
       \end{split}
  \]
  The last equality uses the fact that $\widetilde{\boldsymbol{U}}$ and $\widetilde{\boldsymbol{V}}$ are orthonormal basis matrices, therefore 
  \[
  \|\widetilde{\boldsymbol{U}} \boldsymbol{A}\| = \|\boldsymbol{A}\|, \quad \|\boldsymbol{B} \widetilde{\boldsymbol{V}}^\top\| = \|\boldsymbol{B}\| 
\]  
for any $\boldsymbol{A}$, $\boldsymbol{B}$ with suitable size.

Due to the fact that $\Omega$ follows from Model \ref{mod:sampling1}, entries of $\mathcal{P}_{\Omega}(\boldsymbol{N})$ can be written as $[\mathcal{P}_{\Omega}(\boldsymbol{N})]_{i,j} = \delta_{i,j} N_{i,j}$, where $\delta_{i,j}$'s are i.i.d.\ Bernoulli random variables such that 
  \[
       \delta_{i,j} = \left\{ \begin{array}{ll}
         1 & \textrm{with probability }p\\
         0 & \textrm{with probability }1-p.
       \end{array} \right.
  \]
  And $N_{i,j}$'s are i.i.d.\ centered sub-exponential random variables. Moreover, $\delta_{i,j}$'s and $N_{i,j}$'s are mutually independent. Therefore, 
  \[
       \begin{split}
        \|\boldsymbol{P}_{\widetilde{\boldsymbol{U}}}\mathcal{P}_{\Omega}(\boldsymbol{N})\boldsymbol{P}_{\widetilde{\boldsymbol{V}}}\|  =& \|  \widetilde{\boldsymbol{U}}^\top \mathcal{P}_{\Omega}(\boldsymbol{N}) \widetilde{\boldsymbol{V}} \|  = \left\|    \widetilde{\boldsymbol{U}}^\top  \left(\sum_{i,j} \delta_{i,j}N_{i,j} \boldsymbol{e}_{i}\boldsymbol{e}_j^\top\right)  \widetilde{\boldsymbol{V}}  \right\| = \left\| \sum_{i,j} \delta_{i,j}N_{i,j} \widetilde{\boldsymbol{U}}_{i,\cdot} \widetilde{\boldsymbol{V}}_{j,\cdot}^\top  \right\|.
       \end{split}
  \]

Now let 
\[
	\boldsymbol{Q}_{i,j} \coloneqq \delta_{i,j}N_{i,j}  \left[ \begin{array}{cc} \boldsymbol{0} & \widetilde{\boldsymbol{U}}_{i,\cdot} \widetilde{\boldsymbol{V}}_{j,\cdot}^\top   \\ \widetilde{\boldsymbol{V}}_{j,\cdot}\widetilde{\boldsymbol{U}}_{i,\cdot}^\top & \boldsymbol{0}  \end{array} \right].
\]
Therefore, 
\[
	\|\boldsymbol{P}_{\widetilde{\boldsymbol{U}}}\mathcal{P}_{\Omega}(\boldsymbol{N})\boldsymbol{P}_{\widetilde{\boldsymbol{V}}}\|  = \|  \widetilde{\boldsymbol{U}}^\top \mathcal{P}_{\Omega}(\boldsymbol{N}) \widetilde{\boldsymbol{V}} \| = \left\|\sum_{i,j}\boldsymbol{Q}_{i,j}\right\|
\]
and $\mathbb{E}[\boldsymbol{Q}_{i,j}] = \boldsymbol{0}$. By following the symmetrization argument in \citet[Example 6.14]{wainwright2019high}, without loss of generality, we can assume that $N_{i,j}$'s are symmetric random variable, i.e., $N_{i,j} \overset{d}{=} -N_{i,j}$. Now we want to verify the Bernstein's condition \citep[Definition 6.10]{wainwright2019high} for $\boldsymbol{Q}_{i,j}$'s. For $k\geqslant 3$, 
\[
	\mathbb{E}\left[ \boldsymbol{Q}_{i,j}^k \right] =  \mathbb{E} \left[\delta_{i,j}^k N_{i,j}^k  \left[ \begin{array}{cc} \boldsymbol{0} & \widetilde{\boldsymbol{U}}_{i,\cdot} \widetilde{\boldsymbol{V}}_{j,\cdot}^\top   \\ \widetilde{\boldsymbol{V}}_{j,\cdot}\widetilde{\boldsymbol{U}}_{i,\cdot}^\top & \boldsymbol{0}  \end{array} \right]^k \right]=  p \mathbb{E} [N_{i,j}^k ] \left[ \begin{array}{cc} \boldsymbol{0} & \widetilde{\boldsymbol{U}}_{i,\cdot} \widetilde{\boldsymbol{V}}_{j,\cdot}^\top   \\ \widetilde{\boldsymbol{V}}_{j,\cdot}\widetilde{\boldsymbol{U}}_{i,\cdot}^\top & \boldsymbol{0}  \end{array} \right]^k.
\]
Due to the symmetry of $N_{i,j}$, $\mathbb{E}[N_{i,j}^k] = 0$ when $k\geqslant 3$ is odd, therefore, $\mathbb{E} [ \boldsymbol{Q}_{i,j}^k  ]  = 0$. For $k\geqslant 2$ even, we have
\[
\begin{split}
\left[ \begin{array}{cc} \boldsymbol{0} & \widetilde{\boldsymbol{U}}_{i,\cdot} \widetilde{\boldsymbol{V}}_{j,\cdot}^\top   \\ \widetilde{\boldsymbol{V}}_{j,\cdot}\widetilde{\boldsymbol{U}}_{i,\cdot}^\top & \boldsymbol{0}  \end{array} \right]^k = & \left[ \begin{array}{cc} (\widetilde{\boldsymbol{U}}_{i,\cdot} \widetilde{\boldsymbol{V}}_{j,\cdot}^\top\widetilde{\boldsymbol{V}}_{j,\cdot}\widetilde{\boldsymbol{U}}_{i,\cdot}^\top)^{k/2} & \boldsymbol{0}    \\ \boldsymbol{0}  & (\widetilde{\boldsymbol{V}}_{j,\cdot}\widetilde{\boldsymbol{U}}_{i,\cdot}^\top\widetilde{\boldsymbol{U}}_{i,\cdot} \widetilde{\boldsymbol{V}}_{j,\cdot}^\top)^{k/2} \end{array} \right]\\
 =& \|\widetilde{\boldsymbol{U}}_{i,\cdot}\|_2^k \|\widetilde{\boldsymbol{V}}_{j,\cdot}\|_2^k  \left[ \begin{array}{cc} \frac{1}{\|\widetilde{\boldsymbol{U}}_{i,\cdot}\|_2^2}\widetilde{\boldsymbol{U}}_{i,\cdot}  \widetilde{\boldsymbol{U}}_{i,\cdot}^\top  & \boldsymbol{0}    \\ \boldsymbol{0}  & \frac{1}{\|\widetilde{\boldsymbol{V}}_{j,\cdot}\|_2^2}\widetilde{\boldsymbol{V}}_{j,\cdot}  \widetilde{\boldsymbol{V}}_{j,\cdot}^\top  \end{array} \right],
\end{split} 
\]
which is a positive semidefinite matrix. And due to the fact that $N_{i,j}$'s satisfy the Bernstein condition, for $k\geqslant 2$,
\[
	\mathbb{E}[N_{i,j}^k] \leqslant \frac{1}{2}k!\nu^2b^{k-2}.
\]
Therefore, for $k\geqslant 3$ even, 
\[
	\mathbb{E}\left[ \boldsymbol{Q}_{i,j}^k \right] \preceq \frac{1}{2}k!\nu^2b^{k-2} p \|\widetilde{\boldsymbol{U}}_{i,\cdot}\|_2^k \|\widetilde{\boldsymbol{V}}_{j,\cdot}\|_2^k  \left[ \begin{array}{cc} \frac{1}{\|\widetilde{\boldsymbol{U}}_{i,\cdot}\|_2^2}\widetilde{\boldsymbol{U}}_{i,\cdot}  \widetilde{\boldsymbol{U}}_{i,\cdot}^\top  & \boldsymbol{0}    \\ \boldsymbol{0}  & \frac{1}{\|\widetilde{\boldsymbol{V}}_{j,\cdot}\|_2^2}\widetilde{\boldsymbol{V}}_{j,\cdot}  \widetilde{\boldsymbol{V}}_{j,\cdot}^\top  \end{array} \right].
\]
And we also have 
\[
\begin{split}
	\mathbb{V}\left[ \boldsymbol{Q}_{i,j} \right] =  & \mathbb{E}\left[ \boldsymbol{Q}_{i,j}^2 \right]  = p \mathbb{E} [N_{i,j}^2 ] \left[ \begin{array}{cc} \boldsymbol{0} & \widetilde{\boldsymbol{U}}_{i,\cdot} \widetilde{\boldsymbol{V}}_{j,\cdot}^\top   \\ \widetilde{\boldsymbol{V}}_{j,\cdot}\widetilde{\boldsymbol{U}}_{i,\cdot}^\top & \boldsymbol{0}  \end{array} \right]^2\\
	 =& p\nu^2 \|\widetilde{\boldsymbol{U}}_{i,\cdot}\|_2^2 \|\widetilde{\boldsymbol{V}}_{j,\cdot}\|_2^2  \left[ \begin{array}{cc} \frac{1}{\|\widetilde{\boldsymbol{U}}_{i,\cdot}\|_2^2}\widetilde{\boldsymbol{U}}_{i,\cdot}  \widetilde{\boldsymbol{U}}_{i,\cdot}^\top  & \boldsymbol{0}    \\ \boldsymbol{0}  & \frac{1}{\|\widetilde{\boldsymbol{V}}_{j,\cdot}\|_2^2}\widetilde{\boldsymbol{V}}_{j,\cdot}  \widetilde{\boldsymbol{V}}_{j,\cdot}^\top  \end{array} \right].
\end{split}
\]
Therefore, for $k\geqslant 3$,
\[
	\mathbb{E}\left[ \boldsymbol{Q}_{i,j}^k \right]  \preceq \frac{1}{2}k! b^{k-2} \|\widetilde{\boldsymbol{U}}_{i,\cdot}\|_2^{k-2} \|\widetilde{\boldsymbol{V}}_{j,\cdot}\|_2^{k-2}\mathbb{V}\left[ \boldsymbol{Q}_{i,j} \right].
\]
Therefore, $\boldsymbol{Q}_{i,j}$ satisfies Bernstein condition with parameter $b\|\widetilde{\boldsymbol{U}}_{i,\cdot}\|_2  \|\widetilde{\boldsymbol{V}}_{j,\cdot}\|_2 \leqslant b \sqrt{\frac{\mu_{\widetilde{U}} \mu_{\widetilde{V}}s_1s_2 }{n_1n_2}}$. Furthermore, 
\[
\begin{split}
	\frac{1}{n_1n_2} \left\|\sum_{(i,j)\in[n_1]\times [n_2]} \mathbb{V}\left[ \boldsymbol{Q}_{i,j} \right]\right\| =& \frac{1}{n_1n_2} p\nu^2  \left\|\sum_{(i,j)\in[n_1]\times [n_2]} \left[ \begin{array}{cc}  \|\widetilde{\boldsymbol{V}}_{j,\cdot}\|_2^2 \widetilde{\boldsymbol{U}}_{i,\cdot}  \widetilde{\boldsymbol{U}}_{i,\cdot}^\top  & \boldsymbol{0}    \\ \boldsymbol{0}  &\|\widetilde{\boldsymbol{U}}_{i,\cdot}\|_2^2 \widetilde{\boldsymbol{V}}_{j,\cdot}  \widetilde{\boldsymbol{V}}_{j,\cdot}^\top  \end{array} \right]  \right\| \\
	=& \frac{1}{n_1n_2} p\nu^2  \left\| \left[ \begin{array}{cc}  \|\widetilde{\boldsymbol{V}}\|_F^2 \widetilde{\boldsymbol{U}}^\top \widetilde{\boldsymbol{U}}  & \boldsymbol{0}    \\ \boldsymbol{0}  &\|\widetilde{\boldsymbol{U}}\|_F^2 \widetilde{\boldsymbol{V}}^\top \widetilde{\boldsymbol{V}}   \end{array} \right]  \right\|\\
	\leqslant & \frac{1}{n_1n_2} p\nu^2 (s_1+s_2).
\end{split}
\]
Where the last equality uses the fact that $\widetilde{\boldsymbol{U}}^\top \widetilde{\boldsymbol{U}} = \boldsymbol{I}$, $\widetilde{\boldsymbol{V}}^\top \widetilde{\boldsymbol{V}} = \boldsymbol{I}$. Then by \citet[Theorem 6.17]{wainwright2019high}, for all $t>0$,
\[
	\mathbb{P}\left[  \frac{1}{n_1n_2}\left\|\sum_{i,j}\boldsymbol{Q}_{i,j}\right\| \geqslant t \right]\leqslant 2 (n_1+n_2) \exp\left( -\frac{n_1n_2 t^2}{2\left(\frac{1}{n_1n_2} p\nu^2 (s_1+s_2)+ b \sqrt{\frac{\mu_{\widetilde{U}} \mu_{\widetilde{V}}s_1s_2 }{n_1n_2}} t\right)} \right).
\]
Therefore, by choosing $t$ as 
\[
	t = C_w \frac{1}{n_1n_2}\left( \sqrt{p\nu^2 (s_1+s_2)\log (n_1+n_2)}+ b\sqrt{\frac{\mu_{\widetilde{U}} \mu_{\widetilde{V}}s_1s_2 }{n_1n_2}} \log (n_1+n_2)\right)
\]
with absolute constant $C_w$ sufficiently large, say $C_w = 10$, then
\[
\begin{split}
	&\mathbb{P} \left[ \|\boldsymbol{P}_{\widetilde{\boldsymbol{U}}}\mathcal{P}_{\Omega}(\boldsymbol{N})\boldsymbol{P}_{\widetilde{\boldsymbol{V}}}\| \geqslant C_w \left( \sqrt{p\nu^2 (s_1+s_2)\log (n_1+n_2)}+ b\sqrt{\frac{\mu_{\widetilde{U}} \mu_{\widetilde{V}}s_1s_2 }{n_1n_2}} \log (n_1+n_2)\right)
	  \right] \\
	  =&   \mathbb{P} \left[ \left\|\sum_{i,j}\boldsymbol{Q}_{i,j}\right\| \geqslant C_w \left( \sqrt{p\nu^2 (s_1+s_2)\log (n_1+n_2)}+ b\sqrt{\frac{\mu_{\widetilde{U}} \mu_{\widetilde{V}}s_1s_2 }{n_1n_2}} \log (n_1+n_2)\right)
	  \right]\\
	  \leqslant & (n_1+n_2)^{-3}.
	  \end{split}
\]
\end{proof}

\end{appendix}

\end{document}